\documentclass[12pt]{article}

\usepackage{amssymb,amsmath,latexsym}
\voffset
-1.5cm

\usepackage{a4,amscd,graphicx,epsfig}
\usepackage[all]{xy}
\everymath{\displaystyle}

\newtheorem{teo}{Theorem}[section]
\newtheorem{lem}[teo]{Lemma}
\newtheorem{cor}[teo]{Corollary}
\newtheorem{exa}[teo]{Example}
\newtheorem{prop}[teo]{Proposition}
\newtheorem{defi}[teo]{Definition}

\newtheorem{remark}[teo]{Remark}
\newtheorem{remarks}[teo]{Remarks}

\newcommand{\ms}{\mathbb{S}}
\newcommand{\mr}{\mathbb{R}}
\newcommand{\mc}{\mathbb{C}}

\newcommand{\mh}{\mathbb{H}}

\newcommand{\mn}{\mathbb{N}}

\newcommand{\mx}{\mathbb{X}}

\newcommand{\md}{\mathbb{D}}

\newcommand{\Aa}{{\mathcal A}}

\newcommand{\Cc}{{\mathcal C}}
\newcommand{\Dd}{{\mathcal D}}
\newcommand{\Ee}{{\mathcal E}}
\newcommand{\Ff}{{\mathcal F}}
\newcommand{\Gg}{{\mathcal G}}

\newcommand{\Ll}{{\mathcal L}}
\newcommand{\Mm}{{\mathcal M}}
\newcommand{\Nn}{{\mathcal N}}
\newcommand{\Oo}{{\mathcal O}}
\newcommand{\Pp}{{\mathcal P}}
\newcommand{\Qq}{{\mathcal Q}}

\newcommand{\Tt}{{\mathcal T}}
\newcommand{\Uu}{{\mathcal U}}

\newcommand{\Vv}{{\mathcal V}}

\newcommand{\PP}{{\bf P}}

\newcommand{\sG}{\mathfrak s}
\newcommand{\lG}{\mathfrak l}

\newcommand{\hL}{{\widehat{L}}}

\newcommand{\C}{{\mathbb C}}

\newcommand{\Z}{{\mathbb Z}}
\newcommand{\R}{{\mathbb R}}
\newcommand{\X}{{\mathbb X}}

\def\ort#1{#1^\perp}

\def\arctgh{\mathrm{arctgh\, }}

\def\Dim{\emph{Proof : }}
\def\cvd{\nopagebreak\par\rightline{$_\blacksquare$}}

\def\E#1#2{\left\langle #1,#2 \right\rangle}  
\def\fut{\mathrm{I}^+}                         

\def\SL#1#2{\mathrm{SL}(#1,\,\mathbb{#2})}
\def\PSL#1#2{\mathbb P\mathrm{SL}(#1,\,\mathbb{#2})}
\def\OO{\mathrm{O}}
\def\SOO{\mathrm{SO}}

\def\ISO{\mathrm{Iso}}

\def\coom{\mathrm{H}^}
\def\eps{\varepsilon}                        

\def\Ad{\mathrm{Ad}}
\def\tr{\mathrm{tr}}                          

\def\d{\mathrm{d}}
\def\ch{\mathrm{ch\,}}
\def\sh{\mathrm{sh\,}}
\def\tgh{\mathrm{tgh\,}}

\def\arctgh{\mathrm{arctgh\, }}

\title{Wick rotations in 3D gravity:\\
$\Mm\Ll(\mh^2)$-spacetimes}

\author {Riccardo Benedetti$^1$ and Francesco Bonsante$^{1,2}$}

\date {}

\begin{document}

\maketitle

\vspace{0.7cm}
\noindent $(1)$ Dipartimento di Matematica, Universit\`a di Pisa,\\
Largo B. Pontecorvo, 5 , I-56127 Pisa.

\noindent $(2)$ Scuola Normale Superiore di Pisa
\smallskip

\noindent Email: benedett@dm.unipi.it ; f.bonsante@sns.it 
\begin{abstract}
\noindent 
``Ends of hyperbolic $3$-manifolds should support {\it canonical
Wick rotations}, so they realize effective interactions of their
ending globally hyperbolic spacetimes of constant curvature.''

We develop a consistent sector of {\it WR-rescaling theory} in 3D
gravity, that, in particular, concretizes the above guess for many
geometrically finite manifolds.  $\Mm\Ll(\mh^2)$-spacetimes are
solutions of pure Lorentzian 3D gravity encoded by {\it measured
geodesic laminations} of $\mh^2$, possibly invariant by any given
discrete, torsion-free subgroup $\Gamma$ of $PSL(2,\R)$.  The
rescalings which correlate spacetimes of different curvature, as well
as the {\it conformal} Wick rotations towards hyperbolic structures,
are {\it directed} by the gradient of the respective {\it canonical
cosmological times}, and have {\it universal} rescaling functions that
only depend on their value.  We get an insight into the WR-rescaling
mechanism by studying rays of $\Mm\Ll(\mh^2)$-spacetimes emanating
from the {\it static} case. In particular, we determine the
``derivatives'' at the starting point of each ray. We point out the
tamest behaviour of the cocompact $\Gamma$ case against the different
general one, even when $\Gamma$ is of cofinite area, but
non-cocompact. We analyze {\it broken $T$-symmetry} of AdS
$\Mm\Ll(\mh^2)$-spacetimes and related {\it earthquake failure}. This
helps us to figure out the main lines of development in order the
achieve a complete WR-rescaling theory.
\end{abstract}

\noindent \emph{Keywords: hyperbolic $3$-manifold, end, domain of
dependence, Wick rotation, cosmological time, measured geodesic
lamination, bending, earthquake.}

\vspace{0.8 cm}
\section{Introduction}\label{INTRODUCTION}

A basic fact of $3$-dimensional geometry is that the Ricci tensor
determines the Riemann tensor. This implies that the solutions of pure
3D gravity are the Lorentzian or Riemannian $3$-manifolds of {\it
constant curvature}. The sign of the curvature coincides with the sign
of the cosmological constant. We stipulate that all manifolds are {\it
oriented} and that the Lorentzian spacetimes are also {\it
time-oriented}. We could also include in the picture the presence of
{\it world lines} of particles, carrying some concentrated
singularities of the metric. A typical example is given by the {\it
cone manifolds} of constant curvature with cone locus at some embedded
link. The cone angles reflect the ``mass'' of the particles. In the
Lorentzian case we also require that the world lines are of causal
type (see e.g. \cite {BG}(2)). However, in the present paper we will
confine ourselves to the pure gravity case.
\smallskip

Sometimes gravity is studied by considering separately its different
``sectors'', according to the metric signature (Lorentzian or
Euclidean), and the sign of the cosmological constant. By using the
comprehensive term ``3D gravity'', we would suggest to
consider it as a unitary body, where different sectors actually
interact. In that direction, the following is one of our leading
ideas:
\smallskip

\noindent 
``Let $Y$ be any {\it topologically tame} hyperbolic $3$-manifold.
Denote by $Y'$ the manifold obtained by removing from $Y$ the geodesic
core of the Margulis tubes of a given {\it thick-thin} decomposition
of $Y$ (for instance, the one canonically associated to the Margulis
constant).  Assume that $Y'$ is not compact. By hypothesis, $Y'$ is
homeomorphic to the interior of a compact manifold with boundary, say
$W$, and each boundary component of $W$ corresponds to an ``end''
of $Y'$.  Then the ends of $Y'$ should support  {\it canonical Wick
rotations}, so that such a hyperbolic ``universe'' $Y$ concretely
realizes an {\it interaction} of its ending globally hyperbolic
Lorentzian spacetimes of constant curvature.''
\medskip

\noindent
As an application of the results of the present paper,
together with the ones of \cite{BBo}, we can concretize this idea at
least when $Y$ is {\it geometrically finite},  without accidental
parabolic and with {\it incompressible} ends. In fact, 
\cite{BBo} is responsible for
the ends that lie in the thin part of $Y$, while the present paper is
responsible for the (parabolic completion of) non-elementary ends 
of $Y$ itself.  
At the end of this Introduction we outline a few expected features
of a complete {\it WR-rescaling theory} in 3D gravity, that
should allow, in particular, to do it  for an arbitrary tame
hyperbolic $3$-manifold. This should include also an analysis of the
ending invariants of hyperbolic $3$-manifolds via WR.

In the series of papers \cite{BB} one constructs a
family of so called {\it quantum hyperbolic field theories} QHFT, that
is, roughly speaking, of determined finite
dimensional representations into the tensorial category of complex
linear spaces, of a $(2+1)$-bordism category based on compact oriented
$3$-manifolds equipped with $PSL(2,\C)$-flat bundles. There are
evidences that these theories are pertinent to 3D gravity (see also
\cite{Be}).  Tame hyperbolic $3$-manifolds furnish fundamental
examples of bordisms in this category. Our work on WR also arises
from the purpose of clarifying their full ``classical'' 3D gravity
content.
\smallskip

Let us outline the main themes of the paper.

\paragraph{Models for constant curvature geometry} 
A nice feature of 3D gravity consists in the fact that we dispose of
very explicit (local) {\it models} for the manifolds of constant
curvature, which we usually normalize to be $\kappa=0$ or $\kappa=\pm
1$.  In the Riemannian case, these are the models $\mr^3$, $\ms^3$ and
$\mh^3$ of the fundamental $3$-dimensional isotropic geometries: {\it flat},
{\it spherical} and {\it hyperbolic}, respectively.  These (in
particular the hyperbolic geometry) are the central objects of
Thurston's geometrization program, which dominates the $3$-dimensional
geometry and topology on the last decades. For the Lorentzian
signature, we have the $3$-dimensional {\it Minkowski}, {\it de
Sitter} and {\it anti de Sitter} spacetimes, respectively. We shall
denote them by $\mx_\kappa$.  Thus we can adopt the very convenient
technology of $(\X,\Gg)$-{\it manifolds}, i.e.  manifolds equipped
with (maximal) {\it special atlas} (see e.g. \cite{Thu} or {Chapter B
of \cite{BP} for more details). We recall that $\X$ denotes the model
manifold, $\Gg$ is the group of isometries of $\X$ (which possibly
preserve the orientation). A special atlas has charts with values
onto open sets of $\X$, and any change of charts is given by the
restriction to each connected component of its domain of definition
of some element $g\in \Gg$. For every $(\X,\Gg)$-manifold $M$, a very
general analytic continuation-like construction, gives us pairs
$(d,h)$, where $ d: \widetilde{M} \to \X$ is a {\it developing map}
defined on the universal covering of $M$, $ h: \pi_1(M)\to \Gg$ is a
{\it holonomy representation} of the fundamental group of $M$.
Moreover, we can assume that $d$ and $h$ are {\it compatible}, that
is, for every $\gamma \in \pi_1(M)$ we have $ d(\gamma (x))=
h(\gamma)(d(x))$, where we consider the natural action of the
fundamental group on $\widetilde{M}$, and the action of $\Gg$ on $\X$,
respectively.  The map $d$ is a local isometry, and it is unique up to
post-composition with elements $g\in \Gg$.  The holonomy
representation $h$ is unique up to conjugation by $g$. A
$(\X,\Gg)$-structure on $M$ lifts to a locally isometric structure on
$\widetilde{M}$, and these share the same developing maps. In many
situations it is convenient to consider this lifted structure, and
keep track of the isometric action of $\pi_1(M)$ on $\widetilde{M}$.
This technology is very flexible, and applies to the more general
situation where $\Gg$ is any group of real analytic transformations of
an analytic model manifold $\X$ (not necessarily a group of
isometries).  For example, we will consider {\it projective
structures} on surfaces, that are in fact $(S^2,
PSL(2,\C))$-structures, where $S^2 = {\bf P}^1(\C)$ is the Riemann
sphere, and $PSL(2,\C)$ is naturally identified with the group of
complex automorphisms of $S^2$.

\paragraph {Wick rotation} 
This is a very basic procedure for interplaying Riemannian and
Lorentzian geometry, including the {\it global causal structure} of
Lorentzian spacetimes. The simplest example applies to $\R^{n+1}$
endowed with both the standard Minkowski metric 
$-dx_0^2+\dots +dx_{n-1}^2 + dx_{n}^2$, 
and the Euclidean metric $dx_0^2+\dots + dx_{n}^2$. 
By definition (see below), these are related via a Wick
rotation directed by the vector field $\partial/\partial
x_{0}$. Sometimes one refers to it as ``passing to the imaginary
time''. More generally we have:
\begin{defi}\label{WR}{\rm  
Given a $n+1$ smooth manifold $Y$ equipped with a Riemannian metric
$g$ and a Lorentzian metric $h$, then we say
that {\it $g,\ h$ are related via a {\rm rough} Wick rotation directed by
$v$}, if:
\smallskip

(1) $v$ is a nowhere vanishing $h$-timelike and future directed
vector field on $Y$;
\smallskip

(2) For every $y\in Y$, the $g$- and $h$-orthogonal spaces to $v(y)$
coincide and we denote them by $\ort{<v(y)>}$.
\smallskip

\noindent The positive function defined on $Y$ by $\beta(y)=
-||v(y)||_h/||v(y)||_g$ is called the {\it vertical rescaling
  function} of the Wick rotation.

\noindent A Wick rotation is said {\it conformal} if there is also a
positive {\it horizontal rescaling function $\alpha$} such that, for
every 
$y\in Y$, 
$$h|_{\ort{<v(y)>}} \ =\ \alpha(y)g|_{\ort{<v(y)>}} \ . $$
}
\end{defi}
In fact, all metrics $g,\ h$ as above are canonically related by a
rough Wick rotation: we use $g$ to identify $h$ to a field of linear
automorphisms $h_y \in {\rm Aut}(TY_y)$, and we take as $v(y)$ the
field of $g$-unitary and $h$-future directed eigenvectors of $h_y$,
with negative eigenvalues.  Call $v_{(g,h)}$ this canonical vector
field associated to the couple of metrics $(g,h)$. Any other field $v$
as in Def. \ref{WR} is of the form $v=\lambda v_{(g,h)}$, for some
positive function $\lambda$. If we fix a nowhere vanishing vector
field $v$, and two positive functions $\alpha, \ \beta$ 
on $Y$ , then the conformal Wick rotation directed by $v$ and with
rescaling functions $(\alpha,\beta)$ establishes a bijection, say
$W_{(v,\beta,\alpha)}$, between the set of Riemannian metrics on
$Y$, and the set of Lorentzian metrics which have $v$ as a
timelike, future directed field. In particular, the couple $(g,v)$
encodes part of the global causality of
$h=W_{(v,\beta, \alpha)}(g)$. Clearly $W_{(v,\beta, \alpha)}^{-1}=
W_{(v,\beta^{-1}, \alpha^{-1})}$.
\smallskip

\noindent
From now on {\it we will consider only conformal WR}, so we do not
longer specify it. The couples $(g,h)$ related via a WR, and such that
{\it both $g$ and $h$ are solutions of pure gravity}, are of
particular interest, especially when the support manifold $Y$ has a
non trivial topology.

\paragraph {Rescaling directed by a vector field} 
This is a simple operation (later simply called ``rescaling'') on
metrics, which preserves the signature, and is strictly related to the
WR.  Let $k$ and $k'$ be either Riemannian or Lorentzian metrics on
$Y$.  Let $v$ be a nowhere vanishing vector field on $Y$,
$(\alpha,\beta)$ be rescaling functions as in the definition of
conformal WR. In the Lorentzian case we assume that $v$ is timelike. 
 Then $k$ and $k'$ are {\it related via a rescaling
directed by $v$, with rescaling functions $(\alpha,\beta)$}, if
\smallskip

(1) For every $y\in Y$, the $k$- and $k'$-orthogonal spaces to $v(y)$
coincide and we denote them by $\ort{<v(y)>}$.
\smallskip

(2) $k'$ coincides with $\beta k$ on the line bundle $<v>$ spanned by $v$.
\smallskip

(3) $k'$ coincides with $\alpha k$ on $\ort{<v>}$.

Again, rescalings which relate different solution of pure gravity,
possibly with different cosmological constant, are of particular
interest.

\paragraph {Canonical cosmological time} We refer
to \cite{A} for a general and careful treatment of this matter.  Here
we limit ourselves to recall that on any arbitrary spacetime $Y$ it is
defined the so-called {\it cosmological function}. Roughly speaking,
this is the proper time that every event $y\in Y$ has been in
existence, that is the upper bound of the Lorentzian lengths of the
past directed causal paths starting at $y$. In general, the
cosmological function can be very degenerate (for example, on the
Minkowski space it is constantly equal to $+\infty$). We say that a
spacetime has {\it (canonical) cosmological time} if the cosmological
fuction actually is a global time on the spacetime. In such a case,
the cosmological time coincides with the {\it finite} Lorentz distance
of every event $y$ from the {\it initial singularity} of the
spacetime. This distance is realized as the Lorentzian length of a
past directed causal curve starting at $y$, and its past limit point
defines a point $r(y)$ on the initial singularity.

\paragraph{$\Mm\Ll(\mh^2)$-spacetimes} 
An aim of this paper is also to provide a geometric explanation of the
pervasive occurrence of {\it measured geodesic laminations} on
hyperbolic surfaces (see Section \ref{laminations}) in the study of
both hyperbolic $3$-manifolds and globally hyperbolic
$(2+1)$-spacetimes of (arbitrary) constant curvature. In fact we have
placed the space $\Mm\Ll(\mh^2)$ of measured geodesic laminations on
the hyperbolic plane $\mh^2$ (the most general ones indeed) at the
centre of our discussion.
\smallskip

Let us roughly summarize only a few results obtained in the paper. 

\noindent
Given any measured geodesic lamination $\lambda=(\Ll,\mu)\in
\Mm\Ll(\mh^2)$, this encodes a maximal globally
hyperbolic spacetime $\Uu_\lambda^\kappa$ of constant curvature equal
to $\kappa$, for every $\kappa = 0,\pm 1$ - by definition they are
called $\Mm\Ll(\mh^2)$-{\it spacetimes} - and a hyperbolic
$3$-manifold $M_\lambda$, all homeomorphic to $\mh^2\times \mr$.
We could say that all these spacetimes and hyperbolic manifolds are
different ``materializations'' in 3D gravity of the same fundamental
structure $\Mm\Ll(\mh^2)$. In fact we will show that there
are natural correlations between them, given either by canonical
rescalings or WR, with {\it universal rescaling functions}.  
A delicate point to stress is that this happens in
fact {\it up to determined ${\rm C}^1$ diffeomorphism}. Note that any
${\rm C}^1$ isometry between Riemannian metrics induces at least an
isometry of the underlying {\it length spaces}. In the Lorentzian
case, it preserves the global causal structure.

More precisely we will show the following facts. Each spacetime
$\Uu_\lambda^\kappa$ has canonical cosmological time, say
$T^\kappa_\lambda$, and they  share the same initial
singularity. This is a non trivial metric space with a rich geometry
``dual'' to the geometry of the lamination $\lambda$. For $\kappa =
0,-1$, the developing map of $\Uu_\lambda^\kappa$ is an 
embedding onto a convex domain of $\mx_\kappa$. This does
not hold for $\kappa=1$ and $M_\lambda$, because the corresponding
developing map is in general not injective.  $T^0_\lambda$ is a
${\rm C}^1$ submersion and each level surface of $T^0_\lambda$
inherits a complete ${\rm C}^1$ Riemannian metric. For $\kappa = 0,1$,
$T^\kappa_\lambda(\Uu_\lambda^\kappa)=]0,+\infty[$, while for $\kappa
= -1$ it is equal to $]0,\pi[$. We adopt the following notations. For
every subset $X$ of $]0,+\infty[$, $\Uu_\lambda^\kappa(X) =
(T^\kappa_\lambda)^{-1}(X)$; for every $a\in
T^\kappa_\lambda(\Uu_\lambda^\kappa)$, $\Uu_\lambda^\kappa(a) =
\Uu_\lambda^\kappa(\{a\})$ denotes the corresponding level surface of
the cosmological time. Sometimes we will also use the notation
$\Uu_\lambda^\kappa(\geq a)$ instead of
$\Uu_\lambda^\kappa([a,+\infty[)$, and so on. Then we will explicitely
define:
\begin{enumerate}

\item A canonical rescaling which relates $\Uu_\lambda^0$
to $\Uu_\lambda^{-1}(]0,\pi/2[)$, such that:

- it is directed by the
gradient of $T^0_\lambda$; 

- the rescaling functions depend only on the value of $T^0_\lambda$.
This means that they are constant on each level surface $\Uu^0_\lambda(a)$
and their value only depends on $a$. We stress that they do not depend
on $\lambda$, so we will call them (as well as all the other ones occuring
throughout the paper) {\it universal rescaling functions}; 

- the inverse rescaling satisfies the same properties w.r.t. the
cosmological time $T^{-1}_\lambda$ restricted to
$\Uu_\lambda^{-1}(]0,\pi/2[)$.

- the rescaling extends to an identification of the respective initial
  singularity.

\item  A canonical rescaling which relates $\Uu_\lambda^0(]0,1[)$
 to $\Uu_\lambda^1$, such that: 

- it is directed by the
gradient of $T^0_\lambda$, restricted to $\Uu_\lambda^0(]0,1[)$; 

- the rescaling functions depend only on the value of $T^0_\lambda$;

- the inverse rescaling satisfies the same properties w.r.t. the
cosmological time $T^{1}_\lambda$.

- the rescaling extends to an identification of the respective initial
  singularity.

\item A canonical WR which converts $\Uu_\lambda^0(]1,+\infty[)$ into
the hyperbolic $3$-manifold $M_\lambda$ such that:

- it is directed by the
gradient of $T^0_\lambda$, restricted to $\Uu_\lambda^0(]1,+\infty[)$; 

- the rescaling functions depend only on the value of $T^0_\lambda$.

- This WR can be transported onto the slab
$\Uu_\lambda^{-1}(]\pi/4,\pi/2[)$ via the above first rescaling.  Then
the composition of WR with a developing map of $M_\lambda$ extends
continuously to the boundary of this slab, with values on
$\bar{\mh}^3 = \mh^3 \cup S^2_\infty$. The initial boundary
component $\Uu^{-1}_\lambda (\pi/4)$ corresponds to a complete end
of $M_\lambda$ and it is sent in $S^2_\infty$. In fact, this is the
developing map of a {\it asymptotic projective structure} on
$\Uu^{-1}_\lambda(\pi/4)$, and its intrinsic spacelike metric is
the {\it Thurston metric} for it. The restriction to
the final boundary component $\Uu^{-1}_\lambda (\pi/2)$ is a
locally isometric immersion onto a {\it pleated surface} in
$\mh^3$, having the measured lamination $\lambda$ as {\it
bending lamination}. Note that Wick rotations {\it cut the
initial singularity off}.

\item In the projective Klein model, the hyperbolic space $\mh^3$ and
the de Sitter space $\mx_1$ can be realized as opposite connected
components of $\PP^3(\mr)\setminus Q$, where $Q$ is a suitable
quadric.  In fact $Q$ realizes the canonical boundary at infinity
$S^2_\infty$ of $\mh^3$. It turns out that the above rescaling on
$\Uu_\lambda^0(]0,1[)$, and the WR on $\Uu_\lambda^0(]1,+\infty[)$
``fit well'' at $Q$, in the sense that they glue at $\Uu_\lambda^0(1)$
and give rise to an immersion of the whole of $\Uu_\lambda^0$ in
$\PP^3(\mr)$.
\end{enumerate}

\noindent
In the special case when the transverse measure $\mu$ of $\lambda$ is
equal to zero we say that the spacetimes are {\it static}. In
particular, $\Uu^0_\lambda$ coincides with the future $I^+(0)$ of the
origin of $\mx_0$, while $M_\lambda$ coincides with an open hyperbolic
half-space in $\mh^3$. This special case is also characterized by the
fact that the initial singularity of each spacetime
$\Uu_\lambda^\kappa$ consists of {\it one} point, and that the
cosmological times $T_\kappa$ as well as the WR or rescalings are
{\it real analytic} maps.

\paragraph {$\Gamma$-invariant constructions}
When the lamination is invariant under the action of a given discrete,
torsion-free group $\Gamma$ of isometries of $\mh^2$, that is when
$\lambda$ is the pullback of a measured geodesic lamination on the
complete hyperbolic surface $F=\mh^2/\Gamma$, then there is a natural
$\Gamma$-invariant version of the above facts. There are faithful
representations of $\Gamma$ onto groups $\widetilde{\Gamma}$ of
isometries of each $\Uu^\kappa_\lambda$ or $M_\lambda$, respectively.
The groups $\widetilde{\Gamma}$ act nicely and give rise to the
quotient spacetimes $\Uu^\kappa_\lambda/\widetilde{\Gamma}$ or the
hyperbolic $3$-manifold $M_\lambda/\widetilde{\Gamma}$, all
homeomorphic to $F\times \mr$. All constructions and structures
descend to the quotient spaces.  Strictly speaking, these spacetimes
and hyperbolic manifolds are well defined {\it up to isometry
homotopic to the identity}. We will establish procedures to select
representatives in the respective Teichm\"uller-like equivalence
classes; so we often confuse the class and its
representative. However, especially in discussing {\it convergence} of
spacetimes (see later), working up to reparametrization becomes an
important point.

\paragraph { Cocompact $\Gamma$-invariant case} 
The case of cocompact groups $\Gamma$ is of particular interest since
Mess's paper \cite{M}. This case concerns the globally hyperbolic
spacetimes homeomorphic to $S\times \mr$, for some compact surface $S$
of genus $g\geq 2$, that have been intensively investigated also in
the physics literature (see for instance \cite{W, Ca2, tH}).
We refer also to \cite{Bo1},\cite{Bo2} for some {\it volume computation}
in this case.

\noindent
In a sense, the cocompact $\Gamma$-invariant case has {\it tamest
features}. Throughout the paper, we will focus these peculiar
features, against the different phenomena that arise in general,
even when $F=\mh^2/\Gamma$ is of {\it finite area}, but non compact.
In section \ref{3cusp} we show simple examples that illustrate these
general phenomena. A key point here is that we can work with geodesic
laminations on $F$ that not necessarily have compact support.

\paragraph{The other side of $\Uu^{-1}_\lambda$ - (Broken) $T$-symmetry} 
As stated above, the WR-rescaling mechanism only concerns the past of
$\Uu^{-1}_\lambda(\pi/2)$ in the AdS spacetime
$\Uu^{-1}_\lambda$. What about the other side of this spacetime?  In
order to answer this question, it is useful to consider the spacetime
$(\Uu^{-1}_\lambda)^*$ obtained just by {\it reversing the time
orientation}. That is we study the behaviour of the AdS
$\Mm\Ll(\mh^2)$-spacetimes with respect to the $T$-{\it symmetry}. 

Recall (see Section \ref{AdS}) that the isometry group
of the spacetime $\mx_{-1}$ is isomorphic to $PSL(2,\R)\times
PSL(2,\R)$. So the AdS holonomy of a quotient of
$\Uu^{-1}_\lambda$ (when $\lambda$ is $\Gamma$-invariant), is given by a
ordered pair $(\Gamma_L,\Gamma_R)$ of representations of $\Gamma$ with
values in $PLS(2,\R)$.  In terms of the holonomy, the above spacetime
involution simply corresponds to
$$ (\Gamma_L,\Gamma_R)\to (\Gamma_R,\Gamma_L) \ .$$
\medskip

\noindent
{\it T-symmetry in the cocompact $\Gamma$-invariant case.}  Assume
that $F=\mh^2/\Gamma$ is a compact hyperbolic surface of genus $g\geq
2$, and that $\lambda$ is $\Gamma$-invariant. 
It is known since \cite{M}, that in such a case  $(\Gamma_L,\Gamma_R)$
is a pair of faithful cocompact representations of the same genus,
and that they determine the AdS spacetime. The behaviour of
$\Uu^{-1}_\lambda$ and $(\Uu^{-1}_\lambda)^*$ is {\it $T$-symmetric} in
the following sense (see Subection \ref{Tsym}):
\medskip

\noindent There are a compact hyperbolic surface $F^*=\mh^2/\Gamma^*$
of genus $g$, and a $\Gamma^*$-invariant measured geodesic lamination
$\lambda^*$ such that
$$(\Uu^{-1}_\lambda)^* = \Uu^{-1}_{\lambda^*} \ .$$
In other words: 
\smallskip

{\it Even $(\Uu^{-1}_\lambda)^*$ is $\Mm\Ll(\mh^2)$-spacetime.} 
\smallskip

\noindent
We can also say that the initial and {\it final} singularities of
$\Uu^{-1}_\lambda$ have the same kind of structure.
\medskip

\noindent
{\it Broken $T$-symmetry in the general case.}  The above $T$-symmetry
does no longer hold in general, even when $F/\Gamma$ is of finite area
but non compact:
\smallskip

{\it In general, $(\Uu^{-1}_\lambda)^*$ is no longer
a $\Mm\Ll(\mh^2)$-spacetime}. 
\smallskip

\noindent
Its initial singularity is not necessarily of the same kind.  Just for
using a somewhat suggestive terminology, let us qualify as ``naked"
the initial singularities of $\Mm\Ll(\mh^2)$-spacetimes.  Then, in
Section \ref{3cusp} we will show examples of broken $T$-symmetry given
by AdS spacetimes having naked initial singularity, whereas the final
singularity is ``censured" by BZT black holes (see \cite{Ca}).

\noindent 
In fact, the flat spacetimes $\Uu^0_\lambda$ are characterized by the
property of being regular domains in $\mx_0$ with {\it surjective
Gauss map} (see Section \ref{flat}). This has as AdS counterpart that
the spacelike surface $\Uu^{-1}_\lambda(\pi/2)$ is {\it complete}. In
general, this property is {\it not preserved up to time orientation
reversing}. We will see that $\Uu^{-1}_\lambda$ is determined by its
{\it curve at infinity} $c_\lambda$, that is the trace of the surface
$\Uu^{-1}_\lambda(\pi/2)$ on the ``boundary'' (diffeomorphic to
$S^1\times S^1$) of $\mx_{-1}$. It is an intriguing problem to
characterize AdS $\Mm\Ll(\mh^2)$-spacetimes (and broken $T$-symmetry)
in terms of $c_\lambda$.  This also depends on the subtle relationship
between these spacetimes and Thurston's ``Earthquake Theorem''
\cite{Thu2}, {\it beyond the cocompact $\Gamma$-invariant case}
already depicted in \cite{M}. We get some partial result in that
direction: we show that $c_\lambda$ is the graph of a homeomorphism of
$S^1$ iff the lamination $\lambda$ {\it generates earthquakes}. 
On the other hand, we show examples of homeomorphism
of $S^1$ such that its graph is {\it not} the $c_\lambda$ of any
$\Uu^{-1}_\lambda$ (that is such a graph determines a spacetime of
more general type). See Sections \ref{AdS} and \ref {3cusp}.

\paragraph {Along a ray of measured laminations}
As usual, let $\lambda=(\Ll,\mu)$ be a (possibly $\Gamma$-invariant)
measured geodesic lamination on $\mh^2$.  We can consider the ray of
($\Gamma$-invariant) laminations $t\lambda = (\Ll, t\mu)$, $t\in
[0,+\infty[$. So we have the corresponding $1$-parameter families of
spacetimes $\Uu^\kappa_{t\lambda }$ and hyperbolic manifolds
$M_{t\lambda}$, emanating from the static case ($t=0$).  The study of
these families gives us interesting information about the WR-rescaling
mechanism (Section \ref{der}). We study the ``derivatives'' at $t=0$
of the spacetimes $\Uu^\kappa_{t\lambda }$, and of holonomies and
``spectra'' of the quotient spacetimes.  In particular, let us denote
by $\frac{1}{t}\Uu^\kappa_{t\lambda }$ the spacetime obtained by
rescaling the Lorentzian metric of $\Uu^\kappa_{t\lambda }$ by the
constant factor $1/t^2$. So $\frac{1}{t}\Uu^\kappa_{t\lambda }$ has
constant curvature $\kappa_t=t^2\kappa$.  All the constructions made
under the normalization $\kappa =0, \pm 1$, apply straightforwardly
to every $\kappa_t$, with the obvious modifications. Then, we will
prove (using a suitable notion of convergence) that
$$ \lim_{t\to 0}\ \ \frac{1}{t}\Uu^\kappa_{t\lambda } =
 \Uu^0_\lambda \ .$$ 

Assume now that $F=\mh^2/\Gamma$ is compact. Hence we have also the
 family $(\Uu^{-1}_{t\lambda})^*= \Uu^{-1}_{\lambda_t^*}$ as in the
 previous paragraph. In such a case (see Section \ref{flat}), the set
 of $\Gamma$-invariant measured laminations has a $\R$-linear
 structure, and it makes sense to consider $-\lambda$. Then we will
 show
$$ \lim_{t\to 0}\ \ \frac{1}{t} \Uu^{-1}_{\lambda_t^*} = 
\Uu^0_{-\lambda} \ .$$

\paragraph{$\Mm\Ll(\mh^2)$-spacetimes and beyond}
The broken $T$-symmetry shows that we have developed so far only a
``sector'' of a complete WR-rescaling theory in 3D gravity.  An
immediate generalization should consist in considering {\it arbitrary}
flat regular domains in $\mx_0$, with not necessarily surjective Gauss
map.  In fact they have in general cosmological time (see Section
\ref{flat}), and the WR-rescaling formulas we have obtained for the
$\Mm\Ll(\mh^2)$-spacetimes, formally apply. However, several points
have to be clarified such as:
\smallskip

- the structure of the initial singularities and of
the (suitably generalized) dual laminations;
\smallskip

- the characterization of the AdS spacetimes 
(the hyperbolic $3$-manifolds) obtained via canonical
rescaling (canonical WR). And so on.  
\smallskip

\noindent
{\it Measured laminations on straight convex sets} defined in
\cite{Ku} seems to furnish a good generalization. Note that all
examples \ref{earthex} and in Section \ref{3cusp} present such a kind
of laminations. Moreover, it is not hard to see that the
procedure to construct a regular domain encoded by $\lambda \in
\Mm\Ll(\mh^2)$ (Section \ref{flat}), holds for these generalized
laminations.
   
We expect that {\it all} $\ (\Uu^{-1}_\lambda)^*$ should arise in such
a more general framework.  In clarifying this perspective of future
work, we have also profited of recent discussions with T. Barbot.

A completion of the WR-rescaling theory is also necessary in order to
get canonical WR on ends of arbitrary tame hyperbolic
$3$-manifolds (see the beginning of this Introduction).  We notice
that broken $T$-symmetry and WR on ends of arbitrary
hyperbolic $3$-manifolds seem to be strictly related facts.  For
instance, the completion of the hyperbolic $3$-manifolds obtained via
canonical WR on the future side of the spacetimes of Section
\ref{3cusp} are homeomorphic to a handlebody (whence having a
compressible end). These examples also show that spacetimes
of more general type can arise as ``limit'' of genuine
$\Mm\Ll(\mh^2)$-spacetimes (for instance, when $t \to +\infty$,
along certain rays of measured laminations). 
\smallskip

As a further step towards a complete theory, in \cite{BBo} we have
analyzed the simplest flat regular domain in $\mx_0$, having the most
degenerate Gauss map, that is the future, say $\Pi$, of a spacelike
geodesic line. Although this domain is very elementary, the resulting
sector of the WR-rescaling theory is far to be trivial. It is
remarkable that this can be developed in explicit and self-contained
way, eventually obtaining a complete agreement with the results of the
present paper (for instance, concerning the universal rescaling
functions). The WR-rescaling mechanism on $\Pi$ extends to so
called $\Qq\Dd$-{\it spacetimes} . These represent the 3D gravity
materialization of {\it holomorphic quadratic differentials} on
$\Omega = \C,\ \mh^2$ (possibly invariant for the action of a group
$\Gamma$ of conformal transformations).  They are suitable flat
$(\Pi,{\rm Isom}^+(\Pi))$-spacetimes (already pointed out in
\cite{BG}(3)), and the dS or AdS ones obtained via canonical
rescaling.  In particular, the quotient spacetimes of $\Pi$ with
compact space, realize all {\it non-static} globally hyperbolic flat
spacetimes with {\it toric} Cauchy surface.  Via WR we get the non-complete
hyperbolic structures on $(S^1\times S^1) \times \R$ that occur in
Thurston's {\it Hyperbolic Dehn Filling} set-up.  By canonical AdS
rescaling of the quotient spacetimes of $\Pi$ homeomorphic to
$(S^1\times \R)\times \R$, we get so called {\it BTZ black holes} (see
\cite{Ca}). In general, $\Qq\Dd$-spacetimes present world-lines of
{\it conical singularities}, and the corresponding developing maps are
{\it not} injective. So they are also a first step towards a
generalization of the theory in presence of ``particles''.
\smallskip

{\bf Acknowledgment:} Apparently our work rests on two bases: the
world discovered by Thurston, and the Lorentzian facets of that world
revealed by Mess in the germinal paper \cite{M}.

\section{Measured Geodesic Laminations}\label{laminations}

\paragraph{Laminations}
A {\it geodesic lamination} $\Ll$ on a complete Riemannian surface $F$
is a {\it closed} subset $L$ (the {\it support} of the lamination),
which is foliated by geodesics. More precisely, $L$ is covered by
boxes $B$ with a product structure $B \cong [a,b]\times [c,d]$, such
that $L\cap B$ is of the form $X\times [c,d]$, and for every $x\in X$,
$\{x\}\times [c,d]$ is a geodesic arc. Moreover, the product
structures are compatible on the intersection of two boxes.  Each
geodesic arc can be extended to a complete geodesic, i.e.  admitting
an arc length parametrization defined on the whole real line
$\R$. Either this parametrization is injective and we call its image a
{\it geodesic line} of $F$, or its image is a simple closed geodesic.
In both cases, we say that they are {\it simple} (complete) geodesics
of $F$. So, these simple geodesics make a partition of $L$, and are
called the {\it leaves} of $\Ll$.  The leaves together with the
connected components of $F\setminus L$ make a {\it stratification} of
$F$.

In this section we are interested in the geodesic laminations on the
hyperbolic plane $\mh^2$.  In this case we can prove that if a closed
subset $L$ is the disjoint union of complete geodesics, then that
geodesics give rise to a foliation of $L$ in the above sense.  In fact
let $L$ be the union of disjoint geodesics $\{l_i\}_{i\in I}$.  Fix a
point $p_0\in L$ and consider a geodesic arc $c$ transverse to the
leaf $l_0$ through $p_0$.  It is not hard to see that there exists a
neighbourhood $K$ of $p$ such that if a geodesic $l_i$ meets $K$ then
it cuts $c$. Orient $c$ arbitrarily and orient any geodesic $l_i$
cutting $c$ in such a way that respective positive tangent vectors at
the intersection point form a positive base.  Now for $x\in L\cap K$
define $v(x)$ the unitary positive tangent vector of the leaf through
$x$ at $x$.  The following lemma assures that $v$ is a $1$-Lipschitzian
vector field on $L\cap K$ (see \cite{Ep-M} for a proof).
\begin{lem}
Let $l,l'$ be disjoint geodesics in $\mh^2$. Take $x\in l$ and $x'\in
l'$ and unitary vectors $v,v'$ respectively tangent to $l$ at $x$ and
to $l'$ at $x'$ pointing in the same direction. Let $\tau(v')$ the
parallel transport of $v'$ along the geodesic segment $[x,x']$ then
\[
     || v-\tau(v') ||< d_{\mh}(x,x')
\]
where $d_{\mh}$ is the hyperbolic distance.
\end{lem}
\cvd 
\smallskip

\noindent
Thus there exists a $1$-Lipschitz vector field $\tilde v$ on
$K$ extending $v$. The flow $\Phi_t$ of this field allows us to
find a box around $p_0$.  Indeed for $\eps$ sufficiently small the
map
\[
  F: c\times (-\eps, \eps)\ni (x,t)\mapsto\Phi_t(x)\in\mh^2
\]
creates a box around $p_0$.
\smallskip

Let $\Gg$ denote the space of geodesics of $\mh^2$. It is well known
that $\Gg$ is homeomorphic to an open Moebius band (as every geodesic
is determined by the unordered pair of its distinct endpoints
belonging to $S^1_\infty$, i.e. the natural boundary at infinity of
$\mh^2$).  We say that a subset $K$ of $\Gg$ is a set of disjoint
geodesics if the geodesics in $K$ are pairwise disjoint.  Given a
geodesic lamination $\Ll$ on $\mh^2$ it is easy to see that the
subset of $\Gg$ made by the leaves of $\Ll$ is a closed set of disjoint
geodesics. Conversely, it follows from above discussion that a closed set
of disjoint geodesics gives rise to a geodesic lamination of
$\mh^2$. For simplicity, in what follows we denote by $\Ll$ both
the lamination and the set of the leaves (considered as a subset of $\Gg$).

\paragraph{Transverse measures} 
Given a geodesic lamination $\Ll$ on a complete Riemannian surface
$F$, a rectifiable arc $k$ in $F$ is {\it transverse} to the
lamination if for every point $p\in k$ there exists a neighbourhood
$k'$ of $p$ in $k$ that intersects leaves in at most a point and
$2$-strata in a sub-arc.  A {\it transverse measure $\mu$} on $\Ll$ is
the assignment of a {\it Radon} measure $\mu_k$ on each rectifiable
arc $k$ transverse to $\Ll$ (this means that $\mu_k$ assigns a finite
non-negative {\it mass} $\mu_k(A)$ to every relatively compact
Borelian subset of the arc, in a countably additive way) in such a way
that:

(1) The support of $\mu_k$ is $k\cap L$;
\smallskip

(2)  If $k' \subset k$, then $\mu_{k'} = \mu_k|_{k'}$;
\smallskip

(3) If $k$ and $k'$ are homotopic through a family of arcs transverse to
$\Ll$, then the homotopy sends the measure $\mu_k$ to $\mu_{k'}$.
\medskip

A {\it measured geodesic lamination on $F$} is a pair $\lambda=(\Ll,\mu)$,
where $\Ll$ is a geodesic lamination and $\mu$ is a transverse measure
on $\Ll$. From now on we will specialize to measured geodesic laminations
on $\mh^2$.

\begin{remark}
{\rm Notice that if $k$ is an arc transverse to a lamination of
 $\mh^2$ there exists a piece-wise geodesic transverse arc
 homotopic to $k$ through a family of transverse arcs.  Indeed there
 exists a finite subdivision of $k$ in sub-arcs $k_i$ for
 $i=1,\ldots,n$ such that $k_i$ intersects a leaf in a point and a
 $2$-stratum in a sub-arc. If $p_{i-1},p_i$ are the end-points of
 $k_i$ it is easy to see that each $k_i$ is homotopic to the geodesic
 segment $[p_{i-1},p_i]$ through a family of transverse arcs. From this
 remark it follows that a transverse measure on a lamination of
 $\mh^2$ is determined by the family of measures on transverse
 geodesic arcs.}\par
\end{remark}

The simplest example of {\it measured laminations} $\lambda=
(\Ll,\mu)$ on $\mh^2$ is given by any {\it finite} family of disjoint
geodesic lines $l_1,\dots,l_s$, each one endowed with a {\it real
positive weight}, say $a_j$. A relatively compact subset $A$ of an arc
$k$ transverse to $\Ll$ intersects it at a finite number of
points, and we set $\mu_k(A)=\sum_i a_i|A\cap l_i|$. 
\smallskip
 
Now we want to point out a different way to describe a measured
geodesic lamination on $\mh^2$.  Given an open geodesic segment $k$
transverse to $\Ll$ we denote by $\Nn(k)$ the set of geodesics
cutting $k$.  Notice that $\Nn(k)$ is an open set in $\Gg$. Moreover
the map
\[
  L\cap k\ni x\mapsto l(x)\in\Nn(k)\cap\Ll 
\]
is a homeomorphism.  Thus let $\mu^*_k$ denote the direct image of the
measure $\mu_k$.  If $k$ and $k'$ are transverse geodesics arcs we see
that there exists proper geodesic segment $k_0\subset k$ and
$k_0'\subset k'$ such that
\[
  \Nn(k)\cap\Nn(k')\cap\Ll=\Nn(k_0)\cap\Ll=\Nn(k_0')\cap\Ll \ .
\]
Notice that $k_0$ is homotopic to $k_0'$ through a family of arcs
transverse to $\lambda$. Thus it follows that
\[
\mu^*_k|_{\Nn(k)\cap\Nn(k')}=\mu^*_{k_0}=\mu^*_{k_0'}=
\mu^*_{k'}|_{\Nn(k)\cap\Nn(k')} \ .
\]
In particular we can glue the measures $\mu^*_k$ in a measure $\mu^*$
on $\Ll$. It is not hard to see that the support of $\mu^*$ is the
whole $\Ll$.  Summarizing, given a measured geodesic lamination
$\lambda=(\Ll,\mu)$, we have constructed a measure $\mu^*$ on $\Gg$
supported on $\Ll$.  Conversely if $\mu^*$ is a measure on $\Gg$ such
that the support is a set $\Ll$ of disjoint geodesics we can construct
in a similar way a measured geodesic lamination $\lambda=(\Ll,\mu)$.
\smallskip

Let us point out two interesting subsets of $\Ll$ associated to a
measured geodesic lamination $\lambda=(\Ll,\mu)$ on $\mh^2$. The {\it
simplicial part} $\Ll_S$ of $\Ll$ consists of the union of the
isolated leaves of $\Ll$. $\Ll_S$ does not depend on the measure
$\mu$. Note that in general, this is not a sublamination, that is its
support $\Ll_S$ is not a closed subset of $\mh^2$.  The {\it weighted
part} of $\lambda$ depends on the measure and it is denoted by
$\Ll_W=\Ll_W(\mu)$. In fact, it consists of the set of atoms of
$\mu^*$ (we denote by $L_W=L_W(\mu)$ its support).  It is a countable
subset of $\Gg$ but, again, it is not in general a sublamination of
$\Ll$. For instance, consider the set of geodesics $\Ll$ with a fixed
end-point $x_0\in S^1_\infty$. Clearly $\Ll$ is a geodesic lamination
of $\mh^2$ and its support is the whole of $\mh^2$.  In the half-plane
model suppose $x_0=\infty$ so that geodesics in $\Ll$ are parametrized
by $\mr$. In particular let $l_t$ denote the geodesic in $\Ll$ with
end-points $t$ and $\infty$.  If we choose a dense sequence
$(q_n)_{n\in\mn}$ in $\mr$ it is not difficult to construct a measure
on $\Ll$ such that $l_{q_n}$ carries the weight $2^{-n}$.  For that
measure $L_W$ is a dense subset of $\mh^2$.\par
 
\noindent As $L$ is the support of $\mu$, then we have the inclusion
$\Ll_S\subset\Ll_W(\mu)$. In general this is a strict inclusion.  This
terminology mostly refers to the ``dual'' geometry of the initial
singularity of the spacetimes that we will associate to the measured
geodesic laminations on $\mh^2$.  In fact, it turns out that
when $\lambda$ coincides with its simplicial part, then  
the initial singularity is a simplicial metric tree.

\paragraph{Convergence of measured geodesic laminations}
Let us denote by $\Mm\Ll(\mh^2)$ the set of measured geodesic
laminations of $\mh^2$.  Given a measured geodesic lamination
$\lambda=(\Ll,\mu)$ we can consider the positive linear functional
defined on the set of continuous functions on $\Gg$ supported on a
compact set
\[
  T_\lambda: \mathrm C_c(\Gg)\ni f\mapsto\int f\d\mu^* \in \mr \ .
\]
By Riesz Representation Theorem this correspondence is an injection of
measured geodesic laminations of $\mh^2$ into the set of positive
functionals on $\mathrm C_c(\Gg)$ and we can consider the topology on
$\Mm\Ll(\mh^2)$ induced by this map.\\ In fact we need to spell 
the notion of local convergence of a sequence of geodesic laminations
$\lambda_n=(\Ll_n, \mu_n)$.  More precisely let $K$ be a compact of
$\mh^2$ and $\Nn(K)$ denote the compact set of geodesics intersecting
$K$.  We say that $\lambda_n$ converges to $\lambda_\infty =
(\Ll_\infty,\mu_\infty)$ if
\[
    \int_{\Nn(K)} f\d\mu^*_n\rightarrow\int_{\Nn(K)} f\d\mu^*_\infty
\]
for every $f\in\mathrm C^0(\Nn(K))$.
\begin{lem}\label{laminazioni:conv:lem}
\begin{enumerate}
\item
Suppose $\lambda_n\rightarrow\lambda_\infty$ on $K$. By using a jump
function it is not hard to see that for every
$l_\infty \in \Ll_\infty\cap\Nn(K)$ there exists a sequence
$l_n\in\Ll_n\cap\Nn(K)$ such that $l_n\rightarrow l_\infty$.
\item
Suppose $\lambda_n\rightarrow\lambda_\infty$ on a compact domain $K$.
Let $k$ be a geodesic arc contained in $K$ tranverse to $\Ll_\infty$
and to $\Ll_n$ for any $n$.  Then we have
\begin{equation}\label{laminazioni:conv:eq}
   \int_k f\d(\mu_n)_k\rightarrow\int_k f\d(\mu_\infty)_k \ .
\end{equation}
\end{enumerate}
\end{lem}
\Dim
The first statement is evident. Thus let us prove only the second one.
Since $(\mu_n)_k(k)\leq\hat\mu_n(\Nn(K))$ it follows that 
up to a subsequence there exists a
measure $\mu$ on $k$ such that
\[
    \int_k f\d(\mu_n)_k\rightarrow\int_k f\d\mu \ .
\]
In order to conclude we have to show that $\mu=(\mu_\infty)_k$.\\

Given a continuous function $f$ on $k$, we can define a continuous
function $\varphi$ on $\Nn(k)\setminus \{l_0\}$, where $l_0$ is the
geodesic through $k$, by setting $\varphi(l)=f(l\cap k)$.\par If the
support of $f$ is contained in the interior of $k$ the function
$\varphi$ can be extended to $\Nn(K)\setminus l_0$. Moreover, we have
\[
   \int_k f\d(\mu_n)_k=\int_{\Nn(K)\setminus l_0} \varphi\d\hat\mu_n \ .
\]
Since $\{l_0\}$ is not contained in the support of $\mu_\infty$ the last
integral tends to the integral of $\varphi$ with respect to $\hat\mu_\infty$.
Thus if the support of $f$ is contained in ${\rm int}(k)$ then
\[
    \int_k f\d\mu=\int_k f\d(\mu_\infty)_k\ .
\]
In order to conclude it is sufficient to show that if $x$ is an
endpoint of $k$ then $\mu(\{x\})=0$. Let $k_\eps$ a neighbourhood of
$k$ in $l_0$ and fix a function $\varphi_n$ on $k_\eps$ with compact
support in its interior part such that $\varphi_n(x)=1$ and
\[
    \int_{k_\eps}\varphi_n\d(\mu_\infty)_{k_\eps}\leq1/n \ .
\]
Then we have that
\[
    \int_{k_\eps}\varphi_n\d(\mu_j)_{k_\eps}\leq 1/n
\]
for $j$ sufficiently large. In this way we obtain that
\[
     \int_k\varphi_n\d\mu\leq 1/n
\]
and this shows that $\mu(\{x\})=0$.
\cvd

\paragraph{Standard finite approximation}
Given a measured geodesic lamination $\lambda=(\Ll,\mu)$ on $\mh^2$,
we want now to construct local approximations by finite
laminations.  More precisely let us fix a geodesic segment $k$
transverse to the lamination.  Denote by $U_\eps(k)$ the
$\eps$-neighbourhood of $k$ and by $k_\eps$ the intersection of
$U_\eps(k)$ with the geodesic containing $k$.  If $\eps$ is
sufficiently small we see that every leaf of $\Ll$ that intersects
$U_\eps(k)$ must intersect $k_\eps$.  Fix $n$ and subdivide
$k_\eps$ into the union of intervals $c_1,\ldots, c_r$ such that $c_i$
has length less than $1/n$ and the end-points of $c_i$ are not in
$L_W(\mu)$.  For every $c_i$ let us set
$a_j=\mu_{k_\eps}(c_j)$. If $a_j>0$, then choose a leaf $l_j$ of $\Ll$
that cuts $c_j$. Thus consider the finite lamination $\Ll_n=\{l_j |
a_j>0\}$ and associate to every $l_j$ the weight $a_j$.  In such a way
we define a measure $\mu_n$ transverse to $\Ll_n$. Moreover, it is easy
to prove that $\lambda_n=(\Ll_n,\mu_n)$ tends to $\lambda=(\Ll,\mu)$
on $U_\eps(k)$ as $n\rightarrow+\infty$.  We call such a sequence a
{\it standard approximation} of $\lambda$.

\paragraph{$\Gamma$-invariant  measured geodesic laminations} 
Any discrete torsion-free group $\Gamma \subset \ISO(\mh^2)$ naturally
acts on the space of geodesics $\Gg$. So we can consider the
$\Gamma$-invariant geodesic laminations on $\mh^2$. More precisely,
consider the complete hyperbolic surface $F= \mh^2/\Gamma$; the
natural projection $p:\mh^2 \to F$ is a (locally isometric) universal
covering map.  Then every geodesic lamination on $F$ lifts to a
$\Gamma$-invariant geodesic lamination on $\mh^2$, and this
establishes a bijective correspondence between geodesic laminations on
$F$ and $\Gamma$-invariant geodesic laminations on $\mh^2$.

If $\lambda=(\Ll,\mu)$ is a measured geodesic lamination on $F$, then
it lifts to a measured geodesic lamination $\tilde\lambda=(\tilde\Ll,
\tilde\mu)$ on $\mh^2$ which is $\Gamma$-invariant. This means that
$\tilde\Ll$ is $\Gamma$-invariant and for every tranverse arc $k$ and
every $\gamma\in\Gamma$ we have
\[
    \mu_{\gamma(k)}=\gamma_*(\mu_k).
\]
Conversely if $\tilde\lambda$ is a $\Gamma$-invariant measured geodesic
lamination on $\mh^2$, then it induces a measured geodesic lamination on
$F=\mh^2/\Gamma$. Notice that a measured lamination $\lambda$ is
$\Gamma$-invariant if and only if the action of $\Gamma$ on $\Gg$
preserves the measure $\mu^*$.

\paragraph{ Cocompact $\Gamma$}
The case of {\it cocompact} Fuchsian groups $\Gamma$ is particularly
interesting and has been widely investigated. 
\smallskip

The simplest example of measured geodesic lamination on a compact
hyperbolic surface $F=\mh^2/\Gamma$ is a finite family of disjoint,
weighted simple closed geodesics on $F$.  This lifts to a
$\Gamma$-invariant measured lamination of $\mh^2$ made by a countable
family of weighted geodesic lines, that do not intersect each
other on the {\it whole} of $\overline{\mh}^2= \mh^2 \cup S^1_\infty$.  The
measure is defined like in the case of a finite family of weighted
geodesics. We call this special laminations {\it weighted multi-curves}.  
\smallskip

When $F=\mh^2/\Gamma$ is compact, the $\Gamma$-invariant measured
geodesic laminations $\lambda=(\Ll,\mu)$ on $\mh^2$ have particularly
good features, that do not hold in general. We limit ourselves to
remind a few of them:
\begin{enumerate}

\item
The lamination $\Ll$ is determined by its support $L$.
The support $L$ is a no-where dense set of null area.

\item 
The simplicial part $\Ll_S$ and the weighted part $\Ll_W$ 
actually coincide; moreover $\Ll_S$ is the maximal weighted
multi-curve sublamination of $\lambda$.

\item
Let us denote by $\Mm\Ll(F)$ the space of measured geodesic
laminations on $F$. It is homeomorphic to $\mr^{6g-6}$, $g\geq 2$
being the genus of $F$. Any homeomorphism $f: F\to F'$ of hyperbolic
surfaces, induces a natural homeomorphism $f_\Ll: \Mm\Ll(F)\to
\Mm\Ll(F')$; if $f$ and $f'$ are homotopic, then $f_\Ll=f'_\Ll$. This
means that $\Mm\Ll(F)$ is a topological object which only depends on
the genus of $F$, so we will denote it by $\Mm\Ll_g$.
Varying $[F]$ in the Teichm\"uller space $\Tt_g$, the above
considerations allows us to define a trivial fibre bundle 
$\Tt_g\times \Mm\Ll_g$ over $\Tt_g$ with fibre $\Mm\Ll_g$.
\end{enumerate}

\begin{exa}\label{genlam}
{\rm (1) Mutata mutandis, similar facts as above hold more generally when
$F=\mh^2/\Gamma$ is of {\it finite type} (i.e. it is homeomorphic to
the interior of a compact surface $S$ possibly with non empty
boundary), providing that the lamination on $F$ has {\it compact
support}. However, even when $F$ is of {\it finite area} but non
compact, we can consider laminations that not necessarily have compact
support (see Section \ref{3cusp}).
\smallskip

(2) Let $\gamma$ be either a geodesic line or a horocycle in $\mh^2$.
  Then, the geodesic lines orthogonal to $\gamma$ make, in the
  respective cases, two different geodesic foliations both having the
  the whole of $\mh^2$ as support. We can also define a transverse
  measure $\mu$ which induces on $\gamma$ the Lebesgue one.}  
\end{exa}

\section{Flat $\Mm\Ll(\mh^2)$-spacetimes}\label{flat}
We denote by $\mx_0$ the $3$-dimensional Minkowski space, that is
$\mr^3$ endowed with the flat Lorentzian metric
\[
h_0=-\mathrm dx_0^2+\mathrm dx_1^2+\mathrm dx_2^2 \ .
\]
Geodesics are straight lines and totally geodesic planes are affine
planes. They are classified by the restriction of $h_0$ on them: in
particular they are spacelike (resp. timelike or null) if the
restriction of $h_0$ is a Riemannian form (resp. Lorentzian or
degenerated).

The orthonormal frame
\[
   \frac{\partial\,}{\partial x_0},\,\frac{\partial\,}{\partial
     x_1},\,\frac{\partial\,}{\partial x_2}
\]
gives rise to an identification of every tangent space $T_x\mx_0$ with
$\mr^3$ provided with the Minkowskian form
\[
   \E{v}{w}=-v_0w_0+v_1w_1+v_2w_2 \ .
\]
This (ordered) framing also determines an orientation of $\mx_0$ and a
time-orientation, by postulating that $\frac{\partial\,}{\partial
x_0}$ is a future timelike vector. The isometries of $\mx_0$ coincide
with the affine transformations of $\mr^3$ with linear part preserving
the Minkowskian form. The group $ISO_0(\mx_0)$ of the isometries that
preserve both the orientation of $\mx_0$ and the time orientation,
coincides with $\mr^3\rtimes SO^+(2,1)$, where $SO^+(2,1)$ denotes the
group of corresponding linear parts. In fact, with the terminology
introduced before, we will consider $(\mx_0,ISO_0(\mx_0))$-spacetimes. 

There is a standard isometric embedding of $\mh^2$ into $\mx_0$ which
identifies the hyperbolic plane with the set of future directed
unitary timelike vectors, that is
\[
   \mh^2=\{v\in\mr^3| \E{v}{v}=-1\textrm{ and } v_0>0\}.
\]
Clearly $SO^+(2,1)$ acts by isometries on $\mh^2$, and it is not hard
to show that this action is faithful and induces an isomorphism between
$SO^+(2,1)$ and the whole group of orientation preserving isometries of
$\mh^2$.

\begin{defi}\label{flatreg}
{\rm
A  {\it flat regular domain} $\ \Uu$ in $\mx_0$ is a
convex set that coincides with the intersection of the future of its null
support planes. We also assume that  $\Uu$ has at least two null support
planes, so that $\Uu$ is neither the whole $\mx_0$ nor
the future of a null plane. Note that $\Uu$ is future complete.}
\end{defi} 

Flat regular domains in Minkowski spaces {\it of arbitrary dimension}
(defined in the same way) have been widely studied in \cite{Ba, Bo1,
Bo}. We remind here some general geometric properties of these
domains, established in the cited papers.

\begin{enumerate}  

\item Any flat regular domain $\Uu$ has canonical {\it cosmological
time}, $T$ say (see the definition given in the Introduction). This
cosmological time is a concave $\mathrm C^1$-submersion.

\item For every point $x\in\Uu$ there exists a unique point
$r(p)\in\partial\Uu$ such that $T(p)=|\E{p-r(p)}{p-r(p)}|^{1/2}$.
Moreover $r(p)$ is the unique point such that the plane orthogonal to
$p-r(p)$ at $r(p)$ is a support plane of $\Uu$.\\ The map $r:p\mapsto
r(p)$ is continuous and locally Lipschitzian.  We call it the {\it
retraction} of $\Uu$. The image $\Sigma = \Sigma_\Uu$ of that
retraction is called the {\it initial singularity} of $\Uu$.  It
coincides with the set of points in $\partial\Uu$ admitting a
spacelike support plane.

\item The gradient of $T$ is the unitary timelike vector field
\[
   \nabla_LT(p)=\frac{1}{T(p)}(r(p)-p)\ . 
\]
Notice that $\nabla_L T(p)$ is past directed, so that 
we have the map
\[
   N:\Uu\ni p \mapsto -\nabla_L T(p)\in\mh^2 
\]
which is called  the {\it Gauss map} of $\Uu$.

\item 
The $T$-level surfaces $\Uu(a)=T^{-1}(a)$ are complete Cauchy spacelike
surfaces (hence $\Uu$ is a globally hyperbolic spacetime).  They are
graphs of $1$-Lipschitz convex $\mathrm C^1$-function defined on the
horizontal plane $\{x_0=0\}$. The function $N$ restricted to $\Uu(a)$
coincides in fact with its ordinary Gauss map ($N(x)$ is the future
directed unitary normal vector to $\Uu(a)$ at $x$). One can prove that
this restriction is $1/a$-Lipschitzian w.r.t. the intrinsic distance.

\item By using the above point 2., it easy to show the following 
inequality
\begin{equation}\label{fund-ineq}
   \E{N(x)}{r(x)-r(y)}\geq 0
\end{equation}
for every $x,y\in\Uu$.

\end{enumerate}

\noindent
From now on we specialize our discussion to the 3-dimensional case.

\paragraph {From measured geodesic laminations towards 
flat regular domains}

\noindent Given a measured geodesic lamination $\lambda=(\Ll,\mu)$ on
$\mh^2$ a general construction produces a regular domain $\Uu_\lambda$
(that is $\Uu^0_\lambda$ with the notations of the Introduction) in
$\mx_0$ (see \cite{M}, \cite{Bo1,Bo}).  
\begin{figure}
\begin{center}
\input{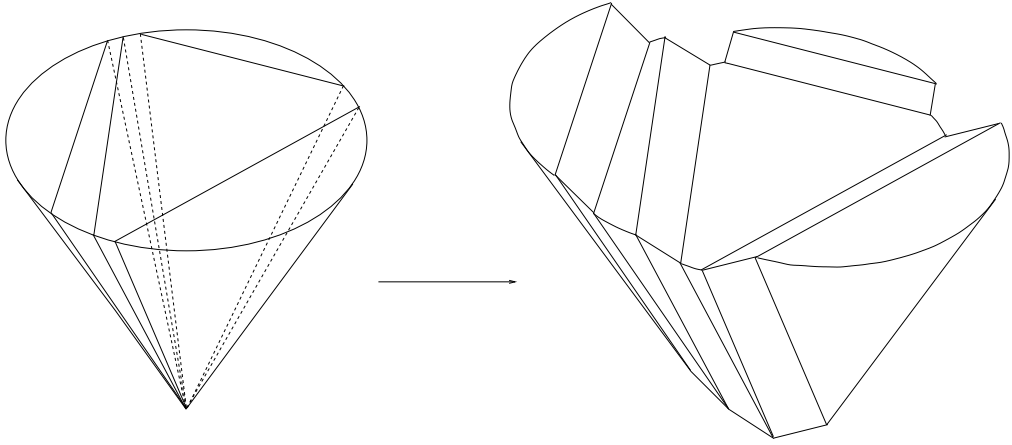}
\caption{{\small An example of regular domain associated to a finite
measured lamination.}}
\end{center}
\end{figure}
Here we sketch this
construction. Fix a base-point $x_0\in\mh^2 \setminus L_W$.  For every
$x\in\mh^2 \setminus L_W$ choose an arc $c$ transverse to $\Ll$ with
end-points $x_0$ and $x$. For $t\in c\cap L$, let $v(t)\in \R^3$ 
denote the unitary spacelike vector tangent at $t$, orthogonal to the
leaf through $t$ and pointing towards $x$. For $t\in c\setminus L$, let us
set $v(t)=0$.  Thus we have defined a function
\[
    v:c\rightarrow\mr^3
\]
that is continuous on the support of $\mu$ \ .
We can define
\[
\rho(x)=\int_{c}v(t)\d\mu(t).
\]
It is not hard to see that $\rho$ does not depend on the path $c$.
Moreover, it is constant on every stratum of $\lambda$ and it is a
continuous function on $\mh^2\setminus L_W$.  The domain $\Uu_\lambda$
can be defined in the following way
\[
  \Uu_\lambda=\bigcap_{x\in\mh^2\setminus L_W}\fut(\rho(x)+\ort{x})\ .
\]
Suppose that $x$ belongs to a weighted leaf $l\subset L_W$, with
weight $A$.  Consider the geodesic ray $c$ starting from $x_0$ and
passing through $x$. Clearly $c$ is a transverse arc. It is
not hard to see that there exists the left and right limits
\[
   \rho_-(x)=\lim_{t\rightarrow x^-}\rho(t)\ , \]
\[   \rho_+(x)=\lim_{t\rightarrow x^+}\rho(t)
\]
and that $\rho_+(x)-\rho_-(x)$ is the spacelike vector with norm equal
$A$ orthogonal to $l$. Notice that $\rho_{\pm}(x)$ depend only on the
leaf through $x$.
\begin{figure}
\begin{center}
\input{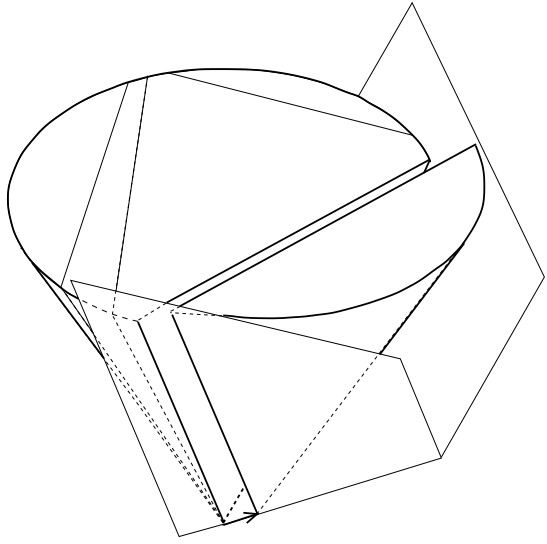}
\caption{{\small The construction of the domain associated to a finite
lamination.}}
\end{center}
\end{figure}
Another description of $\Uu_\lambda$ is the
following one
\begin{eqnarray*}
  \Uu_\lambda=\{ax+\rho(x)|x\in\mh^2-
  L_W, a>0\}\cup\\
  \cup \{x+t\rho_+(x)+(1-t)\rho_-(x)|x\in L_W, t\in[0,1], a>0\}.
\end{eqnarray*}
As before, let us denote by $T=T_\lambda$ the cosmological time, by
$N=N_\lambda$ the Gauss map and by $r=r_\lambda$ the retraction of
$\Uu_\lambda$. Then for $x\in\mh^2 \setminus L_W$ we have
\[
\begin{array}{l}
T(ax+\rho(x))=a \ ,\\
N(ax+\rho(x))=x \ ,\\
r(ax+\rho(x))=\rho(x) \ .
\end{array}
\] 
For $x\in L_W$ we have
\[
\begin{array}{l}
T(ax+t\rho_+(x)+(1-t)\rho_-(x))=a \ ,\\
N(ax+t\rho_+(x)+(1-t)\rho_-(x))=x \ ,\\
r(ax+t\rho_+(x)+(1-t)\rho_-(x))=t\rho_+(x)+(1-t)\rho_-(x) \ .
\end{array}
\] 
In particular, we have shown that 

\begin{lem}\label{Nonto} The Gauss map on $\Uu_\lambda$ is surjective.
\end{lem}

\paragraph{Measured geodesic laminations on the $T$-level surfaces}

\noindent
Let us fix a level surface $\Uu_\lambda(a)=T^{-1}(a)$ of the
cosmological time of $\Uu_\lambda$.  We want to show
that also $\Uu_\lambda(a)$ does carry a natural measured geodesic
lamination on. Consider the closed set $\hL_a=N^{-1}(L)\cap
\Uu_\lambda(a)$.  First we show that $\hL_a$ is foliated by geodesics.
If $l$ is a non weighted leaf of $L$ we have that $\hat
l_a:=N^{-1}(l)\cap\Uu_\lambda(a)=a l+\rho(x)$ where $x$ is any point
on $l$. Since the map $N:\Uu_\lambda(a)\rightarrow\mh^2$ is
$1/a$-Lipschitzian, it is easy to see that $\hat l_a$ is a geodesic of
$\Uu_\lambda(a)$.  If $l$ is a weighted leaf, then $N^{-1}(l)= al +
[\rho_-(x),\rho_+(x)]$ where $x$ is any point on $l$. Thus we have
that $N^{-1}(l)$ is an Euclidean band foliated by geodesics $\hat
l_a(t):=al +t$ where $t\in[\rho_-(x),\rho_+(x)]$.  Thus we have that
\[
  \hL_a= \bigcup_{l\subset L\setminus L_W} \hat l_a 
\cup\bigcup_{l\subset L_W}\bigcup_{t\in
  [\rho_-(x),\rho_+(x)]} \hat l_a(t).
\]
is a geodesic lamination of $\Uu_\lambda(a)$.  Given a rectifiable
arc $c$ transverse to this lamination we know that $r\circ c$ is a
Lipschitz map. Thus we have that $r$ is differentiable almost
every-where and
\[
  r(p)-r(q)=\int_c \dot r(t)\mathrm d t
\]
where $p$ and $q$ are the endpoints of $c$.  It is not hard to see
that $\dot r(t)$ is a spacelike vector. Thus we can define a measure
$\hat\mu_c=\E{\dot r(t)}{\dot r(t)}^{1/2}\mathrm dt$.  If $N(p)$ and
$N(q)$ are not in $L_W$ then
\[
  \hat\mu (c)=\mu(N(c)) \ .
\]
By this identity we can deduce that $\hat\lambda_a=(\hL_a,\hat\mu)$ is
a measured geodesic lamination on $\Uu_\lambda(a)$.  Notice in
particular that the measure $\hat\mu_c$ defined on every rectifiable
transverse arc is absolutely continuous with respect to the Lebesgue
measure of $c$. Moreover, inequality ~(\ref{fund-ineq})
implies that  $\E{\dot r(t)}{\dot r(t)}\leq\E{\dot
c(t)}{\dot c(t)}=1$ (hence the total mass of a transverse path is
bounded by the length of the path). The density of $\hat\mu$ is
bounded by $1$.

\paragraph {From flat regular domains towards measured 
geodesic laminations}
The following theorem is proved in the last chapter of \cite{Bo1}.
\begin{teo}\label{lamversusdom}
The map $\lambda\mapsto\Uu_\lambda$ establishes a bijection between
$\Mm\Ll(\mh^2)$ and the set of flat regular domains of $\mx_0$ with
surjective Gauss map, considered up to translation of $\mx_0$.
\end{teo}
For the sake of completeness, we sketch here the main steps in the
construction of the inverse map.\par Given a flat regular domain $\Uu$
provided with a cosmological time $T$, Gauss map $N$ and retraction
$r$ onto the initial singularity $\Sigma$ it is a very general fact
that the map $N$ restricted to a fiber of $r$ on $r_0\in\Sigma$ is an
isometric embedding onto an ideal convex $\Ff(r_0)$ set of
$\mh^2$.\par Given two points $r_1,r_2\in\Sigma$ we have that
$r_2-r_1$ is a spacelike vector and the goedesic $l_{1,2}$ orthogonal
to that vector can be oriented in such a way that if $x_-$
(resp. $x_+$) is a null vector corresponding to its initial
(resp. end) point then $r_2-r_1, x_-, x_+$ form a positive basis of
$\mr^3$.  Then by~(\ref{fund-ineq}) we have that $\Ff(r_1)$
(resp. $\Ff(r_2)$) is contained in the (closed) half planes on the
left (resp. on the right) of $l_{1,2}$.\par Thus we can define a
geodesic lamination $\Ll=\Ll_\Uu$ on $\mh^2$ by considering the union
of all the geodesics of the form $\Ff(r)$ and the boundary geodesics
of the ideal convex sets of the type $\Ff(r)$.\par We need the
following technical lemma to define the transverse measure
(see~\cite{Bo}).
\begin{lem}  
If $c$ is a geodesic arc in $\mh^2$ transverse
to $\Ll$ then the set $\hat c=N^{-1}(c)\cap\Uu(1)$ 
is a rectifiable arc of $\Uu(1)$.
\end{lem}
Since $N^{-1}(x)\cap\Uu(1)$ is either a point or a segment it is no
hard to see that $\hat c$ is homeomorphic to an arc.  The proof that
$\hat c$ is rectifiable is based on the following remark.  If
$p_1,p_2,p_3$ are points on $\hat c$ such that $N(p_1), N(p_2),
N(p_3)$ are ordered points on $c$ then the geodesics 
in $\mh^2$
corresponding to
$r(p_2)-r(p_1)$ and $r(p_3)-r(p_2)$ are disjoint. Thus the reverse of
the Schwarz inequality holds and so
\[
    |r(p_3)-r(p_1)|\geq|r(p_3)-r(p_2)|+|r(p_2)+r(p_1)|
\]
where $|\cdot|$ denote the Lorentzian norm.
By this inequality it follows that if $p_1,\ldots,p_n$ is a partition of
$\hat c$ then we have
\[
    \sum |p_i-p_{i-1}|\leq \ell(c)+ |r(p_n)-r(p_1)|
\]
where $\ell(c)$ is the length of $c$.\\

Let $t$ be the arc-length parameter of $\hat c$. Since $r$ is a locally
Lipschitzian map then the map
\[
     r(t)=r(\hat c(t))
\]
is Lipshitzian.  Then it is differentiable almost
every-where. Moreover by using again inequality~(\ref{fund-ineq}) we
see that $r'(t)$ is a vector tangent to $\Uu(1)$ at $\hat c(t)$,
orthogonal to the fiber of $r$ and pointing as $\hat c'(t)$.  In
particular it is a spacelike vector. If we consider the measure
$\hat\mu = |r'(t)|\d t$ then $\mu=N_*\hat\mu$ is a tranverse measure
on $c$.  Moreover if $v$ is the field on $L\cap c$ orthogonal to $L$
then we have
\[
     \int_c v\d\mu=\int_{\hat c} \hat v\d\hat\mu=\int_{\hat c}r'(t)\d t.
\]
In particular this shows that $\Uu$ is the domain associated to
$\lambda=(\Ll,\mu)$.
\cvd
\paragraph {Continuous dependence of $\ \Uu_\lambda$}  
We want to discuss now how the construction of $\Uu_\lambda$ depends
continuously on $\lambda$ (see \cite{Bo1, Bo}).
 
Fix a compact domain $K\subset\mh^2$ that contains the base point
$x_0$.  Suppose that $\lambda_n$ is a sequence of measured geodesic
laminations such that $\lambda_n\rightarrow\lambda_\infty$ on $K$.\\
We shall denote by $\Uu_n$ (resp. $\Uu_\infty$) the domain associated
to $\lambda_n$ (resp. $\lambda_\infty$) and by $T_n,r_n,N_n$
(resp. $T_\infty,r_\infty, N_\infty$) the corresponding cosmological
time, retraction and Gauss map.

In fact, we are going to outline a proof of the following proposition

\begin{prop}\label{piatto:conv:prop}
Let $K$ be a compact domain of $\mh^2$ and assume that $\lambda_n$
converge to $\lambda_\infty$ on $K$.  For any couple of positive
number $a<b$ let $U(K;a,b)$ be the set of points $x$ in
$\Uu_{\lambda_\infty}$ such that $a<T_\infty(x)<b$ and $N_\infty(x)\in
K$. We have
\begin{enumerate}
\item
  $U(K;a,b)\subset\Uu_n$ for $n>>0$;
\item
  $T_n\rightarrow T_\infty$ in $\mathrm C^1(U(K;a,b))$;
\item
  $N_n\rightarrow N_\infty$ and $r_n\rightarrow r_\infty$ uniformly on
  $U(K;a,b)$.
\end{enumerate}    
\end{prop}
  
The first simple remark is that for any $x\in K \setminus
(L_\infty)_W$ we have that
\[
   \int_{[x_0,x]}v_n(t)\d\mu_n(t)\rightarrow\rho_\infty(x)
\]
where $v_n(t)$ is the orthogonal field of $L_k$.  Such a field is
$C$-Lipschitz on $L_n\cap[x_0,x]$, for some $C$ that depends only
$K$. Thus we can extend $v_n|_{L_n\cap [x_0,x]}$ to a $C$-Lipschitz
field $\tilde v_n$ on $[x_0,x_]$. Clearly
\[
    \int_{[x_0,x]}v_n(t)\d\mu_n(t)=\int_{[x_0,x]}\tilde v_n(t)\d\mu_n(t) \ .
\]
Possibly up to passing to a subsequence, we have that $\tilde
v_n\rightarrow \tilde v_\infty$ on $\mathrm C^0([x_0,x])$. Moreover it
is not hard to see that $\tilde v_\infty=v_\infty$ on
$L_\infty\cap[x_0,x]$. Thus we have
\[
  \int_{[x_0,x]}v_n(t)\d\mu_n(t)\rightarrow\int_{[x_0,x]}\tilde
  v_\infty\d\mu_\infty(t)=\rho_\infty(x) \ .
\]
By this fact we can deduce the following result.
\begin{lem}\label{piatto:conv1:lem}
Let us take $p\in\Uu_\infty(a)$ such that $N_\infty(p)\in K$.  There
exists a sequence $p_n\in\Uu_n(a)$ such that $p_n\rightarrow
p_\infty$.
\end{lem}
\Dim
First assume that $x=N(p)\notin (L_\infty)_W$. Then we know that 
\[
  p= a x +\rho_\infty(x) \ .
\]
Hence $p_n= a x +\int_{[x_0,x]}v_n(t)\d\mu_n(t)$ works.

Now assume that $x=N(p)\in (L_\infty)_W$ so $p$ lies in a band
$al+[\rho_-(x),\rho_+(x)]$.  We can consider two points $y,z\notin
(L_\infty)_W$ such that $d(y,x)<\eps$ and $d(z, x)<\eps$ and
$||\rho(y)-\rho_-(x)||<\eps$ and $||\rho(z)-\rho_-(x)||<\eps$
(where $||\cdot||$ is the Euclidean norm).  If we
put $q^-=a y+\rho(y)$ and $q^+=z+\rho(z)$ we have that the distance of
$p$ from the segment $[q^-_\eps, q^+_\eps]$ is less than $4\eps$.  Now
let us set $q^-_n=ay+\rho_n(y)$ and $q^+_n=az+\rho_n(y)$ and choose
$n$ sufficiently large such that $||q^\pm-q^\pm_n||<\eps$.  We have
that the distance of $p$ from $[q^-_n, q^+_n||$ is less than
$6\eps$. On the other hand since the support planes for the surface
$\Uu_{\lambda_n}(a)$ at $q^-_n$ and at $q^+_n$ are near it is easy to
see that the distance of any point on $[q^-_n, q^+_n]$ from $\Uu_n(a)$
is less then $\eta(\eps)$ and $\eta\rightarrow 0$ for $\eps\rightarrow
0$.  It follows that we can take a point $p_n\in\Uu_n(a)$ arbitrarily
close to $p$ for $n$ sufficiently large.  \cvd 

Choose coordinates $(y_0,y_1,y_2)$ such that the coordinates of $x_0$
are $(1,0,0)$.  We have that the surface $\Uu_n(a)$
(resp. $\Uu_\infty(a)$) is the graph of a positive function
$\varphi_n^a$ (resp. $\varphi_\infty^a$) defined over the horizontal
plane $H=\{y_0=0\}$.  Moreover we now that $\varphi_n^a$ is a
$1$-Lipschitzian convex function and $\varphi_n^a(0)=a$.  Thus
Ascoli-Arzel\`a Theorem implies that $\{\varphi_n^a\}_{n\in\mn}$ is a
pre-compact family in $\mathrm C^0(H)$.  Up to passing to a
subsequence, there exists a function $\varphi$ on $H$ such that
$\varphi_n^a\rightarrow\varphi$ as $n\rightarrow +\infty$.
Consider the compact domain of $H$
\[
    H(K,a)=\{p\in H | N_\infty(\varphi_\infty^a(p),p)\in K\} \ .
\]
By Lemma~\ref{piatto:conv1:lem} it is easy to check that
$\varphi=\varphi_\infty^a$ on $H(K,a)$.  Thus we can deduce
\begin{equation}\label{piatto:conv2:eq}
  \varphi_n^a|_{H(K,a)}\rightarrow\varphi_\infty^a|_{H(K,a)} \ .
\end{equation}
Fix $b>a>\alpha$. We have that the domain $U(K;a,b)$ is contained in the
future of the portion of surface $N_\infty^{-1}(K)\cap\Uu_\infty(\alpha)$.
By (\ref{piatto:conv2:eq}) we see that $U(K;a,b)$ is contained in the future
of $\Uu_n(\alpha)$ for $n$ sufficiently large. Thus we have
\begin{equation}\label{piatto:conv3:eq}
    U(K;a,b)\subset\Uu_n(\geq\alpha)\qquad\textrm{ for }n>>0 \ .
\end{equation}
Since we are interested in the limit behaviour of functions
$T_n,N_n,r_n$ we can suppose that $U(K;a,b)$ is contained in
$\Uu_n(\geq\alpha)$ for any $n$.\\
Thus we have that $T_n,N_n, r_n$ are defined on $U(K;a,b)$ for any $n$.
Moreover notice that 
\[
\begin{array}{l}
    T_n(\xi,0,0)=\xi \ ,\\
    N_n(\xi,0,0)=(1,0,0) \ ,\\
    r_n(\xi,0,0)=0 \ .
\end{array}
\]
Thus we have that $r_n(p)$ lies in the half-space $H^+=\{y_0>0\}$ for
every $p\in\Uu_n$.  Since $U(K;a,b)$ is compact then there exists a
constant $C$ such that for every $p\in U(K;a,b)$ and for every past
directed vector $v$ such that $p+v$ is $H^+$ then $||v||<C$.  Since
$r_n(p)=p-T_n(p)N_n(p)$ we have that
$$||T_n(p)N_n(p)||<C$$ for every $n\in\mn$ and for every $p\in
U(K;a,b)$.  Since $T_n(p)\geq\alpha$ we can deduce the following
property.
\begin{lem}\label{piatto:Gauss:lem}
The family $\{N_n\}$ is bounded in $\mathrm C^0(U(K;a,b);\mh^2)$.
\end{lem}
\cvd Since $N_n(p)=-\nabla_L T_n(p)$ Lemma~\ref{piatto:Gauss:lem}
implies that the family $\{T_n\}$ is equi-continuous on $U(K;a,b)$. On
the other hand since $||N_n(p)||\geq 1$ we have that $|T_n(p)|<C$ for
every $p\in U(K;a,b)$.  Thus the family $\{T_n\}$ is pre-compact in
$\mathrm C^0(U(K;a,b))$.  On the other hand by using
Lemma~\ref{piatto:conv1:lem} we easily see that $T_n\rightarrow
T_\infty$.\\ Finally the same argument as in Proposition 6.5 of
\cite{Bo} shows that $N_n\rightarrow N_\infty$ in $\mathrm
C^0(U(K;a,b);\mh^2)$.  The proof of Proposition \ref{piatto:conv:prop}
easily follows.

\paragraph {$\Gamma$-invariant constructions} 
Let us assume now that $\Gamma\subset SO^+(2,1)$ is a discrete,
torsion free group of isometries of $\mh^2$, and that $\lambda$
is invariant under $\Gamma$.  We can construct a representation
\[
  f_\lambda:\Gamma\rightarrow\ISO_0(\mx_0)
\]
such that
\begin{enumerate}
\item
$\Uu_\lambda$ is $f_\lambda(\Gamma)$-invariant and the action of
$f_\lambda(\Gamma)$ on it is free and properly discontinuous.\\
\item
The Gauss map $N: \Uu_\lambda\rightarrow\mh^2$ 
is $f_\lambda$-equivariant:
\[
    N(f_\lambda(\gamma)p)=\gamma N(p).
\]
\item
The linear part of $f_\lambda(\gamma)$ is $\gamma$.
\end{enumerate}
In fact, we simply define
\[
  f_\lambda(\gamma)=\gamma\, +\, \rho(\gamma x_0) \ .
\]
Notice that $\tau(\gamma)=\rho(\gamma x_0)$ defines a {\it cocycle}
in $Z^1(\Gamma,\R^3)$; by changing the base point $x_0$,
that cocycle changes by coboundary, so we have a well defined
class in $H^1(\Gamma, \R^3)$ associated to $\lambda$.

Let us consider the hyperbolic surface $F=\mh^2/\Gamma$.
Then $X= X(\lambda,\Gamma)= \Uu_\lambda/f_\lambda(\Gamma)$ is a flat
maximal globally hyperbolic, future complete spacetime homeomorphic to
$F\times \mr$. The natural projection $\Uu_\lambda \to X$ is a locally
isometric universal covering map. The cosmological time $T$ of
$\Uu_\lambda$ descends onto the canonical cosmological time of $X$.

\paragraph{Cocompact $\Gamma$-invariant case}
The case when $\Gamma$ is a cocompact Fuchsian group has been
particularly studied, since the germinal Mess's work \cite{M}. Let us
recall a few known facts that hold in this case (see \cite {M, Bo,
Bo1}\cite{BG}(1)).

Let $X$ be any maximal globally hyperbolic, future complete flat
spacetime which is homeomorphic to $S\times \mr$, $S$ being a compact
surfaces of genus $g\geq 2$. Then:

\begin{enumerate}
\item  The linear part of the holonomy of $X$ is injective 
and its image is a cocompact Fuchsian group $\Gamma=\Gamma (X)$
such that $S$ is homeomorphic to $F=\mh^2/\Gamma$.  

\item
There is a unique $\Gamma$-invariant measured lamination $\lambda$ on
$\mh^2$, such that $X$ is isometric to $X(\lambda,\Gamma):=
\Uu_\lambda/f_\lambda(\Gamma)$.
In fact, the lamination $\lambda$ can be obtained as a
particular case of Theorem \ref{lamversusdom}, after one has proved
that the universal covering of $X$ is isometric to a flat regular
domain of $\mx_0$ with surjective Gauss map.

\item
Hence, up to isometry homotopic to the identity, the 
maximal globally hyperbolic, future complete flat
spacetime structures on $S\times \mr$ are parametrized either by:
\smallskip
 
(a) $T_g\times \Mm\Ll_g$, where $T_g$ denotes the Teichm\"uller
space of hyperbolic structures on $S$, and  $\Mm\Ll_g$ has
been introduced in Section \ref{laminations}.
\smallskip

\noindent or

\smallskip 
(b) the flat Lorentzian holonomy groups $f_\lambda(\Gamma)$'s,
up to conjugation by $ISO_0(\mx_0)$. If we fix $\Gamma$, this
induces an identification between $\Mm\Ll_g$ and the cohomology
group $H^1(\Gamma, \mr^3)$ (where $\mr^3$ is identified with the
group of translations on $\mx_0$). 

\smallskip

\noindent Moreover, $X(\lambda,\Gamma)$ is determined by the {\it
asymptotic states of its cosmological time} in the following
sense. For every $s>0$, denote by $s\ \Uu_\lambda(a)$, the surface
obtained by rescaling the metric on the level surface $\Uu_\lambda(a)
= T_\lambda^{-1}(a)$ by a constant factor $s^2$. Clearly there is a
natural isometric action of the fundamental group $\pi_1(S)\cong
\Gamma$ on each $s\ \Uu_\lambda(a)$. Then

\smallskip
(i) when $a\to +\infty$, then the action of $\Gamma$ on
$(1/a)\Uu_\lambda(a)$ 
converges (in the sense of Gromov) to the action of $\Gamma$ on
$\mh^2$;
\smallskip

(ii) when $a\to 0$, then action of $\Gamma$ on $\Uu_\lambda(a)$
converges to the action on the initial singularity of
$\Uu_\lambda$. This is a {\it real tree}, spacelike embedded into the
boundary of $\Uu_\lambda$ in $\mx_0$.  In fact this is the {\it dual}
real tree of $\lambda$, according to Skora's duality theorem (for the
notions of equivariant Gromov convergence, real tree, Skora duality
see e.g.  \cite{Ot})
\end{enumerate}

\begin{remark}\label{somma}
{\rm It follows from the above discussion, that if we fix a compact
surface $F=\mh^2/\Gamma$ of genus $g$, then we can lift on $\Mm\Ll_g$
the linear structure of $H^1(\Gamma, \mr^3)$. $\Mm\Ll_g$ is a
topological object which does not depend on the choice of $F$ (see
Section \ref{laminations}). This fact holds also for the cone
structure on $\Mm\Ll_g$ obtained by considering the ray of each
lamination $t\lambda= (\Ll,t\mu)$, $t\geq 0$. Moreover, these rays
coincide with the ones of the linear structure on $H^1(\Gamma,
\mr^3)$. However, the full linear structure {\it actually depends} on
the choice of the base surface $F$ (see \cite{Bo2}).}
\end{remark}

In Section \ref{3cusp} we show that
even when $F$ is of finite area, but non compact, we can have
radically different phenomena with respect to ones listed above for
the cocompact $\Gamma$-invariant case.


\section{WR: flat Lorentzian towards hyperbolic \\
 geometry}\label{hyp}
Let us fix a flat regular domain $\ \Uu = \Uu^0_\lambda$,
$\lambda=(\Ll,\mu)$, with surjective Gauss map $N$, according to the
results of Section \ref{flat}. If $T$ denotes its cosmological time,
remind that $\Uu(a) = T^{-1}(a)$, $\Uu(\geq 1) = T^{-1}([1,+\infty[)$,
and so on. The main aim of this section is to construct a local
$\mathrm C^1$-diffeomorphism
\[
   D=D_\lambda:\Uu(> 1)\rightarrow\mh^3
\]
such that the pull-back of the hyperbolic metric is obtained by a WR
of the standard flat Lorentzian metric, directed by the gradient of
the cosmological time of $\Uu$ (restricted to $\Uu(>1)$), and with
rescaling functions that are constant on each level surface of the
cosmological time. The hyperbolic $3$-manifold $M_\lambda$, that we
have mentioned in the Introduction, is just defined by imposing that
$D$ is a developing map of $M_\lambda$.  Later we will discuss the
asymptotic behaviour of $D$ at $\Uu(1)$. In fact, we will see that $D$
continuously extends to a map $\bar{D}$ with values in $\overline{\mh}^3 =
\mh^3 \cup S^2_\infty$, which restricts to a local $\mathrm C^1$-
embedding of $\ \Uu(1)$ into $S^2_\infty$. So, this restriction can be
regarded as a developing map of a {\it projective} structure on
$\Uu(1)$. It turns out that the spacelike metric on $\Uu(1)$ is the
Thurston metric for it. The study of the asymptotic behaviour of $D$
when the cosmological time tends to $+\infty$ will correspond to the
determination of the metric completion of $M_\lambda$.

\subsection{Bending cocycle}
We fix once for ever an embedding of $\mh^2$ into $\mh^3$
as a totally geodesic hyperbolic plane.

In order to construct the map $D$ we have to remind the construction
of {\it bending} $\mh^2$ along $\lambda$. This notion was first introduced
by Thurston in \cite{Thu}. We mostly refer to Epstein-Marden
paper \cite{Ep-M} where bending has been carefully studied. In that
paper a {\it quake-bend} map is more generally associated to every
{\it complex-valued} transverse measure on a lamination
$\lambda$. Bending maps correspond to imaginary valued measures. So,
given a measured geodesic lamination $\lambda=(\Ll,\mu)$ we will look
at the quake-bend map corresponding to the complex-valued measure
$i\mu$.  In what follows we will describe Epstein-Marden construction
referring to the cited paper for rigorous proofs.
\smallskip

Given a measured geodesic lamination $\lambda$ on $\mh^2$, first let us
recall that there is a {\it bending cocycle} associated to it. 
This is a map
\[
    B_\lambda:\mh^2\times\mh^2\rightarrow PSL(2,\mc)
\]
which verifies the following properties:
\begin{enumerate}
\item
$B_\lambda(x,y)\circ B_\lambda(y,z)=B_\lambda(x,z)$ for every $x,y,z\in\mh^2$.
\item
$B_\lambda(x,x)=Id$ for every $x\in\mh^2$.
\item
$B_\lambda$ is constant on the strata of the stratification of $\mh^2$
determined by $\lambda$
\item
If $\lambda_n\rightarrow\lambda$ on a $\eps$-neighbourhood of the
segment $[x,y]$ and $x,y \notin L_W$, then
$B_{\lambda_n}(x,y)\rightarrow B_{\lambda}(x,y)$ .
\end{enumerate}

\noindent
By definition, a {\it cocycle} on an arbitrary set $S$,
 taking values on $PSL(2,\mc)$, is a map
\[
   b:S\times S\rightarrow PSL(2,\mc)
\]
satisfying the above conditions 1. and 2.
\smallskip

If $\lambda$ coincides with its simplicial part (this notion has been
introduced in Section \ref{laminations}), then there is an easy
description of $B_{\lambda}$.\\ If $l$ is an oriented geodesic of
$\mh^3$, let $X_l\in\sG\lG(2,\mc)$ denote the infinitesimal generator
of the positive rotation around $l$ such that $\exp(2\pi X_l)=Id$
(since $l$ is oriented the notion of \emph{positive} rotation is well
defined). We call $X_l$ the standard generator of rotations around
$l$.\\ Now take $x,y\in\mh^2$. If they lie in the same leaf of
$\lambda$ then put $B_\lambda(x,y)=Id$. If both $x$ and $y$ do not lie
on the support of $\lambda$, then let $l_1,\ldots,l_s$ be the
geodesics of $\lambda$ meeting the segment $[x,y]$ and $a_1,\ldots,
a_s$ be the respective weights.  Let us consider the orientation on
$l_i$ induced by the half plane bounded by $l_i$ containing $x$ and
non-containing $y$. Then put
\[
   B_\lambda(x,y)=\exp(a_1 X_1)\circ\exp(a_2 X_2)\circ\cdots\circ\exp(a_s X_s)
\ .
\]
If $x$ lies in $l_1$ use the same construction, but replace $a_1$
by $a_1/2$; if $y$ lies in $l_s$ replace $a_s$ by $a_s/2$.

The following estimate will play an important role in our study. It is
a direct consequence of Lemma 3.4.4 (Bunch of geodesics) of
\cite{Ep-M}. We will use the operator norm on $PSL(2,\R)$.
\begin{lem}\label{hyperbolic:bend:cont:lem}
For any compact set $K$ of $\mh^2$ there exists a constant $C$ with
the following property.  Let $\lambda=(\Ll,\mu)$ be a measured
geodesic lamination on $\mh^2$.  For every $x,y\in K$ and every
geodesic line $l$ of $\Ll$ that cuts $[x,y]$, let $X$ be the standard
generator of the rotation along $l$ and $m$ be the total mass of the
segment $[x,y]$. Then we have
\[
    || B_\lambda(x,y)- \exp(m X)||\leq C m d_\mh(x,y).
\]
\end{lem}
\cvd

\smallskip

Notice that bending cocycle is not continuous on
$\mh^2\times\mh^2$. In fact by Lemma~\ref{hyperbolic:bend:cont:lem} it
follows that it is continuous on $(\mh^2\setminus
L_W)\times(\mh^2\setminus L_W)$ (recall that $L_W$ is the support of
the weighted part of $\lambda$).  Moreover if we take $x$ on a
weighted geodesic $(l,a)$ of $\lambda$ and sequences $x_n$ and $y_n$
converging to $x$ from the opposite sides of $l$ then we have that
\[
B_\lambda(x_n,y_n)\rightarrow \exp(aX)
\]
where $X$ is the generator of rotation around $l$.\\ Now we want to
define a continuous ``lifting'' of the bending cocycle on the level
surface $\Uu(1)$ of our spacetime.  More precisely we want to prove
the following proposition.
\begin{prop}\label{lift:cocycle}
A determined construction produces a {\it continuous} cocycle
\[
    \hat B_\lambda:\Uu(1)\times\Uu(1)\rightarrow PSL(2,\mc)
\]  
such that 
\begin{equation}\label{hyperbolic:bend:ext:eq}
   \hat B_\lambda (p,q)=B_\lambda(N(p), N(q))
\end{equation}
for $p,q$ such that $N(p)$ and $N(q)$ do not lie on $L_W$.
Moreover, the  map $\hat B_\lambda$ is locally Lipschitzian.
For every compact subset of $\ \Uu(1)$, the Lipschitz constant
on $K$ depends only on $N(K)$ and on the diameter of $K$.
\end{prop}    
\Dim Clearly the formula \ref{hyperbolic:bend:ext:eq}, 
defines $\hat B_\lambda$ on 
$\Uu_\lambda \setminus N^{-1}(L_W)$. We claim
that this map is locally Lipschitzian.
This follows from the following general lemma.

\begin{lem}\label{hyperbolic:bend:lip:lem}
Let $(E,d)$ be a bounded metric space, $b: E\times E \to PSL(2,\C)$
be a cocycle on $E$. Suppose there exists $C>0$ such that
\[
   ||b(x,y)-1||< C d(x,y).
\]
Then there exists a constant $H$, depending only on $C$ and on the
diameter $D$ of $E$, such that $b$ is $H$-Lipschitzian.
\end{lem}

{\it Proof of Lemma \ref{hyperbolic:bend:lip:lem}:}  
For $x,x',y,y'\in E$ we have
\[
  || b(x,y)- b(x',y')||= ||b(x,y)-b(x',x)b(x,y)b(y,y')||. 
\]
It is not hard to show that, given three elements
$\alpha,\beta,\gamma\in PSL(2,C)$ such that $||\beta-1||<\eps$ and
$||\gamma - 1||<\eps$, there exists a constant $L_\eps>0$ such that
\[
   ||\alpha-\beta\alpha\gamma||<L_\eps ||\alpha||(||\beta-1||+||\gamma-1||).
\]
Thus, we have that $||b(x,y)-1||< CD$. If
we put $\eps=C D$ we have that
\[
  || b(x,y)-b(x',y')||\leq L_\eps (D+1)C(d(x,x')+d(y,y')) \ .
\]
Thus $H=L_\eps C(D+1)$ works.
\cvd
\smallskip

\noindent 
Fix a compact subset $K$
of $\Uu(1)$ and let $C$ be the constant given by
Lemma~\ref{hyperbolic:bend:cont:lem}. Then for 
$x,x'\in K'= K\setminus N^{-1}(L_W)$ we have
\[
   ||\hat B_\lambda (x,x')-1||\leq ||\exp\hat\mu(x,x') X-1|| + 
C\hat\mu(x,x')d_\mh(N(x),N(x'))
\]
where $X$ is the generator of an infinitesimal rotation around a
geodesic of $\Ll$ cutting the segment $[N(x), N(x')]$ and
$\hat\mu(x,x')$ is the total mass of the measure along a geodesic
between $x$ and $x'$.  Thus if $A$ is the maximum of the norm of
generators of rotations around geodesics cutting $K$ and $M$ is the
diameter of $N(K)$ we get
\begin{equation}\label{hyperbolic:bend:lip3:eq}
   ||\hat B_\lambda(x,x')-1||\leq (A+CM)\hat\mu(x,x')\leq (A+CM) d(x,x') \ .
\end{equation}
By Lemma~\ref{hyperbolic:bend:lip:lem} we have that $\hat B_\lambda$
is Lipschitzian on $K'\times K'$. Moreover, since $A, C,
M$ depend only on $N(K)$, we have that the Lipschitz constant depend
only on $N(K)$ and the diameter of $K$.\\ In particular $\hat
B_\lambda$ extends to a locally Lipschitz cocycle on the closure of
$\Uu(1)\setminus N^{-1}(L_W)$ in $\ \Uu(1)$.  Notice that this
closure is obtained  by removing from $\Uu(1)$ the interior part of the
bands corresponding to leaves of $\Ll_W$.  Fix a band
$\Aa\subset\Uu(1)$ corresponding to a weighted leaf $l$. We have that
$\Aa=\{x+u|x\in l\textrm{ and } u\in [\rho_-,\rho_+]\}$.  For
$p,q\in\Aa$, let us set $p=x+u$ and $q=y+v$.  If $u=v$ then let us
put $\hat B_\lambda(p,q)=1$. Otherwise notice that $v-u$ is a vector
tangent to $\mh^2$ at $x$ and orthogonal to $l$.  Consider the
orientation on $l$ given by a positive $\pi/2$-rotation of $v-u$ in
the tangent space $T_x\mh^2$. Let $X\in\sG\lG(2,\mc)$ be the standard
generator of positive rotation around $l$ . Then for $p,q\in\Aa$ let
us put
\[
   \hat B_{\Aa}(p,q)=\exp(|v-u|X)
\]
where $|v-u|=\E{v-u}{v-u}^{1/2}$.  Notice that $\hat B_{\Aa}$ is a
cocycle. Moreover if $p,q\in\partial\Aa$, then
Lemma~\ref{hyperbolic:bend:cont:lem} implies that
\[
  \hat B_{\Aa}(p,q)= \hat B_\lambda (p,q) \ .
\]
Let us fix $p,q\in\Uu(1)$. If $p$ (resp. $q$) lies in a band $\Aa$
(resp. $\Aa'$) let us take a
point $p'\in\partial\Aa$  (resp. $q'\in\partial\Aa$) 
otherwise put $p'=p$ ($q=q'$). Then let us define
\[
   \hat B_\lambda(p,q)=\hat B_{\Aa}(p,p')\hat B_\lambda(p',q')\hat
   B_{\Aa'}(q',q).
\]
By the above remarks it is easy to see that $\hat B(p, q)$ is
well-defined, that is it does not depend on the choice of $p'$ and $q'$.
Moreover it is continuous.  Now we can prove that there exists a
constant $C$ depending only on $N(K)$ and on diameter of $K$ such that
\begin{equation}\label{hyperbolic:bend:lip:eq}
   ||\hat B_\lambda (p,q)-1||\leq C d(p,q).
\end{equation}
In fact we have found a constant $C'$ that works for
$p,q\in\Uu(1)\setminus N^{-1}(L_W)$.  On the other hand we have that
if $p,q$ lie in the same band $\Aa$ corresponding to a geodesic $l\in
\Ll_W$, then the maximum $A$ of norms of standard generators of rotations
around geodesics in $\Nn(K)$ works.  If $p$ lies in $\Uu(1)\setminus
N^{-1}(L_W)$ and $q$ lies in a band $\Aa$, then consider the geodesic
arc $c$ between $p$ and $q$ and let $q'$ in the intersection of $c$
with the boundary of $\Aa$. Then we have
\begin{eqnarray*}
  ||\hat B_\lambda(p,q)-1||=
||\hat B_\lambda(p,q')\hat B_\lambda(q',q)-1||
\end{eqnarray*}

\begin{eqnarray*}
||\hat B_\lambda(p,q')\hat B_\lambda(q',q)-1|| \leq ||\hat
    B_\lambda(p,q')-1||\,||\hat B_\lambda(q'q)||+||B_\lambda(q',q)-1||
\end{eqnarray*}

\begin{eqnarray*}
||\hat
    B_\lambda(p,q')-1||\,||\hat B_\lambda(q'q)||+||B_\lambda(q',q)-1||
   < (C'Ad(q,q')+A) d(p,q).
\end{eqnarray*}

\noindent
Thus, if $D$ is the diameter of $K$, then the constant $C''=A(C'D+1)$
works.  In the same way we can found a constant $C'''$ working for
$p,q$ that lie in different bands. Thus the maximum $C$ between
$C',C'',C'''$ works. Notice that $C$ depends only on $N(K)$ and on the
diameter of $K$. Proposition \ref{lift:cocycle} is now
proved.
\cvd

\begin{remark}\label{pippo} \emph{
By using Lemma~\ref{hyperbolic:bend:cont:lem}, we can deduce
that for a fixed $K$ in $\Uu(1)$ there exists a constant $C$ depending
only on $N(K)$ and on the diameter of $K$ such that for every
transverse arc $c$ contained in $K$ with end-points $p,q$ we have
\[
       ||\hat B_\lambda (p,q)-\exp(\hat\mu(c)X)||\leq C \hat\mu(c)
         d_\mh(N(p),N(q))
\]
where $X$ is the standard generator of a rotation around a geodesic of
$\lambda$ cutting the segment $[N(p),N(q)]$.
}\end{remark}

Let us extend now $\hat B_\lambda$ on the whole $\Uu\times\Uu$. If
$p\in \Uu$ we know that $p=r(p)+T(p)N(p)$. Let us denote by
$r(1,p)=r(p)+N(p)$ and put
\[
    \hat B_\lambda (p,q)=\hat B_\lambda(r(1,p),r(1,q)) \ .
\] 
Proposition~\ref{lift:cocycle} immediately extends to
the whole of $\Uu$.
\begin{cor}\label{hyperbolic:bend:lip:cor}
The map 
\[
  \hat B_\lambda:\Uu\times\Uu\rightarrow PSL(2,\mc)
\]
is locally Lipschitzian (with respect to the Euclidean distance on $\Uu$).
Moreover the Lipschitz constant on $K\times K$ depends only on $N(K)$, on
the diameter of $r(1,K)$ and on the maximum $M$ and minimum $m$ of $T$ on $K$.
\end{cor}
\Dim
It is sufficient to show that the map $p\mapsto r(1,p)$ is locally
Lipschitzian. Now take $p,q$ in a compact set $K\subset\Uu$.
We have that $p=r(1,p)+(T(p)-1)N(p)$ and $q= r(1,q)+(T(q)-1)N(q)$.
Thus we have
\[
  r(1,p)-r(1,q)=p-q + (N(p)-N(q)) + T(q)N(q)-T(p)N(p) \ .
\]
Since $N(K)$ is compact there exists $C$ such that $||N(p)||<C$  and 
$||N(p)-N(q)||< C |N(p)-N(q)|$ for $p,q\in K$.
Now if we set $a=T(p)$ and $b=T(q)$ we have that
$|N(p)-N(q)|<1/b|p_b-q|$ where $p_b=r(p)+bN(p)$.
It follows that 
\[
   |N(p)-N(q)|<1/m(||p-q||+||p-p_b||)=1/m(||p-q||+|T(p)-T(q)|) \ .
\]
Hence we obtain the following inequality
\[
  || r(1,p)-r(1,q)||\leq ||p-q|| + C'||p-q|| + C''|T(q)-T(p)| \ .
\]
Now since $N$ is the Lorentzian gradient of $T$ we have  that
\[
  |T(p)-T(q)|\leq C||p-q||
\]
and so the Lipschitz constant of $r_1$ is less than $1+C'+CC''$.
\cvd
\smallskip

In the last part of this subsection we will show that if
$\lambda_n\rightarrow \lambda$ on a $\eps$-neighbourhood $K_\eps$ of a
compact set $K$, then $\hat B_{\lambda_n}$ tends to $\hat B_\lambda$ on
$N^{-1}(K)$.  More precisely, for $a<b$ let $U(K;a,b)$ denote the
subset of $\Uu_\lambda$ of the points in $N^{-1}(K)$ with cosmological
time greater than $a$ and less than $b$. We know that
$U(K;a,b)\subset\Uu_{\lambda_n}$, for $n$ sufficiently large.  Then
we can consider the map $\hat B_n$ given by the restriction of $\hat
B_{\lambda_n}$ to $U(K;a,b)$.
\begin{prop}\label{hyperbolic:bend:conv:prop}
The sequence $\{\hat B_n\}$ converges to the map $\hat B=\hat
B_{\lambda}$ uniformly on $U(K;a,b)$.
\end{prop}
\Dim For $n$ sufficiently large we have $N_n(p)\in K_\eps$ for $p\in
U(K;a,b)$. Now the diameter of $N_n^{-1}(K_\eps)\cap\Uu_n(1)$ is
bounded by ${\rm diam}( K_\eps) +\mu_n^*(\Nn(K_\eps))$.  Thus we see
that there exists a constant $C$ such that every $\hat B_n$ is
$C$-Lipschitzian on $U(K;a,b)$ for $n$ sufficiently large. It follows
that the family $\{\hat B_n\}$ is pre-compact in $\mathrm
C^0(U(K;a,b)\times U(K;a,b); PSL(2,\C))$.\\ So it is sufficient to
prove that if $\hat B_n$ converges to $\hat B_\infty$ then $\hat
B_\infty=\hat B$. Clearly we have that $\hat B_\infty$ is a cocycle
and it is sufficient to show that $\hat B(p_0,q)= \hat
B_\infty(p_0,q)$. First suppose that $N(q)\notin L_W$. We can take
$q_n\in U(K;a,b)$ such that $q_n\rightarrow q$ and $N_n(q_n)$ are not
in $(L_W)_n$. Thus we have
\[
   \hat B_n(p_0, q_n)= B_n(N_n(p_0), N_n(q_n)) \ .
\]
By using Proposition~3.11.5 of \cite{Ep-M} we see that $B_n(N_n(p_0),
N_n(q_n))$ converges to $B(N(p_0), N(q))=\hat B(p_0,q)$.  Thus we have
that $\hat B(p,q)=\hat B_\infty(p,q)$ for $p,q$ lying in the closure
of $N^{-1}(\mh^2\setminus L_W)$.  Now take a point $q$ in a Euclidean
band $\Aa$. In order to conclude it is sufficient to show that $\hat
B_\infty(p,q)=\hat B(p,q)$ for $p\in\partial \Aa$ such that
$N(p)=N(q)$.  Now notice that $r_n(p)$ is different from $r_n(q)$ for
$n$ sufficiently large so $[N_n(p), N_n(q)]$ intersects the lamination
$\lambda_n$. Choose for every $n$ a leaf $l_n$ intersecting $[N_n(p),
N_n(q)]$ and let $X_n$ be the standard generator of the rotation
around $l_n$.\\ Now consider the path $c_n(t)=r_n(1, tp+(1-t)q)$. It
is not hard to see that $c_n$ is a transverse arc in $\Uu(1)$ so that
a measure $\hat\mu_n$ is defined on it. Moreover the total mass $m_n$
of $\hat\mu_n$ is
\[
     m_n=\int_0^1|\dot r_n(t)|\d t \ .
\]
By Remark~\ref{pippo} there exists a constant $C$ such that
\[
   |\hat B_n(p,q)-\exp(m_nX_n)|<Cd_\mh(N_n(p), N_n(q))\ .
\]  
On the other hand since $N_n(p)$ and $N_n(q)$ converge to $N(p)=N(q)$
it is not difficult to show that $X_n$ tends to the generator of the
rotation around the leaf through $N(q)$.  In order to conclude it is
sufficient to show that $m_n$ converges to $|r(1,p)-r(1,q)|=|p-q|$.
Now we know that
\[
tp+(1-t)q=r_n(t)+T_n(t)N_n(t)
\]
so deriving in $t$ we get
\begin{equation}\label{hyperbolic:bend:conv:eq}
p-q=\dot r_n(t)-\E{N_n(t)}{p-q}N_n(t)+T_n(t)\dot N_n(t) \ .
\end{equation}
Now we have that $N_n(t)\rightarrow N(p)$ thus $\E{N_n(t)}{p-q}$ tends
to $0$.  We will prove that $\dot N_n(t)$ tends to $0$ in
$L^2([0,1];\mr^3)$ so $\dot r_n(t)$ tends to $p-q$ in
$L^2([0,1];\mr^3)$. From this fact it is easy to see that
$m_n\rightarrow |p-q|$.  

Since the images of $N_n$ are all contained in a compact set
$\overline K_\eps\subset\mh^2$ there exists $C$ such that
\[
    \int_0^1||\dot N_n(t)||^2\d t\leq C \int_0^1 |\dot N_n(t)|^2\d t \ .
\]
On the other hand by taking the scalar product of both the hands of
equation~(\ref{hyperbolic:bend:conv:eq}) with $\dot N(t)$ we obtain
\[
 \E{p-q}{\dot N_n(t)}=\E{\dot N_n(t)}{\dot r_n(t)}+ T_n(t)|\dot N_n(t)|^2 \ .
\]
By inequality~(\ref{fund-ineq}) we can deduce $\E{\dot N_n(t)}{\dot
  r_n(t)}\geq 0$  so
\[
  \E{p-q}{\dot N_n(t)}\geq a|\dot N_n(t)|^2 \ .
\]
By integrating on $[0,1]$ we obtain that $\dot N_n$ tends to $0$ in
$L^2([0,1];\mr^3)\ $.  \cvd

\subsection{The Wick Rotation}
We are going to construct the local $\mathrm C^1$-diffeomorphism
\[
   D=D_\lambda:\Uu(>1)\rightarrow\mh^3
\]
with the properties outlined at the beginning of this Section.

Let $B=B_\lambda$ be the bending cocycle, and $\hat B = \hat
B_\lambda$ be the map defined on the whole of $\Uu \times \Uu$,
that continuously lifts $B$, as we have done above.

\noindent 
Fix $x_0$ a base point of $\mh^2$ 
($x_0$ is supposed not to be in $L_W$). The {\it bending} of $\mh^2$ along
$\tilde\lambda$ is the map 
\[
    F=F_\lambda:\mh^2\ni x\mapsto B(x_0,x)x\in\mh^3
\] 
It is not hard to see that $F$ verifies the following properties:
\begin{enumerate}
\item
It does not depend on $x_0$ up to post-composition of elements of
$PSL(2,\mc)$.
\item
It is a $1$-Lipschitz map.
\item
If $\lambda_n\rightarrow\lambda$ then $F_{\lambda_n}\rightarrow F_\lambda$
with respect to the compact open topology.
\end{enumerate}
\begin{remark}\emph{
Roughly speaking, if we bend $\mh^2$ taking $x$ fixed then $B(x,y)$ is
the isometry of $\mh^3$ that takes $y$ on the corresponding point of
the bending surface.  }\end{remark}

\noindent
Since both $\mh^3$ and $\mh^2 \subset \mh^3$ 
are oriented, the normal bundle is oriented too. Let
$v$ denote the normal vector field on $\mh^2$ that is positive
oriented with respect to the orientation of the normal bundle.
Let us take $p_0 \in N^{-1}(x_0)$ and for
$p\in\Uu(>1)$ consider the geodesic ray $c_p$ of $\mh^3$
starting from $F(N(p))$ with speed vector equal to $w(p)=\hat
B(p_0,p)_*(v(N(p)))$.  Thus $D$ is defined in the
following way:
\[
    D(p)=c_p(\arctgh(1/T(p)))=
\exp_{F(N(p))}\left(\arctgh\left(\frac{1}{T(p)}\right)w(p)\right) \ .
\]
\begin{teo}\label{hyperbolic:WR:teo}
The map $D$ is a local $\mathrm C^1$-diffeomorphism such that the
pull-back of the hyperbolic metric is equal to the Wick Rotation of
the flat Lorentz metric, directed by the gradient $X$ of the
cosmological time with universal rescaling functions (constant on the
level surfaces of the cosmological time):
    \[
       \alpha =  \frac{1}{T^2-1} \ , \qquad\qquad \beta=\frac{1}{(T^2-1)^2} \ .
    \]
\end{teo}

\begin{remark}\emph{
Before proving the theorem we want to give some heuristic motivations
for the rescaling functions we have found. Suppose $\lambda$ to be a
finite lamination.  If the weights of $\lambda$ are sufficiently small
then $F_\lambda$ is an embedding onto a bent surface of $\mh^3$. In
that case the map $D$ is a homeomorphism onto the non-convex
component, say $\Ee$ of $\mh^3\setminus F_\lambda(\mh^2)$.  The distance
$\delta$ from the boundary is a $\mathrm C^1$-submersion. Thus the
level surfaces $\Ee(a)$ give rise to a foliation of $\Ee$.  The map
$D$ satisfies the following requirement:}
\begin{enumerate}
\item
\emph{
The foliation of $\Uu$ by $T$-level surfaces is sent to that foliation of
$\Ee$.
}
\item
\emph{
The restriction of $D$ on a level surface $\Uu(a)$ is a dilation by a factor
depending only on $a$.
}
\end{enumerate}
\emph{ The first requirement implies that $\delta(D(x))$ depends only
on $T(x)$ that means that there exists a function
$f:\mr\rightarrow\mr$ such that $\delta(D(x))=f(T(x))$.}\par

\noindent
\emph{ Denote by $\mh(\lambda)$ the surface obtained by replacing
every geodesic of $\lambda$ by an Euclidean band of width equal to the
weight of that geodesic.  We have that $\Uu(t)$ is isometric to the
surface $t\mh(\lambda/t)$, where we use the same notation used in the
introduction and in the last paragraph of Section \ref{flat}.
On the other hand the surface $\Ee(\delta)$ is
isometric to $\ch\delta\mh(\lambda\tgh t)$.  Now it is not difficult
to see that $\mh(a\lambda)$ and $\mh(b\lambda)$ are related by a
dilation if and only if $a=b$.  Thus by comparing $\Ee(\delta(t))$
with $\Uu(t)$ we can deduce that }
\[
    t=1/\tgh\delta(t)
\]
\emph{
so $t>1$. Moreover the dilation factor is}
\[
    (\alpha(t))^{1/2}=\frac{\ch\delta(t)}{t}=\frac{1}{(t^2-1)^{1/2}}.
\]
\emph{
In order to compute the vertical rescaling factor notice that the hyperbolic
gradient of $\delta$ is unitary. Now, by requiring that $D$ induces a
Wick-rotation directed by the gradient $X$  of $T$, we obtain that
$D_*X=\lambda\nabla\delta$ for some function $\lambda$.
Thus we have}
\[
   \lambda=g(\nabla\delta, D_*X)=X(D^*(\delta))=\mathrm
   d(\arctgh(1/T))[X]=-\frac{1}{T^2-1}\d T(X)=\frac{1}{T^2-1}.
\]
\end{remark}

\noindent
We will prove the theorem by analyzing progressively more complicated
cases.  First we will prove it in a very special case when $\Uu$ is
the future of a spacelike segment. Then, we will deduce the theorem
under the assumption that the lamination $\lambda$ consists of a finite
number of weighted geodesic lines. Finally, by using the standard
approximation (see Section \ref{laminations}), we will obtain the full
statement.\\

Let $\Uu_0$ be the future of the segment $I=[0,\alpha_0v_0]$, where
$v_0$ is a unitary spacelike vector and $0<\alpha_0<\pi$.  If $l_0$
denotes the geodesic of $\mh^2$ orthogonal to $v_0$, the
measured geodesic lamination corresponding to $\Uu_0$ is simply
$\lambda_0=(l_0,\alpha_0)$.  We denote by $P_\pm$ the components
of $\mh^2 \setminus l_0$ in such a way that $v_0$ is outgoing from
$P_-$. It is easy to see that in this case the map
$D_0:\Uu_0\rightarrow\mh^3$ is a homeomorphism. We are going to point
out suitable $\mathrm C^1$-coordinates on $\Uu_0$ and on the image of
$D_0$ respectively, such that $D_0$ can be easily recognized with
respect to these coordinates.
\paragraph{Coordinates on $\Uu_0$}
\begin{figure}[h!]
\begin{center}
\input{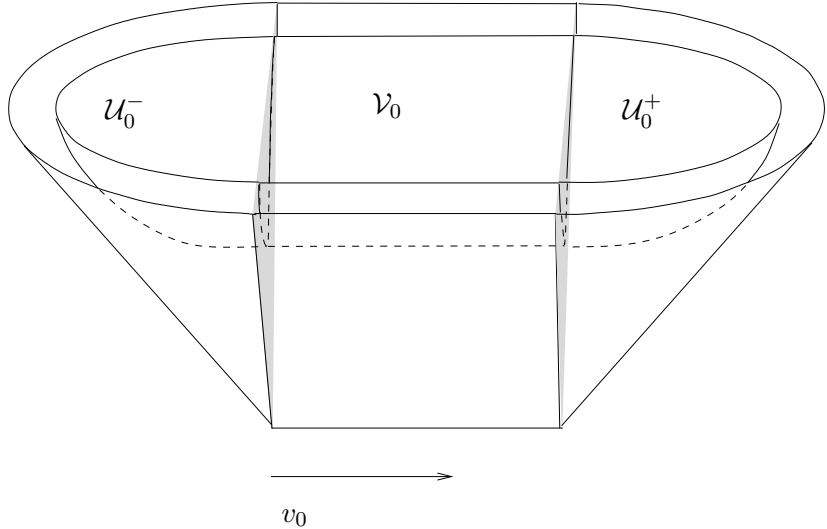}
\caption{{\small The domain $\Uu_0$ and its decomposition. Also a level
    surface $\Uu_0(a)$ is shown.}}
\end{center}
\end{figure}
As usual, let $T$ be the cosmological time, $N$ denote the Gauss map
of the level surfaces of $T$ and $r$ denote the retraction on the
singularity $I$.  We have a decomposition of $\Uu_0$ in three pieces
$\Uu_0^-,\Uu_0^+,\Vv$ defined in the following way:
\[
\begin{array}{l}
\Uu_0^-=r^{-1}(0)=N^{-1}(P_-)\\
\Vv=r^{-1}(0,\alpha_0 v_0)=N^{-1}(l_0)\\
\Uu_0^+=r^{-1}(\alpha_0 v_0)=N^{-1}(P_+) \ .
\end{array}
\]

We denote by $\Uu_0^+(a),\Uu_0^-(a),\Vv(a)$ the intersections of
corresponding domains with the surface $\Uu_0(a)$. The Gauss map on
$\Uu_0^+(a)$ (resp. $\Uu_0^-(a)$) is a diffeomorphism onto $P^+$
(resp. $P^-$) that realizes a rescaling of the metric by a constant 
factor $a^2$.  Instead, the
parametrization of $\Vv$ given by
\[
  (0,\alpha_0)\times l_0\ni(t, y)\mapsto a y+t v_0\Vv(a)
\]
produces two orthogonal geodesic foliations on $\Vv$. The parametrization
restricted to horizontal leaves is an isometry, whereas on the
on vertical leaves it acts as a rescaling of factor $a$. Thus $\Vv(a)$ is a
Euclidean band of width $a$.\\ 

Now we introduce $\mathrm C^1$
coordinates on $\Uu_0$.  We denote by $l_a$ the boundary of
$\Uu_0^{-}(a)$ and by $d_a$ the intrinsic distance of $\Uu_0(a)$. 
Fix a point $z_0$ on $l_0$ and denote by $O_a\in l_a$
the point such that $N(O_a)=z_0$.

For every $x\in\Uu_0(a)$ there is a unique point $\pi(x)\in l_a$ such
that $d_a(x,l_a)=d_a(x,\pi(x))$.  Then we consider coordinates $T,\zeta,
u$, where $T$ is again the cosmological time, and $\zeta, u$ are
defined in the following way
\[
\begin{array}{l}
\zeta(x)=\eps(x) d_{T(x)}(x,l_{T(x)})/T(x)\\
u(x)=\eps'(x)d_{T(x)}(\pi(x), O_{T(x)})/T(x) 
\end{array}
\]
where $\eps(x)$ (resp. $\eps'(x)$ ) is $-1$ if $x\in\Uu_0^{-1}$
(resp. $\pi(x)$ is on the left of $O_{T(x)}$) and is $1$ otherwise.

Choose affine coordinates of Minkowski space $(y_0,y_1,y_2)$ such that
$v_0=(0,0,1)$ and $z_0=(1,0,0)$. Thus the parametrization induced by
coordinates $T,\zeta,u$ is given by
\[
 (T,u,\zeta)\mapsto\left\{\begin{array}{ll} T(\ch u\ch \zeta,\ \sh
                     u\sh\zeta,\ \sh\zeta) & \textrm{ if }\zeta<0\\
                     T(\ch u,\ \sh u,\ \zeta) & \textrm{ if
                     }\zeta\in[0,\alpha_0/T]\\ T(\ch u\ch\zeta',\ \sh
                     u\sh\zeta',\ \sh\zeta') & \textrm{otherwise}
                     \end{array}\right.
\]
where we have put $\zeta'=\zeta-\alpha_0/T$.\par
Denote by $X$ the Lorentzian gradient of $T$. 
With respect to the frame
\[
\begin{array}{lll}
f_1=X & f_2=\frac{\partial\,}{\partial\zeta} &
f_3=\frac{\partial\,}{\partial u}
\end{array}
\]
the matrix $(h_0)_{ij}$ of the flat Lorentzian metric $h_0$ is 
diagonal, and we have: $(h_0)_{11}=-1$,  $(h_0)_{22}= T^2$,
\[
(h_0)_{33}= \left\{\begin{array}{ll}
             T^2\ch^2\zeta & \textrm{ if }\eta<0\\
             T^2           & \textrm{ if }\eta\in[0,\alpha_0/T]\\
             T^2\ch^2\zeta'& \textrm{ otherwise}
             \end{array}\right.
\]
where $\zeta'=\zeta-\alpha_0/T$.
We will also adopt the compact notation
\[
 h_0(T,\zeta, u)=\left\{\begin{array}{ll}
                           -(f^1)^2+T^2((f^2)^2+\ch^2\zeta(f^3)^2) &
                           \textrm{ if }\eta<0\\
                           -(f^1)^2 + T^2((f^2)^2+(f^3)^2) & 
                           \textrm{ if }\eta\in[0,\alpha_0/T]\\
                           -(f^1)^2+ T^2((f^2)^2+\ch^2(\zeta')(f^3)^2) &
                           \textrm{ otherwise.}
                            \end{array} \right.
\]
\paragraph{Hyperbolic Coordinates} 
We consider on $l_0$ the orientation induced by $P_-$.  Let
$X_0\in\sG\lG(2,\mc)$ be the standard generator of the rotation around
$l_0$ and consider the bent surface along $l_0$ that is the
$P=P_-\cup\exp(\alpha_0X_0)(P_+)$. The surface bounds a convex set in
$\mh^3$.  Let $\Ee_0$ denote the exterior of this convex set in
$\mh^3$.  It is not hard to see that the image of
$D_0:\Uu_0\rightarrow\mh^3$ is contained in $\Ee_0$. Now we want to
define $\mathrm C^1$-coordinates on $\Ee_0$.  We know that the
distance from $P$ is a $\mathrm C^{1,1}$ function denoted by
$\delta_0$. Moreover there is a Lipschitz retraction
\begin{figure}
\begin{center}
\input{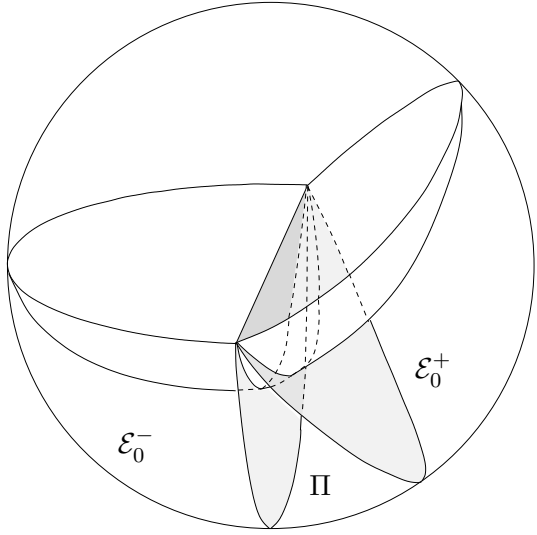}
\caption{{\small The domain $\Ee_0$ and its decomposition. A level
    surface $\Ee_0(a)$ is also shown.}}
\end{center}
\end{figure}
\[
   r:\Ee_0\rightarrow P
\]
such that $\delta(x)$ is the distance between $x$ and $r(x)$.  Let
$F^-$ and $F^+$ respectively denote the intersection of $\Ee_0$ with
the orthogonal planes to $\mh^2$ and $\exp(\alpha_0 X_0)(\mh^2)$ along
$l_0$.  We have that $\delta$ and $r$ are smooth on the complementary
regions of $F^-\cup F^+$. More precisely $\Ee_0-F^--F^+$ is the union
of $3$ dihedral angles. $P^{\pm}$ is on the boundary of one of them,
say $\Ee_0^{\pm}$, whereas denote by $\Pi$ the third one.  We have
that on $\Ee_0^-$ the function $\delta$ coincides with the distance
from $\mh^2$, on $\Ee_0^+$ it coincides with the distance from
$\exp(\alpha_0X_0)(\mh^2)$ and on $\Pi$ it coincides with the distance
from $l_0$.

Now consider the level surface $\Ee_0(a)=\delta^{-1}(a)$. It is a
$\mathrm C^1$-surface moreover the intersections with $\Ee_0^+$,
$\Ee_0^-$ and $\Pi$ give rise to a decomposition of $\Ee_0(a)$ in
three pieces, say $\Ee_0^+(a)$, $\Ee_0^-(a)$ and $\Pi(a)$.  The
retraction on $\Ee_0^+(a)$ and $\Ee_0^-(a)$ is a scaling of a factor
$\ch a$ so they are isometric to $\ch\alpha_0\cdot P_\pm$.
 On the other hand the retraction on $\Pi(a)$ is not
injective. Consider for any $x\in\Pi(a)$ the angle $\alpha(x)$ that
the ray $[r(x),x]$ forms with the normal to $\mh^2$. Then we have that
the map
\[
\Pi(a)\ni x\mapsto (r(x),\alpha(x))\in [0,\alpha_0]\times l_0 
\]
is a diffeomorphism. Moreover it induces a geodesic orthogonal product
structure on $\Pi(a)$. The map $r$ restricted to vertical leaves is a
scaling of a factor $1/\ch a$ whereas the map $\alpha$ restricted to
horizontal leaves is a scaling of a factor $1/\sh a$.  Thus we have that
$\Pi(a)$ is a Euclidean band of width equal to $[0,\alpha_0\sh a]$.\\
Now let $d_a$ denote the intrinsic distance on the surface $\Ee_0(a)$
and $l_a$ be boundary geodesic of $\Ee_0^-(a)$. Let us consider the
point $z_a\in l_a$ such that $r(z_a)=z_0$. (Notice that $a\mapsto z_a$
parameterizes the geodesic starting from $z_0$ orthogonal to
$\mh^2$. Now it is easy to see that for every $x\in\Ee_0(a)$ there is
a unique point $\pi(x)\in l_a$ such that
$d_a(x,l_a)=d_a(x,\pi(x))$. Thus we consider on $S_a$ the coordinates
$\delta,\eta,s$ where $\delta$ is the distance function and the others
are defined in the following way
\[
\begin{array}{l}
\eta(x)=\eps(x)d_{\delta(x)}(x, l_{\delta(x)})/\ch(\delta(x))\\
s(x)=\eps'(x)d_{\delta(x)}(\pi(x),z_{\delta(x)})/\ch(\delta(x))
\end{array}
\]
where $\eps(x)$ (resp. $\eps'(x)$) is $-1$ if $x\in\Ee_0^-$
(resp. $\pi(x)$ is on the left of $z_{\delta(x)}$) and is $1$
otherwise.  In the hyperboloid model of $\mh^3$ let us identify
$\mh^2$ to the hyperbolic plane orthogonal to $(0,0,0,1)$ and $l_0$ to
the geodesic with end-points $(1,-1,0,0)$ and $(1,1,0,0)$.  The
parametrization of $\Ee_0$ induced by $\delta,\zeta,s$ is the
following one
\[
(s,\eta, \delta)\mapsto\left\{\begin{array}{ll}
                          \ch\delta\left(\ch\eta\ch s,\ \ch\eta\sh s,\
                          \sh\eta,\ 0\right)\ +\ \sh\delta(0,0,0,1) \\
                          \textrm{if }\eta\leq 0\ ;\\ \ch\delta\left(\ch
                          s,\ \sh s,\ 0,0\right)\ + \
                          \sh\delta\left(0,\ 0,\
                          -\sin\frac{\eta}{\tgh\delta},\,
                          \cos\frac{\eta}{\tgh\delta}\right) \\
                          \textrm{if }\eta\in[0,\alpha_0\tgh \delta]\ ;\\

                          \ch\delta\left(\ch\eta'\ch s,\ \ch\eta'\sh
                          s,\ \sh\eta'\cos\alpha_0,\
                          \sh\eta'\sin\alpha_0\right)\ + & \\
                          +\sh\delta(0,0,-\sin\alpha_0,\
                          \cos\alpha_0)\\ \textrm{if
                          }\eta\geq\alpha_0\tgh \delta.
                           \end{array}\right.
\] 
where $\eta'=\eta-\alpha_0\tgh \delta$.\par
Let $Y$ denote the gradient of $\delta$. With
respect to the frame
\[
\begin{array}{lll}
e_1=Y & e_2=\frac{\partial\,}{\partial\eta} & e_3=\frac{\partial\,}{\partial
  s}
\end{array}
\]
the hyperbolic metric is written in this way
\[
 g(\delta,\eta, s)=\left\{\begin{array}{ll}
                           (e^1)^2+\ch^2\delta((e^2)^2+\ch^2\eta(e^3)^2) &
                           \textrm{ if }\eta<0\\
                           (e^1)^2 + \ch^2\delta((e^2)^2+(e^3)^2) & 
                           \textrm{ if }\eta\in[0,\alpha_0\tgh \delta]\\
                           (e^1)^2+\ch^2\delta((e^2)^2+\ch^2(\eta')(e^3)^2) &
                           \textrm{ otherwise}
                            \end{array}\right.
\]
where $\eta'=\eta-\alpha_0\tgh\delta$.

\paragraph{Local map}
Denote by $\hat B_0$ the cocycle associated to the domain $\Uu_0$.
An easy computation shows that $\hat B_0((T,u,\zeta),(1,0,0))$ is the identity
if $\zeta<0$, and is the positive rotation around $l_0$ of angle
$\min(T\zeta,\alpha_0)$ otherwise.
Thus if we put $d={\rm arcth} 1/T$ we can deduce
\[
    D_0(T,u,\zeta)=\left\{\begin{array}{ll} \ch d\left(\ch\zeta\ch u,\
                          \ch\zeta\sh u,\ \sh\zeta,\ 0\right)\ +\ \sh
                          d(0,0,0,1) \\ \textrm{if }\zeta\leq 0\ ;\\
                          \ch d\left(\ch u,\ \sh u,\ 0,0\right)\ + \
                          \sh d\left(0,\ 0,\ -\sin T\zeta,\, \cos
                          T\zeta \right) \\ \textrm{if
                          }\zeta\in[0,\alpha_0/T]\ ;\\ \ch
                          d\left(\ch\zeta'\ch u,\ \ch\zeta'\sh u,\
                          \sh\zeta'\cos\alpha_0,\
                          \sh\zeta'\sin\alpha_0\right)\ + & \\ +\sh
                          d(0,0,-\sin\alpha_0,\ \cos\alpha_0)\\
                          \textrm{if }\zeta\geq\alpha_0/T.
                           \end{array}\right.
\]
where $\zeta'=\zeta-\alpha_0/T$.
With respect to the coordinates
$(T,\zeta,u)$ and $(\delta,\eta,s)$ the function
$D_0:\Uu_0\rightarrow\mh^3$ takes the form
\[
\left\{ \begin{array}{l}
         \delta(D_0(p))={\rm arcth}(1/T(p))\\
         \eta(D_0(p))=\zeta(p)\\
          s(D_0(p))=u(p)\ .
        \end{array}\right.
\]
Thus $D_0$ is a $\mathrm C^1$ homeomorphism, $\mathrm C^2$-almost
every-where with second derivative locally bounded.
Now we want to compute the pull-back $D_0^*(g)$ of the hyperbolic metric.
First notice that
\[
\begin{array}{l} 
(D_0)_*(\frac{\partial}{\partial\zeta})=\frac{\partial}{\partial\eta}\\
(D_0)_*(\frac{\partial}{\partial u})=\frac{\partial}{\partial s}\ .
\end{array}
\]
Now notice that the flow lines of $X$ are transformed in the flow lines of $Y$
so $(D_0)_*(X)=f Y$. In order to compute $f$ notice that
\begin{eqnarray*}
g(Y,(D_0)_*(X))=(D_0)_*(X)(\delta)=X(D_0^*(\delta))=\\
=X({\rm arcth}(1/T))=1/(1-T^2)X(T)=1/(T^2-1)\ .
\end{eqnarray*}
Thus we have
\[
(D_0)_*(X)=\frac{1}{T^2-1} Y.
\]
Finally, an easy computation shows that
\[
  (D_0)^*(g)=W_{(X,\frac{1}{(T^2-1)^2}, \frac{1}{T^2-1})}(h_0)
\]

\paragraph{Finite laminations} 
Suppose that $\lambda$ is a finite lamination on $\mh^2$.
We want to reduce this case to the previous one. In fact,
we will show that for any $p\in\Uu_\lambda$ there exist a
small neighbourhood $U$ in $\Uu_\lambda$, an isometry $\gamma$ of $\mx_0$ and
an isometry $\sigma$ of $\mh^3$ such that
\begin{enumerate}
\item
$\gamma(U)\subset\Uu_0$.
\item
   $\gamma$ preserves the cosmological times, that is
  \[
     T(\gamma(p)) = T_\lambda(p)
  \] 
  for every $p\in U$
\item 
 We have
\begin{equation}\label{1}
    \sigma\circ D_\lambda(p)=D_0\circ \gamma(p)
\end{equation}
for every $p\in U$.
\end{enumerate}

\noindent
First suppose that $p$ does not lie in any Euclidean band.  Fix 
$\eps>0$ so that the disk $B_\epsilon(x)$ in $\mh^2$, with center at
$x=N(p)$ and ray equal to $\eps$, does not intersect any geodesic of
$\lambda$.
Thus we can choose an isometry $\gamma_0$ of $\mh^2$ such that the
distance between $z=\gamma_0(x)$ and $z_0$ is less than $2\eps$.  Now
let us set $U=N^{-1}(B_\eps(x))$ and $\gamma =\gamma_0 - r(p)$
(where $r$ denotes here the retraction of $\Uu_\lambda$).  Finally
let us set $\sigma=\hat B_\lambda(p_0,p)$.  In fact we have
\[
  D(\xi)= \hat B_\lambda(p_0, p)\circ D_0\circ \gamma(\xi)
\]
for $\xi\in U$.\\

\noindent
If $p$ lies in the interior of a band $\Aa=r^{-1}[r_-,r_+]$
corresponding to a weighted leaf $l$, let $I$ be an interval $[s_-,
s_+]$ contained in $[r_-,r_+]$ centered in $r(p)$ of length less than
$\alpha_0$ and set $U=r^{-1}(I)$.  Let $q$ be a point in $U$ such
that $r(q)=s_-$ and $N(q)=N(p)$.  Let $\gamma_0$ be an isometry of
$\mh^2$ which sends $N(q)$ onto $z_0$ and $l$ onto $l_0$; set
$\gamma =\gamma_0-r(q)$.  Now for $\xi\in U$ we have
\[
  D_\lambda(\xi)= \hat B_\lambda(p_0,q)D_0(\gamma \xi).
\]

\noindent
Finally suppose that $p$ lies on the boundary of a band
$\Aa=r^{-1}[r_-,r_+]$ corresponding to the weighted leaf $l$.  Without
lost of generality we can suppose that $r(p)=r_-$.  Now let us fix a
neighbourhood $U$ of $p$ that does not intersect any other Euclidean
band and such that $r(U)\cap [r_-, r_+]$ is a proper interval of
length less than $\alpha_0$.  Let $\gamma_0\in PSL(2,\mr)$ send $N(p)$
onto $z_0$ and $l$ onto $l_0$; set  $\gamma = \gamma_0-r(p)$.  Also in this
case we have
\[
  D_\lambda(\xi)= \hat B_\lambda(p_0,p)D_0(\gamma \xi).
\]
\paragraph{General case}
Before extending this result to the general case we need some remarks
about the regularity of $D_\lambda$, when $\lambda$ is finite. Clearly
it is a smooth function outside the boundaries of the Euclidean
bands. Moreover, we need the following estimate. We use the notations
of the proof of Proposition \ref{hyperbolic:bend:conv:prop}.
\begin{lem}\label{hyperbolic:WR:lem}
Fix a bounded domain $K$ of $\mh^2$, a bounded domain $\Dd$ of
$\Mm\Ll(\mh^2)$ and $1<a<b$. For every finite lamination $\lambda$ in $\Dd$
let us set $U_\lambda=U_{\lambda}(K;a,b)$.
Then there exists a constant $C$ depending only on $K$, $\Dd$ and
$a,b$ such that the first and the second derivatives of $D_\lambda$ on
$U_\lambda$ are bounded by $C$.
\end{lem}
\Dim For every point $p$ of $U_\lambda$, the above construction
gives us a neighbourhood $W$, an isometry $\gamma_W$ of $\mx_0$,
and an isometry $\sigma_W$ of $\mh^3$ such that
\[
   D_\lambda=\sigma_W\circ D_0\circ \gamma_W \ .
\]
Moreover, we can choose $W$ so small in such a way that $\gamma_W$ is
contained in $U_0(B_{2\alpha_0}(z_0);a,b)$. Let us fix a constant $C'$
such that first and second derivatives of $D_0$ are bounded by $C'$ on
this set.  Moreover it is easy to see that the families $\{\gamma_W\}$
and $\{\sigma_W\}$ are bounded sets. Hence there exists a constant
$C''$ such that every first and second derivative of both $\gamma_W$ and
$\sigma_W$ are bounded by $C''$.  Thus we see that first and second
derivatives of $D_\lambda$ are bounded by $C=27(C'')^2C'$.  
\cvd

\noindent We can finally prove Theorem \ref{hyperbolic:WR:teo} in the
general case.
Take a point $p\in\Uu_\lambda$ and consider a sequence of standard
approximations $\lambda_n$ of $\lambda$ on a neighbourhood $K$ of the
segment $[N(p_0), N(p)]$.  It is not hard to see that $D_{\lambda_n}$
converges to $D_\lambda$ on $U(K;a,b)$.  On the other hand by
Lemma~\ref{hyperbolic:WR:lem} we have that $D_{\lambda_n}$ is a
pre-compact family in $\mathrm C^1(U(K;a,b); \mh^3)$.  Thus it follows
that the limit of $D_{\lambda_n}$ is a $\mathrm C^1$-function.  
Finally, as $D_{\lambda_n}$ ${\rm C}^1$-converges to $D_\lambda$,
the cosmological time of $\Uu_{\lambda_n}$  ${\rm C}^1$-converges
on $U(K;a,b)$ to the one of $\Uu_\lambda$, and the pull-back on 
$\Uu_{\lambda_n}$ of the hyperbolic metric is obtained via the
determined WR, the same conclusion holds on $\Uu$.
\cvd

\subsection{On the geometry of $M_\lambda$}
Let $M_\lambda$ be the hyperbolic $3$-manifold arising by performing
the Wick Rotation described in Theorem~\ref{hyperbolic:WR:teo}. So
$M_\lambda$ consists of the domain $\Uu_\lambda(>1)$ endowed
with a determined hyperbolic metric, say $g_\lambda$. We are going to
study some geometric properties of $M_\lambda$. As usual, $T$ denotes
the cosmological time of the spacetime $\Uu^0_\lambda$, and $N$
its Gauss map.\\

\paragraph{Completion of $M_\lambda$}
Let us denote by $\delta$ the length-space-distance on $M_\lambda$
associated to $g_\lambda$. We want to determine the metric completion
$\overline M_\lambda$ of $(M_\lambda,\delta)$. 
The following theorem summarizes the main features
of the geometry of  $\overline M_\lambda$. The rest of the paragraph
is devoted to prove it. 

\begin{teo}\label{hyperbolic:compl:teo}
(1) The completion of $M_\lambda$ is 
$\overline M_\lambda=M_\lambda\cup\mh^2$ endowed with the
distance $\overline\delta$ 
\[
\begin{array}{ll}
   \overline\delta (p,q)=\delta(p,q) & \textrm{ if }p,q\in M_\lambda\\
   \overline\delta (p,q)=d_\mh(p,q)  & \textrm{ if }p,q\in\mh^2\\
   \overline\delta (p,q)=\lim_{n\rightarrow+\infty}\delta(p,q_n) 
& \textrm{ if } p\in M_\lambda\textrm{and }q\in\mh^2
\end{array}
\]
where $(q_n)$ is any sequence in $\Uu_\lambda$ such that $T(q_n)=n$ and
$N(q_n)=q$. The copy of $\mh^2$ embedded into $\overline M_\lambda$
is called the {\rm hyperbolic boundary}
$\partial_h M_\lambda$ of $M_\lambda$. 
\smallskip

(2) The developing map $D_\lambda$ continuously extends to a map
defined on $\overline M_\lambda$. Moreover, the restriction of
$D_\lambda$ to the hyperbolic boundary $\partial_h 
M_\lambda$ coincides with the bending map $F_\lambda$.
\smallskip

(3) Each level surface of the cosmological time $T$ restricted to
$\Uu(>1)$ is also a level surface in $\overline M_\lambda$
of the distance function $\Delta$ from its hyperbolic boundary 
$\partial_h M_\lambda$. Hence the inverse WR is directed
by the gradient of $\Delta$.
\smallskip

(4) $\overline M_\lambda$ is a topological manifold with boundary,
 homeomorphic to $\mr^2\times[0,+\infty)$. Moreover, $\overline
 M_\lambda(\Delta\leq\eps)$ is a collar of $\mh^2=\partial_h M_\lambda$.
\end{teo}

For simplicity, in what follows we denote by $\delta$ both the
distance on $M_\lambda$ and the distance on $M_\lambda\cup\mh^2$.

We are going to establish some auxiliary results.
\begin{lem}\label{hyperbolic:compl:lem}
The map $N:M_\lambda\rightarrow\mh^2$  is $1$-Lipschitzian.\\
\end{lem}
\Dim Let $p(t)$ be a $\mathrm C^1$-path in $M_\lambda$. We have to
show that the length of $N(t)=N(p(t))$ is less than the length of
$p(t)$ with respect to $g_\lambda$.  (Since $N$ is locally
Lipschitzian with respect to the Euclidean topology $N(t)$ is a
Lipschitzian path in $\mh^2$.)\\ By deriving the equality
\[
   p(t)= r(t)+ T(t)N(t)
\]
we get
\[
  \dot p(t) =\dot r(t) + \dot T(t) N(t) + T(t)\dot N(t)
\]
As $\dot r$ and $\dot N$ are orthogonal to $N$ (that up to the sign is the
gradient of $T$) we have
\begin{equation}\label{hyperbolic:compl:eq}
  g_\lambda(\dot p(t), \dot p(t))= 
\frac{\dot T(t)^2}{(T(t)^2-1)^2}+\frac{1}{T(t)^2-1}\E{\dot r(t) + T(t)\dot
  N(t)}{\dot r(t) + T(t)\dot N(t)}.
\end{equation}
By using inequality ~(\ref{fund-ineq}), we see that
$\E{\dot r(t)}{\dot N(t)}\geq 0$, then we have
\[
   g_\lambda(\dot p(t), \dot p(t))\geq 
\frac{T(t)^2}{T(t)^2-1}\E{\dot N(t)}{\dot
   N(t)}\geq\E{\dot N(t)}{\dot N(t)} \ .
\]
\cvd
\begin{lem}\label{hyperbolic:compl2:lem}
Let $(q_n)$ be a sequence in $M_\lambda$.  Then it is a Cauchy
sequence if and only if either it converges to a point
$q_\infty\in\Uu_\lambda(>1)$ or $N(q_n)$ is a Cauchy sequence and
$T(q_n)$ converges to $+\infty$.
\end{lem}
\Dim
Consider the function of $M_\lambda$
\[
\Delta=\arctgh (1/T) \ .
\]
By using equation~(\ref{hyperbolic:compl:eq}) we can easily see that
$\Delta$ is $1$-Lipschitz function.  So both $N$ and $\Delta$ extend
to continuous functions of $\overline M_\lambda$ and if $(q_n)$ is a
Cauchy sequence in $M_\lambda$ then $N(q_n)$ and $\arctgh(1/T(q_n))$
are Cauchy sequences.  In particular either $T(q_n)$ converges to
$a>1$ or to $+\infty$. In the first case it is not difficult to show
that $q_n$ runs in a compact set of $\Uu_\lambda(>1)$ so it converges
to a point in $\Uu_\lambda(>1)$.\par This proves the only if part.
Now suppose $(q_n)$ to be a sequence such that $N(q_n)\rightarrow
x_\infty$ and $T(q_n)\rightarrow +\infty$. We have to show that
$(q_n)$ is a Cauchy sequence.\par Let us introduce the following
notation: for $p\in\Uu_\lambda$ and $a>0$ let $p_a$ denote the point
on the orbit through $p$ of the flow of $N$ such that $T(p_a)=a$, that
is
\[
   p_a=r(p)+aN(p).
\]
We need the following statement proved in~\cite{Bo} Prop. 7.1 .
\vspace{.2cm}\\ \emph{ Denote by $d_a$ the intrinsic metric of the
surface $\Uu_\lambda(a)$. Moreover for every point $p\in\Uu_\lambda$
and $a>0$ let $p_a=r(p)+aN(p)$.  Then for every compact subset $K$ of $\mh^2$
and for every $\eps>0$ there exists $M>0$ such that for  $p,q\in
N^{-1}(K)$ we have
\[
      \left|\frac{1}{a} d_a(p_a,q_a)-d_\mh(N(p),N(q))\right|<\eps
\]
for every $a>M$.
}
\vspace{.2cm}\\
If $\delta_a$ denote the distance of the surface $\Uu(a)$ with respect to
the metric $g_\lambda$ we see that 
\[
    \delta_a = \frac{1}{\sqrt{a^2-1}}d_a \ .
\]
Thus given a compact set $K$  of $\mh^2$ containing $N(q_n)$ 
and $\eps>0$ we can find $M$ such that
if $a>M$ then 
\[
    |\delta_a(p_a,q_a)-d_\mh(N(p),N(q))|<\eps \ .
\]
for every $p,q$ such that $N(p),N(q)\in K$.
Let us set $a_n=T(q_n)$ and fix $N$ such that $a_n>M$ if $n>N$.
Now for $n,m>N$ we have
\[
   \delta (q_n,q_m)\leq\delta(q_n,(q_m)_{a_n})+\delta(q_m,(q_m)_{a_n})\ .
\]
Now the first term of this sum is less than
$\eps+d_\mh(N(q_n),N(q_m))$ whereas the last term is less than the
length of the arc $c(t)=r(q_m)+(ta_n+(1-t)a_m)q_m$.  Since the length
of $c_m$ is less than $\arctgh(1/M)$ we get that $(q_n)$ is a Cauchy
sequence.  \cvd 
\smallskip

\emph{Proof of statements (1) (2) and of
Theorem~\ref{hyperbolic:compl:teo}:} From
Lemma~\ref{hyperbolic:compl:lem} we have that the map
$N:M_\lambda\rightarrow\mh^2$ extends to a map on $\overline
M_\lambda$.  Moreover Lemma~\ref{hyperbolic:compl2:lem} implies that
the map $N$ restricted to $\partial M_\lambda=\overline
M_\lambda-M_\lambda$ is injective.  We want to prove that it is an
isometry.  Take $x,y\in\mh^2$ and $p,q\in M_\lambda$ such that
$N(p)=x$ and $N(q)=y$.  Paths $(p_a)_{a>1}$ and $(q_a)_{a>1}$ converge
as $a\rightarrow+\infty$ to points $q_\infty$ and $p_\infty$ of
$\partial M_\lambda$.  Moreover we have that
\[
   \delta(q_\infty, p_\infty)=\lim_{a\rightarrow+\infty}\delta(q_a,p_a)\ .
\]
Now $\delta(q_a,p_a)<\delta_a(q_a,p_a)$ and we know that
$\delta_a(q_a, p_a)\rightarrow d_\mh(x,y)$.  So we deduce that
$\delta(q_\infty, p_\infty)\leq d_\mh (x,y)$.  On the other hand since
$N$ is $1$-Lipschitzian the other inequality holds.
\cvd
\smallskip

We are going now to prove statement (3) of the theorem.  In fact we
have
\begin{cor}
The function $\Delta$ is $\mathrm C^1$. Moreover the following formula
holds
\[
   \Delta(p)=\arctgh(1/T(p)).
\]
For every point $p\in M_\lambda$ the unique point realizing $\Delta$
on the boundary is $N(p)$ and the geodesic joining $p$ to $N(p)$ is
parametrized by the path
\[
    c:[T(p),+\infty)\ni t\mapsto r(p)+tN(p)\in M_\lambda.
\]
\end{cor}
\Dim
If $p(t)$ is a $\mathrm C^1$-path it
we have seen that
\[
   g_\lambda(\dot p(t),\dot p(t))\geq (\dot T(t))^2/(T^2-1)^2
\]
and the equality holds if and only if $\dot r(t)=0$ and $\dot N(t)=0$.
Thus we obtain $\Delta(p)\geq\arctgh(1/T(p))$. The path $c$ has
hyperbolic length equal to $\arctgh(1/T(p))$ so $\Delta(p)= \arctgh(1/T(p))$.
Moreover if $p(t)$ is a geodesic realizing the distance $\Delta$ we have that
$\dot r=0$ and $\dot N=0$ so $p$ is a parametrization of $c$.
\cvd
\smallskip
  
Finally we want to show that $\overline M_\lambda$ is a manifold with
boundary $\partial M_\lambda=\mh^2$ homeomorphic to $\mh^2\times
[0,+\infty)$.  Notice that it is sufficient to show that the for every
$\eps>0$ the set $\overline M_\lambda(\Delta\leq\eps)$ is homeomorphic
to $\mh^2\times [0,\eps]$.  Unfortunately the map
\[
  \overline M_\lambda(\Delta\leq\eps)\ni x\mapsto (N(x),
  \Delta(x))\in\mh^2\times[0,\eps]
\]
works only if $L_W$ is empty. Otherwise it is not injective.
Now the idea to avoid this problem is the following.
Take a point $z_0\in M_\lambda$ and consider the surface
\[ 
\mh(z_0)=\{x\in\fut(r(z_0))|\E{x-r(z_0)}{x-r(z_0)}= -T(z_0)^2\}.
\]
It is a spacelike surface of $\Uu_\lambda(>1)$ (in fact $\mh(z_0)$ is
contained in $\Uu_\lambda(>a)$ for every $a<T(z_0)$).  Denote by $v$
the Gauss map of the surface $\mh(z_0)$.  It sends the metric of
$\mh(z_0)$ to the hyperbolic metric multiplied by a factor $1/T(z_0)$.
Now we have an embedding
\[
  \varphi:\mh(z_0)\times [0,+\infty)\ni (p,t)\mapsto p+tv(p)\in\Uu_\lambda(>1)
\]
that parameterizes the future of $\mh(z_0)$.  Clearly if we cut the
future of $\mh(z_0)$ from $M_\lambda$ we obtain a manifold
homeomorphic to $\mr^2\times (0,+1]$.  Thus in order to prove that
$\overline M_\lambda$ is homeomorphic to $\mr^2\times[0,1]$ is
sufficient to prove the following proposition.
\begin{prop}\label{hyperbolic:compl:top:prop}
The map $\varphi$ extends to a map
\[
   \varphi:\mh(z_0)\times [0,+\infty]\mapsto\overline M_\lambda
\]
that is an embedding on a neighbourhood of $\partial M_\lambda=\mh^2$ in
$\overline M_\lambda$ such
that
\[
   \varphi(p,+\infty)=v(p) \ .
\]
\end{prop}
\Dim It is not hard to show that a fundamental family of neighbourhoods
of a point $v_0\in\mh^2=\partial M_\lambda$ in $\overline M_\lambda$
is given by
\[
    V(v_0;\eps,a)=\{x\in\Uu_\lambda| d_\mh(N(x),v_0)\leq\eps\textrm{ and
    }T(x)\geq a\}\cup \{v\in\mh^2|d_\mh(v,v_0)\leq\eps\} \ .
\]
Now we claim that for any compact $H\subset\Uu_\lambda$, compact
$K\subset\mh^2$, and $\eps,a>0$ there exists $M>0$ such that
\[
     p_0 + t v_0\in V(v_0;\eps,a)
\]
for every $p_0\in H$, $v_0\in K$ and $t\geq M$.\\ 
Before proving the claim let us show that it implies the extension of
$\varphi$ defined in the statement of this proposition is continuous.
Indeed if $p_n$ is a 
sequence in $\mh(z_0)$ converging to $p_\infty$ 
then $v_n=v(p_n)$ is a convergent sequence in $\mh^2$ with
limit $v_\infty=v(p_\infty)$.
Thus if $t_n$ is any divergent sequence of positive numbers by the claim it
follows that $p_n+t_nv_n\rightarrow v_\infty$ in $\overline M_\lambda$.\\

Now, let us prove the claim.  Let us set $p(t)=p_0+t v_0$ and denote by
$r(t)$, $N(t)$, $T(t)$ the retraction, the Gauss map and the
cosmological time computed at $p(t)$.  Notice that $T(t)>t+T(0)$ and
since $p_0$ runs in a compact set there exists $m$ that does not
depend on $p_0$ and $v_0$ such that $ T(t)>t+m$.\\ On the other hand
by deriving the identity
\[
   p(t)=r(t)+T(t)N(t)
\]
we obtain 
\begin{equation}\label{hyperbolic:compl:top:eq1}
   v_0=\dot p(t)=\dot r(t)+T(t)\dot N(t)+\dot T(t)N(t) \ .
\end{equation}
By taking the scalar product with $\dot N$ we obtain
\[
  \E{v_0}{\dot N}=\E{\dot r}{\dot N}+T\E{\dot N}{\dot N}>0\ .
\]
Since  $\ch(d_\mh(v_0, N(t))=-\E{v_0}{N(t)}$,
the function 
\[
  t\mapsto d_\mh(v_0, N(t))
\]
is decreasing.  Thus there exists a compact set $L\subset\mh^2$ such that
$N(t)\in L$ for every $p_0\in H$ , $v_0\in K$, $t>0$.  It follows that
there exists a compact set $S$ in $\mx_0$ such that $r(t)\in S$ for every
$t>0$, $p_0\in H$ and $v_0\in K$.  We can choose a point $q\in\mx_0$
such that $S\subset\fut(q)$.  Notice that
\[
   T(t)=\sqrt{-\E{p(t)-r(t)}{p(t)-r(t)}}\leq\sqrt{-\E{p(t)-q}{p(t)-q}} \ .
\]
By using this inequality it is easy to find a constant $M$ (that depends
only on $H$ and $K$) such that
\[
   T(t)\leq t + M.
\]
This inequality can be written in the following way:
\[
   \int_0^{t} (\dot T(s)-1)\mathrm ds \leq M \ .
\]
On the other hand by using identity~(\ref{hyperbolic:compl:top:eq1}) we have
\[
  \ch\d_\mh(v_0, N(t))=-\E{v_0}{N(t)}=\dot T(t)
\]
so $\dot T>1$. It follows that the measure of the set
\[
    I_\eps=\{s| \dot T(s)-1>\eps\}
\]
is less than $M/\eps$. Since $T$ is concave, $I_\eps$ is an
interval (if non-empty) of $[0,+\infty)$ with an endpoint at $0$.
Thus $I_\eps$ is contained in $[0, M/\eps]$.\\ Finally we have proved
that for $t>\max(M/\eps, a)$ we have $T(t)>a$ and
\[
   \ch\d_\mh(v_0, N(t))=-\E{v_0}{N(t)}=\dot T(t)\leq 1+\eps \ .
\]
Thus the claim is proved.\\

In order to conclude the proof we have to show that $\varphi$ is
proper and the image is a neighbourhood of $\mh^2=\partial M_\lambda$
in $\overline M_\lambda$.  For the last statement we will show that
for every $v_0$ and $\eps>0$ there exists $a>0$ such that
$V(v_0;\eps,a)\subset\fut(\mh(z_0))$.  In fact by using that $N$ is
proper on level surfaces we have that there exists a compact set $S$ such
that $r(V(v_0;\eps,a))$ is contained in $S$ for every $a>0$. In
particular it is easy to see that there exist constants $c,d$ such
that if we take $p\in V(v_0;\eps,a)$ we have
\[
  \E{p-r(z_0)}{p-r(z_0)}=-T(p)^2+cT(p)+d \ .
\]
We can choose $a_0$ sufficiently large such that if $T(p)>a_0$ then
$\E{p-r_0}{p-r_0}<-T(z_0)$. So $V(v_0;\eps,a)\subset\fut(\mh_{q_0})$
for $a>a_0$.\\

Finally we have to prove that if $q_n=\varphi(p_n, t_n)$ converges to
a point then $p_n$ is bounded in $\mh(z_0)$.\\ Let us set
$p_n(t)=p_n+tv(p_n)=r_n(t)+T_n(t)N_n(t)$.  Since $N_n(t_n)$ is compact
there exists a compact $S$ such that $r_n(t_n)\in S$.  Thus by arguing
as above we see that there exist $c,d>0$ such that
\[
    \E{q_n-r(z_0)}{q_n-r(z_0)}\leq-T(q_n)^2+cT(q_n)+d \ .
\]
Hence we can find $M>0$ such that
\[
    T(q_n)\leq t_n+M \ .
\]
On the other hand we have
\[
    T(q_n)-t_n>\int_0^{t_n}(-\E{N_n(t)}{v(p_n)}-1)\mathrm ds \ .
\]
Since $-\E{N_n(t)}{v(p_n)}-1=\dot T_n-1$ is a decreasing positive
function we have that for every $\eps>0$
\[
   0<-\E{N_n(t)}{v(p_n)}<1+\eps
\]
for $t_n\geq t\geq M/\eps$.
In particular, for $n$ sufficiently large
we have that $\E{N_n(t_n)}{v(p_n)}<2$. Thus $v(p_n)$ runs is a compact
set. Since $p_n=v(p_n)+r_0$  conclusion follows.
\cvd
\smallskip

The proof of Theorem \ref{hyperbolic:compl:teo} is complete.

\paragraph{Projective boundary of $M_\lambda$}
We have seen that $\overline M_\lambda$ is homeomorphic to
$\mr^2\times [0,+\infty)$. Now let us define $\hat M_\lambda=\overline
M_\lambda\cup\Uu_\lambda(1)$.  Clearly $\hat M_\lambda$ is
homeomorphic to $\mr^2\times[0,+\infty]$.  In this section we will
prove that the map $D_\lambda:\overline M_\lambda\rightarrow\mh^3$ can
be extended to a map
\[
D_\lambda:\hat M_\lambda\rightarrow\overline\mh^3
\]
in such a way that the restriction of $D_\lambda$ on $\Uu_\lambda(1)$
takes value on $S^2_\infty=\partial\mh^3$ and is a $\mathrm C^1$-developing
map for a projective structure on $\Uu_\lambda(1)$.\\ In fact for a
point $p\in\Uu_\lambda(1)$ we know that the image via $D_\lambda$
of the path $(p_t)_{t>1}$ is a geodesic ray in $\mh^2$ starting at 
$F_\lambda(N(p))$.  Thus we define $D_\lambda(p)$ to be
the limit point in $S^2_\infty$ of this geodesic ray.  It is not hard to
prove that such an extension is continuous.  Moreover we will prove
the following theorem. In the statement we use the lamination 
$\hat \lambda$ on the level surface $\Uu_\lambda(1)$ 
defined in Section \ref{flat}. 
\begin{teo}\label{hyperbolic:proj:teo}
The map $D_\lambda:\Uu_\lambda(1)\rightarrow S^2_\infty$ is a local $\mathrm
C^1$-conformal map. In particular it is a developing map for a
projective structure on $\Uu(1)$.\\ The {\rm canonical stratification}
associated to this projective structure coincides with the
stratification induced by the lamination $\hat\lambda$ and its
{\rm Thurston metric} coincides with the intrinsic spacelike surface
metric $k_\lambda$ on $\Uu_\lambda(1)$ .
\end{teo}
\begin{remark}\emph{
Let us remind the notion of {\it canonical stratification} and {\it
Thurston metric} mentioned in the above statement
(see~\cite{Ku, Ap} for details). In general, if $D:\tilde
M\rightarrow S^2_\infty$ is a developing map for a
projective structure of ``hyperbolic type'' on
$M$, let us define an {\it open round ball} in $\tilde M$ to be a
subset $\Delta$ of $\tilde M$ such that the restriction of $D$ on it
is a homeomorphism onto a round ball of $S^2_\infty$. One can see that a
round ball $\Delta$ in $\tilde M$ is maximal (with respect to the
inclusion) if and only if its closure $\overline\Delta$ is not compact
in $\tilde M$. Moreover, $D(\overline\Delta)$ is contained in
$\overline {D(\Delta)}$. Let $\Lambda=\overline{
D(\Delta)}-D(\overline\Delta)$ and denote by $\Delta'$ the pre-image
in $\Delta$ of the convex core of $\Lambda$ with respect to the
hyperbolic metric of $D(\Delta)$.  Now for every point $p\in\tilde M$
there exists a unique maximal ball $\Delta_p$ such that
$p\in\Delta_p'$. So the family $\{\Delta'_p\}$ furnishes a
stratification of $\tilde M$ by hyperbolic ideal convex sets.  We call
it the canonical stratification.  The Thurston metric is a $\mathrm
C^{1,1}$ metric $h$ on $\tilde M$ such that $h(p)$ coincides with the
hyperbolic metric of $\Delta_p$.  }\end{remark}

\Dim The first part of this theorem is proved just as
Theorem~\ref{hyperbolic:WR:teo}.  In fact an explicit computation
shows that $D_\lambda:\Uu(1)\rightarrow S^2$ is a $\mathrm
C^1$-conformal map if $\lambda$ is a weighted geodesic. Thus it
follows that $D_\lambda$ is a $\mathrm C^1$-conformal map if $\lambda$
is a simplicial lamination.  Then by using standard approximations we
can prove that $D_\lambda$ is a $\mathrm C^1$-conformal map. Indeed
$\Uu_\lambda(1)$ can be regarded as the graph of a $\mathrm
C^1$-function $\varphi_\lambda$ defined on the horizontal plane
$H=\{x_0=0\}$.  Moreover if $\lambda_n\rightarrow\lambda$ on a compact
set $K$, then $\varphi_{\lambda_n}$ converges to $\varphi_\lambda$ on
$H(K)=\{x|N(\varphi(x),x)\in K\}$ in $\mathrm C^1$-topology.  Thus by
using parametrizations of $\Uu_{\lambda_n}(1)$ given by
\[
\sigma_{\lambda_n}(x)=(\varphi_{\lambda_n}(x),x)
\]
we obtain maps
\[
   d_{\lambda_n}:H(K) \rightarrow S^2_\infty \ .
\]
The same argument used in Theorem~\ref{hyperbolic:WR:teo} shows that
$d_{\lambda_n}$ converges to $d_\lambda$ on $H(K)$ in $\mathrm
C^1$-topology.  Finally if $k_n$ is the pull-back of $k_{\lambda_n}$
on $H(K)$, we have that $k_n$ converges on $H(K)$ to the pull-back of
$k_\lambda$.  Since $d_{\lambda_n}:(H,g_n)\rightarrow S^2_\infty$ is a
conformal map by taking the limit we obtain that $d_\lambda$ is a
conformal $\mathrm C^1$-map on $H(K)$.\\

The proof of the second part of the statement is more difficult.
Consider the round disk $\md_0$ in $S^2$ such that
$\partial\md_0=\partial\mh^2$ that is the infinite boundary of the
right half-space bounded by $\mh^2$ in $\mh^3$.  Notice that the
retraction $\md_0\rightarrow\mh^2$ is a conformal map (an isometry if
we endow $\md_0$ with its hyperbolic metric). We denote by
$\sigma:\mh^2\rightarrow\md_0$ the inverse map.  With this
notation the map $D_\lambda:\Uu(1)\rightarrow S^2_\infty$ can be expressed in
the following way:
\[
    D_\lambda(p)=\hat B(p_0,p)\sigma(N(p)) \ .
\]
Now for every point $p\in\Uu(1)$ let us consider the round circle
$\md_p=\hat B(p_0,p)(\md_0)$ and define $\Delta_p$ to be the connected
component of $D_\lambda^{-1}(\md_p)$ containing $p$.\\ We claim that if
we consider on $\md_p$ the hyperbolic metric $g_{\md_p}$ then for
$q\in\Delta_p$ we have
\[
    D_\lambda^*(g_{\md_p})(q)=\eta k_\lambda(q)
\]
where $k_\lambda$ is as usual the intrinsic metric of 
$\Uu_\lambda$ and $\eta$ is a positive number such that
\begin{equation}\label{hyperbolic:proj:eq}
    \log\eta>\int_{[N(p),N(q)]}\delta(t)\d\mu(t)+a(p,q) \ .
\end{equation}
where $\delta(t)$ is the distance of $N(q)$ from the piece of $\lambda$
containing $N(t)$ and $a(p,q)$ is defined in the following way:
$a(p,q)=0$ if $N(p)\neq N(q)$ and otherwise
$a(p,q)=\E{r(p)-r(q)}{r(p)-r(q)}^{1/2}$.

Before proving the claim, let us show how we can conclude the proof.
We have that $D_\lambda:\Delta_p\rightarrow\md_p$ strengthens the
lengths.  Thus a classical argument shows that it is a homeomorphism.
Since $\Delta_p$ is not compact in $\Uu_\lambda(1)$ it is a maximal
round ball.  If $F_p$ is the stratum of $\hat\lambda$ through $p$ then
$D_\lambda|_{F_p}=B(p_0,p)\circ N$. Thus the image of
$F_p$ in $\md_p$ is an ideal convex set. Moreover if $\Delta'_p$ is
the stratum corresponding to $\Delta_p$ the same argument shows that
$F_p\subset\Delta_p'$ and in particular $\Delta'_p$ is the stratum
through $p$.\\ Now
\[
(D_\lambda)_{*,p}:T_p\Uu_\lambda(1)\rightarrow T_{D_\lambda(p)}\md_p
\] 
is a conformal map, moreover its restriction on $T_p F_p$ is an
isometry (with respect to the hyperbolic metric of $\md_p$). Thus it
is an isometry and this shows that $k_\lambda$ coincides with the
Thurston metric.\\ Finally we have to show that $F_p=\Delta'_p$. If
$q\notin F_p$, formula~(\ref{hyperbolic:proj:eq}) implies that
\[
(D_\lambda)_{*,q}:T_q\Uu_\lambda(1)\rightarrow T_{D_\lambda(q)}\md_p
\] 
is not an isometry. Thus $\Delta_p$ is different from
$\Delta_q$ and  $q\notin\Delta'_p$.\\

\begin{figure}
\begin{center}
\input{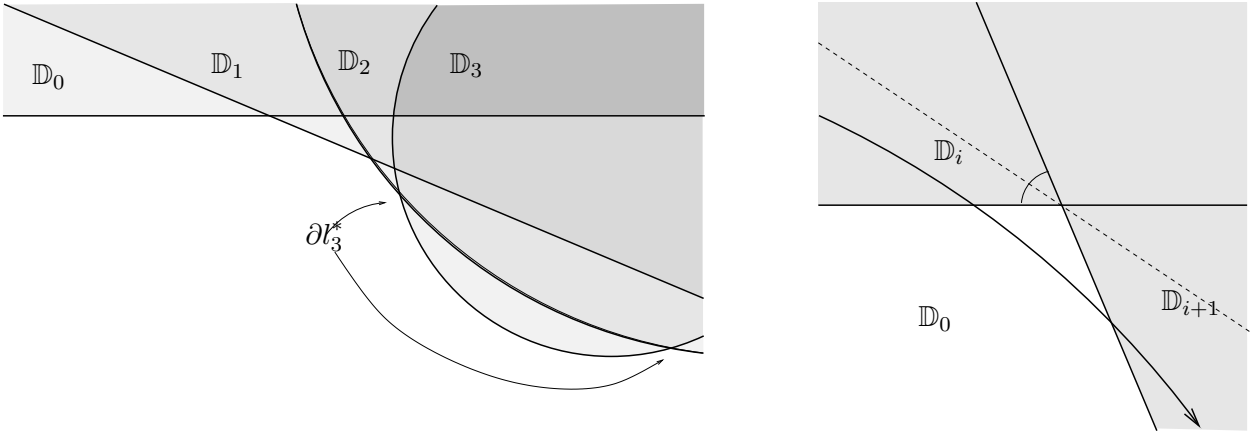}
\caption{{\small On the left it is shown how disks $\md_i$ intersect
    each other. On the right picture shows that if
    (\ref{hyperbolic:proj:eq3}) is not verified then
    (\ref{hyperbolic:proj:eq2}) does not hold}.}\label{proj:fig}
\end{center}
\end{figure}
Now we have to prove the claim.  The following lemma gives the
estimate we need.
\begin{lem}\label{hyperbolic:proj:lem}
Let $l$ be an oriented geodesic of $\md_0$.  Denote by $R\in
 PSL(2,\mc)$ the rotation around the corresponding geodesic of $\mh^3$
 with angle $\alpha$.  Take a point $p\in\md_0$ in the right
 half-space bounded by $l$ and suppose $R(p)\in\md_0$. Then if $g$ is
 the hyperbolic metric on $\md_0$ we have
\[
    R^*(g)(p)= \eta(p) g(p)
\]
where $\eta(p)=(\cos\alpha-\sh(d)\sin\alpha)^{-1}$ where $d$
is the distance from $p$ and $l$.
\end{lem}
\Dim Up to isometries we can identify $\md_0$ with the half-plane
$\{(x,y)|y>0\}$ in such a way that $l=\{x=0\}$ is oriented from $0$
towards $\infty$.  In this coordinates we have
\[
    R(x,y)=(\cos\alpha x+\sin\alpha y, -\sin\alpha x+\cos\alpha y)\ .
\]
Since $p$ is in the right half-plane bounded by $l$ then its
coordinates $(x,y)$ are both positive. Moreover as $R(p)\in\md_0$ we
have that $y/x>\tan\alpha$.  Now by an explicit computation we have
\[
   R^*(g)(p)=  \frac{1}{\cos\alpha y-\sin\alpha x}(\d x^2+\d y^2)
\]
then we see that $\eta(p)=(\cos\alpha-u\sin\alpha)^{-1}$ where
$u=x/y$.  On the other hand by classical hyperbolic formulas we have
that $x/y=\sh d$ where $d$ is the distance of $p$ from $l$.  \cvd

Now let us consider both the cases when $\lambda$ is a simplicial
lamination and $N(p)=N(q)$.  Up to post-composition with an element of
$PSL(2,\mc)$ we can suppose that the point $p$ is the base point so
$\md_p$ is $\md_0$.  Take $q\in\Delta_p$ and consider a path $c$ in
$\Delta_p$ containing $p$ and $q$. It is not hard to see that the
intersection of every stratum of $\hat\lambda$ with $\Delta_p$ is
convex. Thus we can suppose that $c$ intersects every leaf at most
once.\\ Denote by $l_0,l_1,\ldots,l_n$ the leaves intersecting $N(c)$
and let $a_0,\ldots a_n$ be the respective weights with the following
modifications.  If $q$ lies in a Euclidean band denote by $a_n$ the
distance from the component of the boundary of the band that meets
$c$.  In the same way if $p$ lies in a Euclidean band $a_1$ is the
distance of $p_1$ from the component of the boundary hitted by $c$.
Finally if $p$ and $q$ lie in the same Euclidean boundary (that is the
case when $N(p)=N(q)$) then $n=1$ and $a_1$ is the distance between
$r(p)$ and $r(q)$ (that is by definition $a(p,q)$).\par Let us set
$B_i=\exp(a_{1}X_{1})\circ\cdots\circ\exp(a_iX_i)$ where $X_i$ is the
standard generator of the rotation around $l_i$.  Notice that
$B_n=B(p,q)$.\par We want to prove that $q\in\md_i=B_i(\md_0)$.  In
fact we will prove that $\md_i\cap\md_0$ is a decreasing sequence of
sets (with respect to the inclusion).  By the hypothesis on $c$ we
have that $\md_i\cap\md_0\neq\varnothing$. Moreover if we denote by
$X^*_{i+1}$ the standard generator along the geodesic
$l_{i+1}^*=B_i(l_{i+1})$ then
\begin{equation}\label{hyperbolic:proj:eq2}
   \exp(tX^*_{i+1})\md_i\cap\md_0\neq\varnothing
\end{equation}
for $0<t<a_{i+1}$
(in fact there exists a point $q'\in c$ lying on the Euclidean band of
$\Uu(1)$ corresponding to
$l_{i+1}$ with distance from the left side equal to $t$ and $D_\lambda(q')$
lies in the intersection (\ref{hyperbolic:proj:eq2})).
Now by induction we can show that
\begin{equation}\label{hyperbolic:proj:eq3}
\left\{
\begin{array}{l}
  \md_0\cap\md_{i+1}\subset\md_0\cap\md_{i}\\
   \textrm{ the component of }\partial\md_i-\partial l^*_{i}\textrm{
   containing }l^*_{i+1}\textrm{ does not meet }\md_0\ .
\end{array}\right.
\end{equation}
In fact suppose $\md_0\cap\md_{i+1}$ to be not contained in
$\md_0\cap\md_i$.  Since $\md_{i+1}$ is obtained by the rotation along
$l_{i+1}^*$ whose end-points are outside $\md_0$ it is easy to see
that there should exist $t_0<a_{i+1}$ such that
$\exp(t_0X^*_{i+1})\md_i$ does not intersect $\md_0$ (see
Fig.~\ref{proj:fig}).\\

Let $g_i$ denote the hyperbolic metric on $\md_i$. We have that
$D_\lambda^*(g_n)$ is the intrinsic metric on $\Uu(1)$. Moreover we
have that
\[
 g_{i}(D_\lambda(q))=\eta_i g_{i+1}(D_\lambda(q))
\]
with $\eta_i^{-1}=\cos a_i-u_i\sin a_i$ where $u_i=\sh d_i$ where $d_i$ is the
distance of $N(q)$ from $l_i$.
Since $\eta=\prod_{i=0}^{n-1}\eta_i$ we obtain
\[
   -\log\eta=\sum\log(\cos a_i)+\sum\log(1-u_i\tan a_i) \ .
\]
Now $\log\cos(a_i)<-a_i^2/2$ and $\log(1-u_i\tan a_i)>d_ia_i$ so we get
\[
   \log\eta>\sum d_ia_i + a(p,q) \ .
\]
Now if the lamination is not simplicial and $N(p)\neq N(q)$ by using a
sequence of standard approximations we obtain the result.  \cvd.
\begin{cor}
Every level surface $\Uu_\lambda(a)$ is conformally flat. So it carries a
natural complex structure.
\end{cor}
\Dim
The map
\[
   \Uu_\lambda\ni x\mapsto tx\in\Uu_{t\lambda}
\]
rescales the metric by a factor $t^2$. Moreover, it takes
$\Uu_\lambda(1/t)$ onto $\Uu_{t\lambda}(1)$.  \cvd
\subsection{$\Gamma$-invariant constructions}
Assume that the lamination $\lambda$ is invariant for the action of a
discrete group $\Gamma$, that is $\lambda$ is the lifting of a
measured geodesic lamination on the hyperbolic surface $F=\mh^2/\Gamma$.\\ The
following lemma, proved in~\cite{Ep-M} determines the behaviour of the
cocycle $B_\lambda$ under the action of the group $\Gamma$.
\begin{lem}\label{hyperbolic:group:lem1}
Let $\lambda$ be a measured geodesic lamination on $\mh^2$ invariant by the
action of $\Gamma$. Then if $B_\lambda:\mh^2\times\mh^2\rightarrow PSL(2,\mc)$ 
is the cocycle associated to $\lambda$ we have
\[
   B_\lambda(\gamma x,\gamma y)=\gamma\circ B(x,y)\circ\gamma^{-1}
\]
for every $\gamma\in\Gamma$.
\end{lem}
\cvd
Now if we fix a base point $x_0\in\mh^2$ we can consider the bending map 
\[
   F_\lambda:\mh^2\rightarrow\mh^3\ .
\]
If we define
\[
   h_\lambda(\gamma)= B_\lambda(x_0,\gamma x_0)\circ\gamma\in PSL(2,\mc)
\]
Lemma~\ref{hyperbolic:group:lem1} implies that $h_\lambda:\Gamma\rightarrow
PSL(2,\mc)$ is a homomorphism.
Moreover by definition it follows that $F_\lambda$ is
$h_\lambda$-equivariant.\\

On the other hand we have seen that there exists a homomorphism
\[
  f_\lambda:\Gamma\rightarrow\ISO_0(\mx_0)
\]
such that $\Uu_\lambda$ is $f_\lambda$-invariant and the Gauss map is
$f_\lambda$-equivariant that is
\[
    N(f_\lambda(\gamma)(p))=\gamma(N(p)) \ .
\]
By using this fact it is easy to see that
\[
   \hat B_\lambda(f_\lambda(\gamma)p,f_\lambda(\gamma)q)=\gamma\hat
   B_\lambda(p,q)\gamma^{-1}.
\]
hence that
\[
    \hat D_\lambda(f_\lambda(\gamma)p)=h_\lambda(\gamma)(D_\lambda(p)) \ .
\]
In particular we have that the map $\hat D_\lambda$ is a developing
map for a hyperbolic structure on $M_\lambda/f_\lambda(\Gamma)$.  The
completion of such a structure is a manifold with boundary
homeomorphic to $F\times [0,+\infty)$. The boundary is
isometric to $F$.\\

The map $D_\lambda:\Uu_\lambda(1)\rightarrow S^2$ is
$h_\lambda$-equivariant so it is a developing map for a projective
structure on $\Uu_\lambda(1)/\Gamma$.\\ Notice that given a marking
$F\rightarrow\Uu_\lambda(1)/f_\lambda(\Gamma)$, by using the flow of
the gradient of the cosmological time, we obtain a marking
$F\rightarrow\Uu_\lambda(a)/f_\lambda(\Gamma)$. Thus we obtain a path
in the Teichm\"uller-like space of projective structures on $F$ and clearly
an underlying path of conformal structures in the Teich\"uller space
of $F$.

\paragraph{Cocompact $\Gamma$-invariant case}
If the group $\Gamma$ is cocompact, we can relate this construction
with the Thurston parametrization of projective structures on a base
compact surface $S$ of genus $g\geq 2$. In fact it is not hard to see
that the projective structure on $\Uu_\lambda(1)/f_\lambda(\Gamma)$ is
simply the structure associated to $(\Gamma,\lambda)$ in Thurston
parametrization (here we take $S=F$).  We have that the conformal
structure on $\Uu_\lambda(1)/\Gamma$ is the {\it grafting} of
$\mh^2/\Gamma$ along $\lambda$ (see \cite{McM,Sc}).  It follows
that the surface $\Uu_\lambda(a)/f_\lambda(\Gamma)$ corresponds to
$gr_{\lambda/a}(F)$; $a\mapsto[\Uu_\lambda(a)/f_\lambda(\Gamma)]$ is a
real analytic path in the Teichm\"uller space $\Tt_g$.
Such a path has an endpoint in $\Tt_g$
at $F$ as $a\rightarrow +\infty$ and an end-point in Thurston
boundary $\partial\Tt_g$ corresponding to the lamination $\lambda$ (or
equivantely to the dual tree $\Sigma$).


\section{Rescaling: flat towards de Sitter Lorentzian geometry}\label{dS}
\paragraph{On the de Sitter space}
Let us consider the $4$-dimensional Minkowski spacetime
$(\mr^4,\E{\cdot}{\cdot})$ and set
\[
  \hat\mx_1=\{v\in\mr^4|\E{v}{v}=1\} \ .
\]
It is not hard to show that $\hat\mx_1$ is a Lorentzian sub-manifold
of constant curvature $1$. Moreover the group $\OO(3,1)$ acts on it by
isometries.  This action is transitive and the stabilizer of a point
is $\OO(2,1)$. It follows that $\hat\mx_1$ is an isotropic Lorentzian
spacetime and $\OO(3,1)$ coincides with the full isometry group of
$\hat\mx_1$.  Notice that $\hat\mx_1$ is orientable and
time-orientable. In particular $\SOO^+(3,1)$ is the group of
time-orientation and orientation preserving isometries whereas
$\SOO(3,1)$ (resp. $\OO^+(3,1)$) is the group of orientation
(resp. time-orientation) preserving isometries.\par The projection of
$\hat\mx_1$ into the projective space $\mathbb P^3$ is a local
embedding onto an open set that is the exterior of the Klein model of
$\mh^3$ in $\mathbb P^3$ (that is a regular neighbourhood of $\mathbb
P^2$ in $\mathbb P^3$).  We denote this set by $\mx_1$.  Now the
projection $\pi:\hat\mx_1\rightarrow\mx_1$ is a $2$-fold covering and
the automorphism group is $\{\pm Id\}$. Thus the metric on $\hat\mx_1$
can be pushed forward to $\mx_1$. In what follows we consider always
$\mx_1$ endowed with such a metric and we call it the {\it Klein model
of de Sitter spacetimes}. Notice that it is an oriented spacetime
(indeed it carries the orientation induced by $\mathbb P^3$) but it is
not time-oriented (automorphisms of the covering
$\hat\mx_1\rightarrow\mx_1$ are not time-orientation preserving).\par
Since the automorphism group $\{\pm Id\}$ is the center of the
isometry group of $\OO(3,1)$ it is not hard to see that $\mx_1$ is an
isotropic Lorentz spacetime. The isometry group of $\mx_1$ is
$\OO(3,1)/\pm Id$. Thus the projection
$\OO^+(3,1)\rightarrow\ISO(\mx_1)$ is an isomorphism.\par Since $\mx_1$
is isotropic, a Lorentzian metric on a manifold $M$ with constant
curvature $1$ is equivalent to a $(\mx_1, \ISO(\mx_1))$-structure.
Thus for every de Sitter spacetime $M$ we have a developing map
$D:\tilde M\rightarrow\mx_1$ and an holonomy representation
$h:\pi_1(M)\rightarrow\OO^+(3,1)$, which are compatible in the sense
stated in the Introduction. \\

We give another description of $\mx_1$. Given a point $v\in\mx_1$,
notice that $v^\perp$ cuts $\mh^3$ along a totally geodesic plane
$P_+(v)$.  In fact we can consider on $P_+(v)$ the orientation induced
by the half-space $U(v)=\{x\in\mh^3|\E{x}{v}\leq 0\}$.  In this way
$\hat\mx_1$ parameterizes the oriented totally geodesic planes of
$\mh^3$.  If we consider the involution given by changing the orientation 
on the set of the oriented totally geodesic planes of $\mh^3$, then  the
corresponding involution on $\hat\mx_1$ is simply $v\mapsto-v$. In
particular, $\mx_1$ parameterizes the set of (un-oriented) hyperbolic
planes of $\mh^3$.  For $v\in\mx_1$ we denote by $P(v)$ the plane
corresponding to $v$.  For $\gamma\in\OO^+(3,1)$ we have
\[
    P(\gamma x)=\gamma (P(x))
\]
(notice that $\OO^+(3,1)$ is the isometry group of both $\mx_1$ and
$\mh^3$).\\

Just as for $\mh^3$, the geodesics in $\hat\mx_1$ are obtained by
intersecting $\hat\mx_1$ with linear $2$-spaces. Thus geodesics in
$\mx_1$ are projective segments. It follows that, given two points
$p,q\in\mx_1$, there exists a unique geodesic joining them.\par A
geodesic line in $\hat\mx_1$ is spacelike (resp. null, timelike) if
and only if it is the intersection of $\hat\mx_1$ with a spacelike
(resp. null, timelike) plane.  For $x\in\hat\mx_1$ and a vector $v$
tangent to $\hat\mx_1$ at $x$ we have
\[
\begin{array}{ll}
\textrm{if }\E{v}{v}=1  &  \exp_x(tv)=\cos t x + \sin t v\\
\textrm{if }\E{v}{v}=0  &  \exp_x(tv)=x+tv\\
\textrm{if }\E{v}{v}=-1 &  \exp_x(tv)=\ch t x + \sh t v.
\end{array}
\]
This implies that a complete geodesic line in $\mx_1$ is
spacelike (resp. null, timelike) if it is a complete projective line
contained in $\mx_1$ (resp. it a projective line tangent to $\mh^3$,
it is a projective segment with both the end-points in
$\partial\mh^3$).  Spacelike geodesics have finite length equal to
$\pi$. Timelike geodesics have infinite Lorentz length.\\

Take a point $x\in\mh^3$ and a unit vector $v\in T_x\mh^3=x^\perp$.
Clearly we have $v\in\hat\mx_1$ and $x\in T_v\hat\mx_1$. Notice that the
projective line joining $[x]$ and $[v]$ in $\mathbb P^3$ intersects
both $\mh^3$ and $\mx_1$ in complete geodesic lines $c$ and $c^*$.
They are parametrized in the following way
\[
\begin{array}{l}
c(t)=[\ch t x +\sh t v] \ , \\
c^*(t)=[\ch t v + \sh t x] \ . 
\end{array}
\]
We say that $c^*$ is the geodesic dual to $c$. They have the same
end-points on $S^2_\infty$ that are $[x+v]$ and $[x-v]$.  Moreover
if $c'$ is the geodesic \emph{ray} starting from $x$ with speed $v$
the dual geodesic ray $(c')^*$ is the geodesic \emph{ray} on $c^*$
starting at $v$ with the same limit point on $S^2_\infty$ as
$c'$.\\

\paragraph{Canonical Rescaling}
Let us use notations introduced in the previous section.  Given a
measured geodesic lamination $\lambda$ on $\mh^2$ we shall construct a
map
\[
      D^*:\Uu_\lambda(<1)\rightarrow\mx_1
\]
that, in a sense, is the map dual to the map $D$ constructed in the
previous section.  We shall prove that such a map is $\mathrm C^1$ and
the pull-back of the de Sitter metric is a rescaling of the flat
metric of $\Uu_\lambda$.  The idea to construct $D^*$ is very
simple. In fact if $s$ is a geodesic integral line of the gradient of
cosmological time we know that $s_{>1}=s\cap\Uu_\lambda(>1)$ is taken
by $D$ onto a geodesic ray of $\mh^3$. We define $D^*$ on $s_{<1}$ in
such a way that it parameterizes the dual geodesic ray in $\mx_1$.\\
Let us be more precise. Consider the standard inclusion
$\mh^2\subset\mh^3$. Since $\mh^2$ is oriented there is a well-defined
dual point $v_0\in\hat\mx_1$ (that is the positive vector of the
normal bundle).\\ Now let us take the base point $x_0\in\mh^2$ for the
bending map and a corresponding point $p_0\in\Uu_\lambda(1)$. For
$p\in\Uu_\lambda$ let us define
\[
\begin{array}{l}
  v(p)=\hat B_\lambda(p_0,p)v_0\in\hat\mx_1\\  
  x(p)=\hat B_\lambda(p_0,p)N(p)=F_\lambda(N(p))\ .
\end{array}
\]
Thus let us set
\[
    D^*(p)=[\ch\tau(p) v(p) +\sh\tau(p) x(p)]
\]
where we have put $\tau(p)=\arctgh T(p)$.
\begin{teo}\label{DS:rescaling:teo}
The map
\[
   D^*:\Uu_\lambda(<1)\rightarrow\mx_1
\]
is $\mathrm C^1$-local diffeomorphism. The pull-back of the metric of
$\mx_1$ is the rescaling of the metric of $\Uu_\lambda(<1)$ along the
gradient of $T$ with rescaling functions
\[
\begin{array}{ll}
   \alpha = \frac{1}{(1-T^2)^2}\,\qquad & \beta =\frac{1}{1-T^2}\ .
\end{array}
\]
\end{teo}
\Dim The proof of this theorem is quite similar to the proof of
Theorem~\ref{hyperbolic:WR:teo}.  In fact by an explicit computation
we get the result in the case when $\lambda$ is a weighted geodesic.
Thus theorem holds when $\lambda$ is a simplicial lamination. Moreover
by proving the analogous of Lemma~\ref{hyperbolic:WR:lem} and using
standard approximations we obtain the proof of the general case.\\

We shall explicit describe how to make the computation when the
lamination is a weighted geodesic. Then the same arguments of
Theorem~\ref{hyperbolic:WR:lem} work in the same way and we omit
details.\\
\begin{figure}
\begin{center}
\input{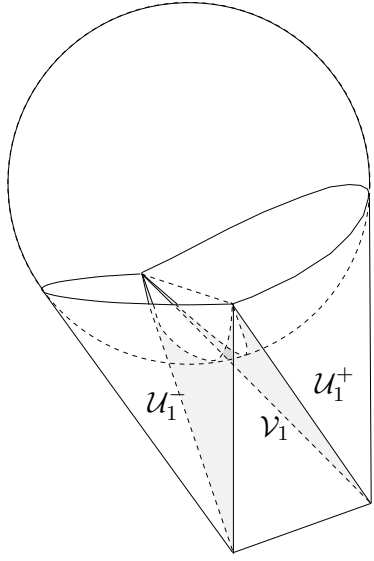}
\caption{{\small The domain $\Uu_1$ and its decomposition.}}
\end{center}
\end{figure}
Thus let us fix a weighted geodesic $(l_0, \alpha_0)$ with $\alpha_0<\pi$. We
have seen that the image of $D_0$ is the exterior $\Ee_0$ of the
convex set bounded by the bent surface $P=P_-\cup\exp (\alpha_0 X)P_+$.
Now it is not hard to show that the image of $D_0^*$ is the set of
points in $\mx_1$ whose dual plane are contained in $\Ee_0$. In fact
for a point $p\in\Uu_\lambda(<1)$ the plane dual to $\hat B(p_0,p)v_0$
separes the exterior of $\Ee_0$ from the plane dual to $D_0^*(p)$. On
the other hand if the plane dual $P(v)$ to $v\in\mx_0$ is entirely
contained in $\Ee_0$ then there exists a unique point $x'$ on it and
$x$ on $\partial\Ee_0$ such that the distance between $P(v)$ and
$\partial\Ee_0$ is equal to the distance between $x$ and $x'$. If $c$
is the geodesic ray starting from $x$ towards $x'$ we have that this
geodesic is orthogonal to $P(v)$ at $x'$ and the normal plane $Q(v)$
to $c$ at $x$ is a support plane for $\mh^3-\Ee_0$. Thus we have that
$c$ is the image through $D$ of an integral line of the gradient of
CT. The dual geodesic ray $c^*$ is contained in the image of
$D^*$. Since $P(v)$ is orthogonal to $c$ at $x'$ we have $v\in
c^*$.\par Denote by $\Uu_{1}$ the image of $D^*_0$. For every point
$v\in\Uu_{1}$ let us set $u(v)$ the dual point of the plane
$Q(v)$. Clearly $u(v)$ is a point on the boundary of $\Uu_{1}$ and
there exists a unique geodesic timelike segment between $v$ and $u(v)$
contained in $\Uu_1$.  Notice that the image of $u$ is the set of
points dual to planes obtained by rotating $\mh^2$ along $l_0$ by
angle less than $\alpha_0$. In particular it is a geodesic spacelike
segment of $\mx_1$.\par Then let us set $\tau(v)$ the proper time of
such a segment.  By making computation it is not hard to see that
$\tau$ is a $\mathrm C^1$-submersion taking values on the whole
$(0,+\infty)$.  The integral line of the gradient of $\tau$ through a
point $v$ is the geodesic ray between $v$ and $u(v)$.\par Given a
point $v\in\Uu_1$ let $n(v)$ be the intersection of the projective
line through $v$ and $u(v)$ with $\partial\Ee_0$ (by the above
discussion it follows that such a point is well-defined).  We denote
by $\Uu_1^\pm(a)=n^{-1}(\mathrm{int}(P^\pm))\cap\Uu_1(\tau=a)$ and
$\Vv_1=n^{-1}(l_0)$.  Notice that we have
\[  
    \begin{array}{l}
    \Uu_1^-(a)=u^{-1}(v_0)\cap\Uu_1(\tau=a) \ ,\\
    \Vv_1(a)=u^{-1}([v_0, \exp (\alpha_0X_0) v_0])\cap\Uu_1(\tau=a) \ ,\\
    \Uu_1^+(a)=u^{-1}(\exp(\alpha_0X_0)v_0)\cap\Uu_1(\tau=a) \ .
    \end{array}
\]
It is not hard to show that the map $n$ restricted to $\Uu_1^{\pm}(a)$
is a dilatation of a factor $(\sh a)^{-1}$ whereas the map
\[
   \Vv_1(a)\ni v\mapsto(u(v),n(v))\in[v_0, \exp (\alpha_0X_0) v_0]\times l_0
\]
sends the metric on $\Vv_1(a)$ onto the metric
\[
     (\ch a)^2\d u + (\sh a)^2\d t
\]
where $u$ and $t$ are the natural parameters on $[v_0, \exp (\alpha_0X_0)
v_0]$ and $l_0$.\par Now denote by $\delta_a$ the intrinsic distance
on the surface $\Uu_1(a)$. For every point $v\in\Uu_1(a)$ denote by
$\pi(v)$ the point on $l_a=\partial\Uu_1^-(a)$ that realizes the
distance of $l_a$ from $v$.  Let us fix a point $p_0\in l_0$ and
denote by $w_a$ the point on $l_a$ such that $n(w_a)=p_0$.  Then
consider the functions
\[
    \begin{array}{l}
     U(v)=\delta_{\tau(v)}(v,l_{\tau(v))}/\sh(\tau(v))\\
     Z(v)=\delta_{\tau(v)}(\pi(v), w_{\tau(v)})/\sh(\tau(v)) \ .
     \end{array}
\]
We have that $(\tau, U, Z)$ are $\mathrm C^1$-coordinates on
$\Uu_1$. Moreover by using coordinates $(T,u,\zeta)$ on $\Uu_0(<1)$
and $(\tau, U, Z)$ on $\Uu_1$ we have
\[  
    \begin{array}{l}
    \tau(D_0^*(T,u,\zeta))=\arctgh(T) \ ,\\
    U(D_0^*(T,u,\zeta))=u \ ,\\
    Z(D_0^*(T,u,\zeta))=\zeta \ .
    \end{array}
\]
Then by making computations in these coordinates it is not hard to
obtain that the pull-back of the de Sitter metric on $\Uu_0(<1)$ is
the rescaling of the flat metric along the gradient of $T$ with
rescaling functions
\[
   \begin{array}{lr}
    \alpha=\frac{1}{(1-T^2)^2}\ ,\qquad  &  \beta=\frac{1}{(1-T^2)} \ .
    \end{array}
\]
\cvd
\begin{cor}
The map $D^*$ extends to a continuous map
\[
   \Uu_\lambda(\leq 1)\cup\Sigma\rightarrow \mx_1\cup S^2 \ .
\]
Moreover $D^*$ restricted to $\Uu_\lambda(1)$ coincides with $D$.
\end{cor}
The extension of $D^*$ on $\Uu_\lambda(1)$ follows by construction. On
the other hand, we see that the cocycle $\hat B$ is induced by a
cocycle
\[
    \overline B:\Sigma\times\Sigma\rightarrow\SOO^+(3,1)
\]
and $D^*$ can be extended on $\Sigma$ by put
\[
  D^*(r)= \overline B(r(p_0),r)v_0 \ .
\]
\cvd
\begin{remark}\emph{
Notice that the above construction allows to identify $\Sigma$ with
the space of maximal round balls of $\Uu_\lambda(1)$.  }
\end{remark}

In what follows, we denote by $\Uu^1_\lambda$ the domain $\Uu_\lambda(<1)$
endowed with de Sitter metric induced by $D^*$.
\begin{prop}\label{desitter:ct:prop}
The cosmological time of $\Uu^1_\lambda$ is
\[
   \tau=\arctgh(T).
\]
Every level surface $\Uu^1_\lambda(\tau=a)$ is a Cauchy surface (so
$\Uu^1_\lambda$ is globally hyperbolic).
\end{prop}
\Dim Let $\gamma:[0,a]\rightarrow\Uu^1_\lambda$ denote a timelike-path
with future end-point $p$ parametrized in Lorentzian arc-length.  Now
as path in $\Uu_\lambda$ we have a decomposition of $\gamma$
\[
   \gamma(t)=r(t)+T(t)N(t) \ .
\]
By computing derivatives we obtain
\[
  \dot\gamma =\dot r+ T\dot N + \dot T N
\]
so the square of de Sitter norm is
\begin{equation}\label{desitter:ct:eq}
   -1=-\frac{\dot T^2}{(1-T^2)^2}+ \frac{|\dot r + T\dot N|^2}{1-T^2} 
\end{equation}
where $|\cdot|$ is the Lorentzian flat norm.
It follows that
\[
   1< \frac{\dot T}{(1-T^2)}
\]
and integrating we obtain
\[
   \arctgh T(p)-\arctgh T(0)> a
\]
that is the de Sitter proper time of $\gamma$ is less than $\arctgh
T(p)$.  On the other hand the path $\gamma(t)=r(p)+tN(p)$ for
$t\in[0,T(p)]$ has proper time $\arctgh T(p)$ so we obtain that the
cosmological time of $\Uu_\lambda^1$ is
\[
   \tau=\arctgh T \ .
\]

Now let $\gamma:(a,b)\rightarrow\Uu^1(\lambda)$ be an inextendable
timelike-curve parametrized in Lorentz arc-length such that
$\gamma(0)=p$. We want to show that the range of $T(t)=T(\gamma(t))$
is $(0,1)$.\par Suppose $\beta=\sup T(t)<1$. Since $T(t)$ is
increasing then $\beta=\lim_{t\rightarrow b}T(t)$.  Then the path
\[
   c(t)=r(t)+N(t)
\]
should be inextendable (otherwise we could extend $\gamma$ in
$\Uu^1_\lambda$).
Now we have
\[
   \dot c =\dot r + \dot N \ .
\]
For $t>0$ we have $T(t)>T(p)=T_0$ so
\[
  T_0 |\dot c| <|\dot r+T\dot N|
\]
Multiplying by the rescaling horizontal factor we have
\[
     \frac{T_0}{\sqrt{1-T^2}}|\dot c|\leq\frac{|\dot r+T\dot
     N|}{\sqrt{1-T^2}}.
\]
Since $T(t)<\beta<1$ it results
\[
 \frac{T_0}{\sqrt{1-\beta^2}}
  |\dot c|\leq\frac{|\dot r+T\dot
     N|}{\sqrt{1-T^2}}.
\]
By looking at equation~(\ref{desitter:ct:eq}) we deduce
\[
 \frac{T_0}{\sqrt{1-\beta^2}}|\dot c|\leq  \frac{\dot T}{1-T^2} \ .
\]
Thus the length of $c$ is bounded. On the other hand since
$\Uu_\lambda(1)$ is complete it follows that $c$ is extendable.  Thus
we have proved that $\sup T(t)=1$.  The same computation applied to
$\gamma(a,0)$ shows that $\inf T(t)=0$.  \cvd

\paragraph{$\Gamma$-equivariant constructions}
Suppose $\lambda$ to be a measured lamination of $\mh^2$ invariant by
the action of a group $\Gamma$ in $PSL(2,\mr)$.  We have seen that
there exists an affine deformation of $\Gamma$
\[
    f_\lambda:\Gamma\rightarrow\ISO_0(\mx_0)
\]
such that $\Uu_\lambda$ is $f_\lambda(\Gamma)$-invariant and the Gauss
map is $f_\lambda$-equivariant.  Moreover in the previous section we
have constructed a representation
\[
    h_\lambda:\Gamma\rightarrow PSL(2,\C)=\SOO^+(3,1)
\]
such that
\[
   D\circ f_\lambda(\gamma)=h_\lambda(\gamma)\circ D 
\qquad\textrm{ for }\gamma\in\Gamma.
\]
Now it is straightforward to see that the same holds changing $D$ by $D^*$.\\

\paragraph{Cocompact $\Gamma$-invariant case}
When $\Gamma$ is cocompact, we know that $f_\lambda(\Gamma)$ acts
properly on $\Uu_\lambda$ and the quotient $Y_\lambda$ is the unique
future-complete maximal globally hyperbolic spacetime diffeomorphic to
$\mh^2/\Gamma\times\mr$ with holonomy $f_\lambda$.  On the other hand
$f_\lambda(\Gamma)$ acts by isometries on the rescaled spacetime
$\Uu^1_\lambda$.  We want to relate this construction with Scannell's
classification of de Sitter spacetimes~\cite{Sc}.\par We have seen that
given a projective surface $F$ it is possible to construct a
hyperbolic manifold $H(F)$ (the $H$-hull). Now Scannell makes a
dual construction that produces a de Sitter spacetime $U(F)$
homeomorphic to $F\times(0,1)$ that satisfies the following conditions
\begin{enumerate}
\item
It is globally hyperbolic.
\item
Its developing map extends to a map
\[
   dev:\tilde F\times (0,1]\mapsto\mx_1\cup\partial\mh^3
\]
such that the restriction of $dev$ to the slice $F\times\{1\}$ is the
developing map for the projective structure on $F$.
\end{enumerate}
We call $U(F)$ the standard spacetime associated with $F$.  In fact by
construction it is not hard to see that
$Y^1_\lambda=\Uu^1_\lambda/f_\lambda(\Gamma)$ is the standard de
Sitter space-time associated to $U_\lambda(1)/f_\lambda(\Gamma)$ (that
we have seen carries a natural projective structure) .  So by Scannell's
classification theorem of de Sitter maximal hyperbolic spacetimes we
obtain the following theorem
\begin{teo}\label{desitter:class:teo}
The  correspondence
\[
   Y_\lambda\rightarrow Y^1_\lambda
\]
induces a bijection between flat maximal globally hyperbolic
spacetimes and de Sitter maximal globally hyperbolic spacetimes.
\end{teo}
\cvd
\begin{remark}\emph{
In general we see that the correspondence
\[
  \Uu_\lambda\mapsto\Uu^1_\lambda
\]
induces a bijection between regular domains with surjective Gauss map
and standard de Sitter spacetime corresponding to projective surface
with pleated locus isometric to $\mh^2$.  }\end{remark}


\section{Rescaling: Flat towards\\ 
Anti de Sitter Lorentzian geometry}
\label{AdS}
The AdS canonical rescaling runs parallel to the WR of Section
\ref{hyp}.  In fact, every spacelike plane $P$ is a copy of $\mh^2$
into the Anti de Sitter space $\mx_{-1}$. So the core of the
construction consists in a suitable {\it bending procedure} of $P$
along any given $\lambda \in \Mm\Ll(\mh^2)$.  However, in details
there are important differences.
\subsection{Anti de Sitter space}\label{generalAdS}
We recall some general features of the AdS local model that we will use
for the rescaling.  In particular, both spacetime and time orientations
will play a subtle role, so it is important to specify them carefully.

Let $\mathrm M_2(\mr)$ be the space of $2\times2$ matrices with real
coefficients endowed with the scalar product $\eta$ induced by the
quadratic form \[ q(A)=-\det A\ .\] The signature of $\eta$ is
$(2,2)$.

\noindent
The group\[ \SL{2}{R}=\{A|q(A)=-1\}\] is a Lorentzian sub-manifold of
$\mathrm M_2(\mr)$, that is the restriction of $\eta$ on it has
signature $(2,1)$.  Given $A,\ B\in\SL{2}{R}$, we have that
\[   q(AXB)=q(X) \qquad\textrm{ for }X\in\mathrm M_2(\mr)\]
Thus,  the left action of $\SL{2}{R}\times\SL{2}{R}$ on $\mathrm M_2{\R}$ 
given by
\[   (A,B)\cdot X=AXB^{-1}\]
preserves $\eta$. In particular, the restriction of $\eta$ on
$\SL{2}{R}$ is a bi-invariant Lorentzian metric, that actually
coincides with its Killing form. 
Notice that for $X,Y\in\sG\lG(2,\mr)$ we have the usual formula.
\[
    \tr XY=2\eta(X,Y) \ .
\]
 We denote by $\hat\mx_{-1}$ the
pair $(\SL{2}{R}, \eta)$.  Clearly $\hat\mx_{-1}$ is an orientable and
time-orientable spacetime.  Hence, the above action is a transitive
isometric action of $\SL{2}{R}\times\SL{2}{R}$ on $\hat\mx_{-1}$.

The stabilizer of $Id\in \hat\mx_{-1}$ is the diagonal group
$\Delta\cong\SL{2}{R}$.  It is not difficult to show that the
differential of isometries corresponding to elements in $\Delta$
produces a surjective representation
\[      \Delta\rightarrow\SOO^+(\sG\lG(2,\mr),\eta_{Id}) \ .\]
It follows that $\hat\mx_{-1}$ is an isotropic Lorentzian spacetime and  the
isometric action on $\hat\mx_{-1}$ induces a surjective representation
\[   \hat\Phi:\SL{2}{R}\times\SL{2}{R}\rightarrow\ISO_0(\hat\mx_{-1}).\]
Since $\ker\hat\Phi=(-Id, -Id)$, we obtain
\[   \ISO_0(\hat\mx_{-1})\cong\SL{2}{R}\times\SL{2}{R}/(-Id,-Id).\]
The center of $\ISO_0(\hat\mx_{-1})$ is generated
by $[Id,-Id]=[-Id,Id]$. Hence, $\eta$ induces
on the quotient\[PSL(2,\R)=\SL{2}{R}/\pm Id\]
an isotropic Lorentzian structure. 
We denote by $\mx_{-1}$ such a spacetime and call it the 
{\it Klein model} of Anti de Sitter spacetime. 
Notice that left and right translations are isometries and
the above remark implies that the induced representation
\[   \Phi:PSL(2,\R)\times PSL(2,\R)\rightarrow\ISO_0(\mx_{-1})\]
is an isomorphism.\\
\smallskip

\paragraph {The boundary of $\mx_{-1}$}
Consider the topological closure $\overline{PSL(2,\R)}$
of $PSL(2,\R)$ in $\mathbb P^3=\mathbb P(\mathrm M_2(\mr))$.  Its
boundary is the quotient of the set
\[    \{X\in\mathrm M_2(\mr)-\{0\}|q(X)=0\}\]
that is the set of rank $1$ matrices. In particular,
$\partial PSL(2,\R)$ is the image of the {\it Segre embedding}
\[   \mathbb P^1\times\mathbb P^1\ni([v],[w])
\mapsto[v\otimes w]\in\mathbb P^3.\] Thus $\partial PSL(2,\R)$ is a
torus in $\mathbb P^3$ and divides it in two solid tori. In
particular, $PSL(2,\R)$ topologically is a solid torus.\par
The action of $PSL(2,\R)\times PSL(2,\R)$ extends on the whole of
$\overline\mx_{-1}$.  Moreover, the action on $\mathbb
P^1\times\mathbb P^1$ induced by Segre embedding is simply
\[   (A,B)(v,w)=(Av,B^*w)\]
where we have set $B^*= (B^{-1})^T$ and considered the natural action of 
$PSL(2,\R)$ on $\mathbb P^1=\partial\mh^2$.
If $E$ denote the rotation by $\pi/2$ of $\mr^2$, it is not hard to show that
\[    EAE^{-1}=(A^{-1})^T\]
for $A\in PSL(2,\R)$.

It is convenient to consider the following modification of Segre
embedding
\[   S:\mathbb P^1\times\mathbb P^1\ni([v],[w])
\mapsto [v\otimes(Ew)]\in\mathbb P^3\] With respect to such a new
embedding, the action of $PSL(2,\R)\times PSL(2,\R)$ on
$\partial\mx_{-1}$ is simply
\[   (A,B)(x,y)=(Ax, By).\]
In what follows, we will consider the identification of the boundary
of $\mx_{-1}$ with $\mathbb P^1\times\mathbb P^1$ given by $S$. \\

The product structure on $\partial\mx_{-1}$ given by $S$ is preserved
by the isometries of $\mx_{-1}$. This allows us to define a {\it
conformal Lorentzian structure} (i.e. a {\it causal structure}) on
$\partial\mx_{-1}$.  More precisely, we can define two foliations on
$\partial\mx_{-1}$.  The left foliation is simply the image of the
foliation with leaves
\[    l_{[w]}=\{([x],[w])|[x]\in\mathbb P^1\}\]
and a leaf of the right foliation is the image of
\[    r_{[v]}=\{([v],[y])|[y]\in\mathbb P^1\}.\]
Notice that left and right leaves are projective lines in $\mathbb P^3$.
Exactly one left and one right leaves
pass through any given point. On the other hand, 
given right leaf and left leaf meet each other at one
point. Left translations preserve leaves of left foliation, whereas
right translations preserve leaves of right foliation.\par

If we orient $\mathbb P^1$ as boundary of $\mh^2$, we have that leaves
of right and left foliations are oriented.  Thus if we take a point
$p\in\partial\mx_{-1}$, the tangent space $T_p\partial\mx_{-1}$ is
divided by the tangent vector of the foliations in four quadrants.  By
using orientation of leaves we can enumerate quadrants in the usual way.
Thus we can consider the $1+1$ cone at $p$ given by
choosing the second and fourth quadrants. We make this choice because
in this way the causal structure on $\partial\mx_{-1}$ is the
``limit'' of the causal structure on $\mx_{-1}$ in the following
sense.\\ Suppose $A_n$ to be a sequence in $\mx_{-1}$ converging to
$A\in\partial\mx_{-1}$, and suppose $X_n\in T_{A_n}\mx_{-1}$ to be a
sequence of timelike vectors converging to $X\in T_A\partial\mx_{-1}$,
then $X$ is non-spacelike with respect to the causal structure of the
boundary.\\ Notice that oriented left (resp. right) leaves are
homologous non-trivial simple cycles on $\partial\mx_{-1}$, so they
determines non-trivial elements of $\coom1(\partial\mx_{-1})$ that we
denote by $c_L$  and $c_R$.

\paragraph{ Geodesic lines and planes}
Geodesics in $\mx_{-1}$ are obtained intersecting projective lines
with $\mx_{-1}$.

A geodesic is timelike if it is a projective line entirely contained
in $\mx_{-1}$; its Lorentzian length is $\pi$.  In this case it is a
non-trivial loop in $\mx_{-1}$ (a core).  Take $x\in\mx_{-1}$ and
$v\in T_x\mx_{-1}$ a unit timelike vector.  If $\hat x$ is a pre-image
of $x$ in $\hat\mx_{-1}$ and $\hat v\in T_{\hat x}\hat\mx_{-1}$ is a
pre-image of $v$, then we have
\[    \exp_x tv = [\cos t \hat x+\sin t \hat v] \ .\]

A geodesic is null if it is contained in a projective line tangent to
$\partial\mx_{-1}$.  Given $x\in\mx_{-1}$, and a null vector $v\in
T_x\mx_{-1}$, if we take $\hat x$ and $\hat v$ as above we have\[
\exp_x tv=[\hat x+t\hat v] \ .\]

Finally, a geodesic is spacelike if it is contained in a projective
line meeting $\partial\mx_{-1}$ at two points; its length is infinite.
Given $x\in\mx_{-1}$ and a unit spacelike vector $v$ at $x$, fixed
$\hat x$ and $\hat v$ as above, we have\[ \exp_x tv=[\ch t \hat x+\sh t
\hat v] \ .\]

Geodesics passing through the identity are $1$-parameter subgroups.
Elliptic subgroups correspond to timelike geodesics, parabolic
subgroups correspond to null geodesics and hyperbolic subgroups are
spacelike geodesics.\\

Totally geodesic planes are obtained intersecting projective planes
with $\mx_{-1}$.

If $W$ is a subspace of dimension $3$ of $\mathrm M_2(\mr)$ and the
restriction of $\eta$ on it has signature $(m_+,m_-)$, then the
projection $P$ of $W$ in $\mathbb P^3$ intersects $PSL(2,\R)$ if and
only if $m_- >0$. In this case the signature of $P\cap\mx_{-1}$ is
$(m_+,m_--1)$.  Since $\eta$ restricted to $P$ is a flat metric we
obtain that
\begin{enumerate}

\item
If $P\cap\mx_{-1}$ is a Riemannian plane, then it is isometric to
$\mh^2$.
\item If $P\cap\mx_{-1}$ is a Loretzian plane, then it is a Moebius
band carrying a de Sitter metric.
\item If $P\cap\mx_{-1}$ is a null plane, then $P$ is tangent to
$\partial\mx_{-1}$.
\end{enumerate}

\noindent
In particular, since every spacelike plane cuts every timelike
geodesics at one point, we obtain that spacelike planes are
compression disks of $\mx_{-1}$.  The boundary of a spacelike plane is
a spacelike curve in $\partial\mx_{-1}$ and it is homologous to
$c_L+c_R$.

\par Every Lorentzian plane is a Moebius band. 
Its boundary is homologous to $c_L- c_R$.

\par Every null plane is a pinched band. Its boundary is the union of one 
right and one left leaf.\\

\paragraph{Duality in $\mx_{-1}$}
The form $\eta$ induces a duality in $\mathbb P^3$ between points and
planes, and between projective lines. Since the isometries of
$\mx_{-1}$ are induced by linear maps of $\mathrm M_2(\mr)$ preserving
$\eta$, this duality is preserved by isometries of $\mx_{-1}$.

\par If we take a point in $\mx_{-1}$ its dual projective planes
defines a Riemannian plane in $\mx_{-1}$ and, conversely, Riemannian
planes are contained in projective planes dual to points in
$\mx_{-1}$.  Thus we have a bijective correspondence between points
and Riemannian planes.  Given a point $x\in\mx_{-1}$, we denote by
$P(x)$ its dual plane and, conversely, if $P$ is a Riemannian plane,
then $x(P)$ denote its dual point. If we take a point $x\in\mx_{-1}$
and a timelike geodesic $c$ starting at $x$ and parametrized in
Lorentzian arc-length we have that $c(\pi/2)\in P(x)$.  Moreover, this
intersection is orthogonal.  Conversely, given a point $y$ in $P(x)$,
there exists a unique timelike geodesic passing through $x$ and $y$
and such a geodesic is orthogonal to $P(x)$.  By using this
characterization, we can see that the plane $P(Id)$ consists of those
elliptic transformations of $\mh^2$ that are the rotation by $\pi$ at
their fixed points.  In this case an isometry between $P(Id)$ and
$\mh^2$ is simply obtained by associating to every $x\in P(Id)$ its
fixed point in $\mh^2$. Moreover, such a map
\[    I:P(Id)\rightarrow\mh^2\] 
is natural in the following sense.  The isometry group of $P(Id)$ is
the stabilizer of the identity, that we have seen to be the diagonal
group $\Delta\subset PSL(2,\R)\times PSL(2,\R)$.  Then we have\[
I\circ(\gamma,\gamma)=\gamma\circ I.\] The boundary of $P(Id)$ is the
diagonal subset of $\partial\mx_{-1}=\mathbb P^1\times\mathbb P^1$ that
is
\[   \partial P(Id)=\{(x,x)\in\partial\mx_{-1}|x\in\mathbb P^1\} \ .\]
The map $I$ extends to $\overline P(Id)$, by sending the
point $(x,x)\in\partial P(Id)$ on $x\in\mathbb P^1=\partial\mh^3$.\\

The dual point of a null plane $P$ is a point $x(P)$ on the
boundary $\partial\mx_{-1}$.  It is the intersection point of the left
and right leaves contained in the boundary of that plane.  Moreover,
the plane is foliated by null-geodesics tangent to $\partial\mx_{-1}$
at $x(P)$.  Conversely, every point in the boundary is dual to the
null-plane tangent to $\partial\mx_{-1}$ at $x$.

\par Finally the dual line of a spacelike line $l$ is a spacelike line
$l^*$.  Actually $l^*$ is the intersection of all $P(x)$ for $x\in l$.
$l^*$ can be obtained by taking the intersection of null planes dual
to the end-points of $l$. In particular, if $x_-$ and $x_+$ are the
end-points of $l$, then the end-points of $l^*$ are obtained by
intersecting the left leaf through $x_-$ with the right leaf through
$x_+$, and the right leaf through $x_-$ with the left leaf through
$x_+$, respectively.

\par There is a simple interpretation of the dual spacelike geodesic
for a hyperbolic $1$-parameter subgroup $l$. In this case $l^*$ is
contained in $P(Id)$ and is the inverse image through $I$ of the axis
fixed by $l$ in $\mh^2$.  Conversely, geodesics in $P(Id)$ correspond
to hyperbolic $1$-parameter subgroups.\\

\paragraph{Orientation and time-orientation of $\mx_{-1}$}
In order to define a time-orientation it is enough to define a time
orientation at $Id$.  This is equivalent to fix an orientation on the
elliptic $1$-parameters subgroups.  We know that such a subgroup
$\Gamma$ is the stabilizer of a point $p\in\mh^2$.  Then we stipulate
that an infinitesimal generator $X$ of $\Gamma$ is future directed if
it is a positive infinitesimal rotation around $p$.

\par A spacelike surface in a oriented and time-oriented spacetime is
oriented by means of the rule: {\it first the normal future-directed
vector field}.  So, we choose the orientation on $\mx_{-1}$ that
induces the orientation on $P(Id)$ that makes $I$ an
orientation-preserving isometry.

\par Clearly orientation and time-orientation on $\mx_{-1}$ induces
orientation and a time-orientation on the boundary. With respect to
these choice we have that left leaves are future-oriented whereas
right leaves are past oriented.

If $e_L$ and $e_R$ are no-where vanishing vector fields on
$\partial\mx_{-1}$ respectively tangent to left and right foliations
then the ordered pair $(e_R(x),e_L(x))$ for a positive basis of
$T_x\partial\mx_{-1}$ for every point $x\in\partial\mx_{-1}$.

\subsection{AdS Bending Cocycle}
The original idea of bending a spacelike plane in $\mx_{-1}$ was
already sketched in \cite{M}.  We go deeply in studying such a
notion and we relate it to the bending cocycle notion of Epstein and
Marden.\\

Let us describe first a {\it rotation around a spacelike geodesic
$l$}. By definition such a rotation is simply an isometry $T$ which
point-wise fixes $l$.  Up to isometries $l$ can be supposed lying on
$P_0=P(Id)$.  We know that the dual geodesic $l^*$ is a hyperbolic
$1$-parameter subgroup.

\begin{lem}\label{adesitt:bend:lem}
Let $l$ be a geodesic contained in $P_0$ and $l^*$ denote its dual
line.  For $x\in l^*$, the pair $(x,
x^{-1})\in PSL(2,\R)\times PSL(2,\R)$ represents a rotation around
$l$. The map
\[   R:l^*\ni x\mapsto (x,x^{-1})\in PSL(2,\R)\times PSL(2,\R)\]
is an isomorphism onto the subgroup of rotations around $l$. 
\end{lem}
\Dim First of all, let us show that the map
\[  \mx_{-1}\ni y\mapsto xyx\in\mx_{-1}\]
fixes point-wise $l$ (clearly $l$ is invariant by this transformation
because so is $l^*$). If $c$ is the axis of $x$ in $\mh^2$, we have
seen that $l$ is the set of rotations by $\pi$ around points in
$c$. Thus it is enough to show that if $p$ is the fixed point of $y$
then \[ xyx(p)=p \ .\] If we orient $c$ from the repulsive fixed point of
$x$ towards the attractive one, $x(p)$ is obtained by translating $p$
along $c$ in the positive direction, in such a way that $d(p,x(p))$ is
the translation length of $x$.  Since $y$ is a rotation by $\pi$ along
$p$, we have that $yx(p)$ is obtained by translating $p$ along $c$ in
negative direction, in such a way $d(p,yx(p))=d(p,x(p))$.  Thus we get
$xyx(p)=p$.\\ Now $R$ is clearly injective. On the other hand, it is
not hard to see that the group of rotations around a geodesic has
dimension at most $1$ (for the differential of a rotation at $p\in l$
fixes the vector tangent to $l$ at $p$).  Thus $R$ is surjective onto
the set of rotations around $l$. \cvd

\begin{cor}
Rotations around a geodesic $l$ act freely and transitively on the
dual geodesic $l^*$. Such action induces an isomorphism between the
set of rotations of $l$ and the set of translations of $l^*$.  \par
Moreover, by duality, we have that rotations around $l$ acts freely
and transitively on the set of spacelike planes containing $l$.  Thus,
given two spacelike planes $P_1,\ P_2$ such that $l\subset P_i$, then
there exists a unique rotation $T_{1,2}$ around $l$ such that
$T_{1,2}(P_1)=P_2$.
\end{cor}  \cvd

Given two spacelike planes $P_1,\ P_2$ meeting each other along a
geodesic $l$, we know that the dual points $x_i=x(P_i)$ lie on the
geodesic $l^*$ dual to $l$.  Then we define the {\it angle between
$P_1$ and $P_2$} as the distance between $x_1$ and $x_2$ along $l^*$.
Notice that:
\smallskip

{\it The angle between two spacelike 
planes is a well-defined number in $(0,+\infty)$}.
\smallskip

\noindent This is a difference with respect to the hyperbolic case,
that shall have important consequences on the result of the bending
procedure.
\begin{figure}
\begin{center}
\input{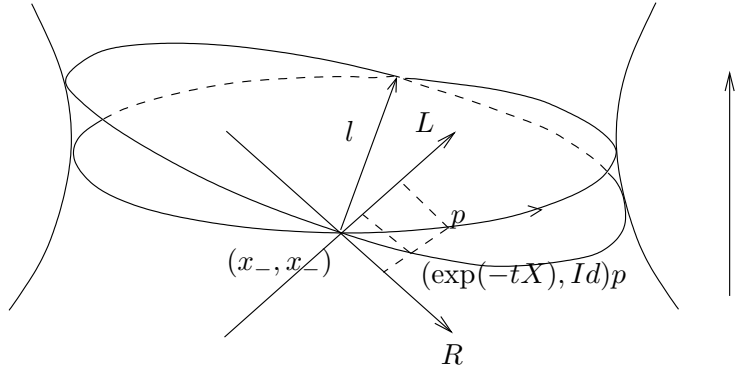}
\caption{{\small $(\exp(-tX),Id)$ rotates planes around $l$ in the
positive sense.}}\label{ADS:rot:fig}
\end{center}
\end{figure}
\begin{cor} 
An isometry $T$ of $\mx_{-1}$ is a rotation around a geodesic if and
only if it is represented by a pair $(x,y)$ such that $x$ and $y$ are
hyperbolic transformations with the same translation length.  \par
Given two spacelike planes $P_1,P_2$ meeting along a geodesic $l$,
let $(x,y)$ be the rotation taking $P_1$ to $P_2$. If $\tau$ is the
translation length of $x$, then the angle between $P_1$ and $P_2$ is
$\tau$.
\end{cor}
\Dim Suppose $(x,y)$ to be a pair of hyperbolic transformations with
the same translation length. Then there exists $z\in PSL(2,\R)$ such
that $zyz^{-1}=y^{-1}$.  Hence $(1,z)$ conjugates $(x,y)$ into
$(x,x^{-1})$. Thus $(x,y)$ is the rotation along the geodesic
$(1,z)^{-1}(l)$ where $l$ is the axes of $(x,x^{-1})$.  \par
Conversely, if $(x,y)$ is a rotation, it is conjugated to a
transformation $(z,z^{-1})$ with $z$ a hyperbolic element of
$PSL(2,\R)$. Thus we obtain that $x$ and $y$ are hyperbolic
transformations with the same fixed points.  \par In order to make the
last check, notice that, up to isometry, we can suppose
$P_1=P(Id)$. Thus, if $(x,x^{-1})$ is the isometry taking $P_1$
onto $P_2$, then the dual points of $P_1$ and $P_2$ are $Id$ and $x^2$
respectively. Now if $d$ is the distance of $x^2$ from $Id$ there exists a
unitary spacelike element $X\in\sG\lG(2,\mr)$ such that
\[
    x^2=\ch d I +\sh d X
\]
Thus we obtain that $\tr x^2=2\ch d$. On the other hand we know that 
$\tr x^2=2\ch u/2$ where $u$ is the translation length of $x^2$. Since
$u=2\tau$ the conclusion follows. \cvd
\smallskip

If we orient a spacelike line $l$, there is a natural definition of
positive rotation around $l$ (depending only on the orientations of
$l$ and $\mx_{-1}$).  Thus, we can induce an orientation on the dual
line $l^*$ by requiring that positive rotations act by positive
translations on $l^*$.  

In particular, if we take an oriented geodesic
$l$ in $P(Id)$, and denote by $X$ the infinitesimal generator
of positive translations along $l$ then it is not difficult to show that the
positive rotations around $l$ are of the form $(\exp(-tX),\exp(tX))$ for
$t>0$. Actually by looking at the action on the boundary we can easily deduce
that both the maps $(\exp(-tX),Id)$ and
$(Id,\exp(tX))$ rotate planes through $l$ in the positive direction
(see Fig.~\ref{ADS:rot:fig}).\\

We can finally define the bending at a measured geodesic lamination.
First, take a finite measured geodesic lamination $\lambda$ of
$\mh^2$.  Take a pair of points $x,y \in\mh^2$ and enumerate the
geodesics in $\lambda$ that cut the segment $[x,y]$ in the natural way
$l_1,\ldots,l_n$.  Moreover, we can orient $l_i$ as the boundary of
the half-plane containing $x$.  With a little abuse, denote by $l_i$ 
also the geodesic in $P(Id)$
corresponding to $l_i$, then let $\beta(x,y)$ be the isometry of
$\mx_{-1}$ obtained by compositionof positive rotations around $l_i$
of angle $a_i$ equal to the weight of $l_i$.  In particular, if $X_i$
denote the unit positive generator of the hyperbolic transformations with
axis equal to $l_i$, then we have
\[\begin{array}{l}\beta_\lambda(x,y)=
(\beta_-(x,y),\beta_+(x,y))\in PSL(2,\R)\times PSL(2,\R)\qquad\textrm
{where}\\\beta_-(x,y)=\exp(-a_1 X_1/2)\circ\exp(-a_2
X_2/2)\circ\ldots\circ \exp(-a_n X_n/2)\\ \beta_+(x,y)=\exp(a_1
X_1/2)\circ\exp(a_2 X_2/2) \circ\ldots\circ\exp(a_n
X_n/2)\end{array}\] with the following possible modifications: $a_1$
is replaced by $a_1/2$ when $x$ lies on $l_1$ and $a_n$ is replaced by
$a_n/2$ when $y$ lies on $l_1$ The factor $1/2$ in definition of
$\beta_\pm$ arises because the length translation of $\exp tX$ is
$2t$. \par Notice that $\beta_-$ and $\beta_+$ are the Epstein-Marden
cocycles corresponding to the {\it real-valued} measured laminations
$-\lambda$ and $\lambda$.  Thus, we can define in general a {\it
bending cocycle}
\[   \beta_\lambda(x,y)=(\beta_-(x,y), \beta_+(x,y))\]
where $\beta_-$ and $\beta_+$ are the Epstein-Marden cocycles associated to 
the real-valued measured geodesic laminations $-\lambda$ and $\lambda$.

\begin{remark}\label{E-M-lambda}
{\rm We stress that the above $-\lambda$ is just obtained from
$\lambda = (\Ll,\mu)$ by taking the {\it negative} ``measure''
$-\mu$. Although this is no longer a measured lamination in the sense
of Section \ref{laminations}, the construction of \cite{Ep-M} does
apply. In Section \ref{der} (and in the Introduction) we have used the
notation $-\lambda$ in a different context and with a different
meaning.}
\end{remark} 

\par The map $\beta_\lambda$ verifies the following properties:
\begin{enumerate}
\item $\beta_\lambda(x,y)\circ\beta_\lambda(y,z)=\beta_\lambda(x,z)$ 
for every $x,y,z\in\mh^2$ (this means that $\beta_\lambda$ is a 
$PSL(2,\R)\times PSL(2,\R)$-valued cocycle);
\item $\beta_\lambda(x,x)=Id$;
\item $\beta_\lambda$ is constant on the strata of the stratification
determined by $\lambda$.
\item If $\lambda_n\rightarrow \lambda$ on a $\eps$-neighbourhood of
the segment $[x,y]$ and $x,y\notin L_W$, then
$\beta_{\lambda_n}(x,y)\rightarrow\beta_{\lambda}(x,y)$.
\end{enumerate}

Take a base point $x_0$ in $P(Id)$. The {\it bending map} of $P(Id)$
with base point $x_0$ is simply
\[\varphi_\lambda: P(Id)\ni x\mapsto\beta_\lambda(x_0,x)x.\]
\begin{prop}
The map $\varphi_\lambda$ is an isometric $\mathrm C^0$-embedding of
$\mh^2$ into $\mx_{-1}$. The closure of its image in
$\overline\mx_{-1}$ is a closed disk and its boundary is an achronal
curve $c_\lambda$ of $\partial\mx_{-1}$.  The convex hull $K_\lambda$
of $c_\lambda$ in $\mx_{-1}$ has two boundary components: the past and
the future boundaries (resp. $\partial_-K_\lambda$
and $\partial_+K_\lambda$).  The map $\varphi_\lambda$ is an isometry
of  $\mh^2$ onto $\partial_+K_\lambda$.
\end{prop}
\Dim By using a sequence of standard approximations, we can show that
$\varphi_\lambda$ is a local isometrical embedding of $\mh^2$ into
$\mx_{-1}$.  Now, by a result proved in~\cite{M}, a local isometrical
embedding of a complete Riemannian surface into $\mx_{-1}$ is a an
isometrical embedding, and the closure of the image in
$\overline\mx_{-1}$ is homeomorphic to a closed disk with an achronal
boundary contained in $\partial\mx_{-1}$.  \par Denote by $c_\lambda$
the boundary of the image of $\varphi_\lambda$.  By Lemma 7.5
of~\cite{M}, we know that there exists a plane $P$ disjoint from
$c_\lambda$.  Thus, we can consider the convex hull $K_\lambda$ of
$c_\lambda$ in $\mr^3=\mathbb P^3-P$ (such a set does not depend on
the choice of the plane $P$).  Since $c$ is achronal, it is not hard
to see that $K_\lambda$ is contained in $\mx_{-1}$, and $\partial K_
\lambda \cap \partial \mx_{-1}=c_\lambda$.  Actually $K_\lambda$ is the convex
hull of $c_\lambda$.  We want to show that the image of
$\varphi_\lambda$ is $\partial_+K$.  For a point $x\in\mh^2$ denote by
$F_x$ the stratum through $x$.  We have that
$\varphi_\lambda(F_x)=\beta(x_0,x)I(F_x)$. Thus $\varphi_\lambda(F_x)$
is the convex hull of its boundary points. On the other hand, that
points are on $c_\lambda$ and so we have that $\varphi_\lambda(x)$ lie
in $K_\lambda$.  \par Let $P_x$ denote the plane $\beta(x_0,x)(P(Id))$
and $Q_x$ denote the plane dual to $\beta_\lambda(x)$. Clearly $P_x$
and $Q_x$ are disjoint (in fact the dual point of $Q_x$ lies in
$P_x$). Thus $\mx_{-1}-(P_x\cup Q_x)$ is the union of two cylinders
$C_-(x),C_+(x)$.  The closure of $C_-(x)$ has past boundary equal to
$Q_x$, whereas the past boundary of $C_+(x)$ $P_x$.  Now it is not
hard to show that the image of $\varphi_\lambda$ is contained in the
closure of $C_-(x)$.  Actually this is evident if $\lambda$ is a
finite lamination and the general statement follows by an usual
approximation argument.  It follows that $c_\lambda$ is contained in
$C_-(x)$ and so $K_\lambda$ too.  Since $\varphi_\lambda(x)$ is in the
future boundary of $C_-(x)$, it follows that it is in the future
boundary of $K_\lambda$. \cvd
\smallskip

Let $\Uu=\Uu^0_\lambda$ be the flat Lorentzian spacetime 
corresponding to $\lambda$, as in Section \ref{flat}.
Just as in the hyperbolic case we want to ``lift'' the bending cocycle
$\beta_\lambda$ to a continuous bending cocycle
\[   \hat\beta_\lambda:\Uu \times\Uu \rightarrow PSL(2,\R)\times 
PSL(2,\R) \ .\]
The following lemma is an immediate consequence of Lemma~3.4.4~(Bunch
of geodesics) of~\cite{Ep-M} applied to real-valued measured geodesic
laminations.

\begin{lem}\label{adesitt:bend:cont:lem}
For any compact set $K$ in $\mh^2$ and any $M>0$, there exists a
constant $C>0$ with the following property. Let
$\lambda=(\Ll,\mu)$ be a measured geodesic lamination of $\mh^2$
such that $\mu^*(\Nn(K))<M$. For every $x,y\in K$ and every geodesic
line $l$ of $\Ll$ that cuts $[x,y]$, let $X$ be the unit
infinitesimal positive generator of hyperbolic group of axis $l$ and
$m$ be the total mass of $[x,y]$. Then we have
\[     ||\beta_\lambda(x,y)-(\exp(mX),\exp(-mX))||<Cmd_\mh(x,y).\]
(We consider on $\ PSL(2,\R)\times PSL(2,\R)$ the product norm of
the norm of $\ PSL(2,\R)$.)
\end{lem} \cvd

\par By using this lemma we can prove the analogous of
Proposition ~\ref{lift:cocycle}. The proof is similar 
so we omit the details.

\begin{prop}\label{adesitt:bend:ext:prop}
A determined construction produces a continuous cocycle
\[\hat\beta_\lambda:\Uu(1)\times \Uu(1)\rightarrow PSL(2,\R)
\times PSL(2,\R)\]
such that
\[    \hat\beta_\lambda(p,q)=\beta_\lambda(N(p),N(q))\]
for $p,q$ such that $N(p)$ and $N(q)$ do not lie on $L_W$.  Moreover,
the map $\hat\beta_\lambda$ is locally Lipschitzian.  For every
compact subset $K$ on $\Uu(1)$, the Lipschitz constant on $K$ depends
only on $N(K)$ and on the measure $\mu^*(\Nn(N(K)))$.
\end{prop}  \cvd

\noindent
Finally we can extend the cocycle $\hat\beta$ on the whole $\Uu$ by
requiring that it is constant along the integral geodesics of the
gradient of the cosmological time $T$. In particular, if we set
$r(1,p)=r(p)+N(q)\in\Uu(1)$, then let us define
\[    \hat\beta(p,q)=\hat\beta(r(1,p),r(1,q)).\]

\begin{cor}\label{adesitt:bend:ext:prop}
The map
\[    \hat\beta_\lambda:\Uu\times \Uu\rightarrow PSL(2,\R)\times PSL(2,\R)\]
is locally Lipschitzian (with respect to the Euclidean distance on
$\Uu$).  Moreover, the Lipschitz constant on $K\times K$ depends only
on $N(K)$,$\mu^*(\Nn(N(K)))$ and the maximum and the minimum of $T$ on
$K$.  \par If $\lambda_n\rightarrow\lambda$ on a $\eps$-neighbourhood
of a compact set $H$ of $\mh^2$, then $\hat B_{\lambda_n}$ converges
uniformly to $\hat B_\lambda$ on $U(H;a,b)$ (that is the set of points
in $\Uu$ sent by Gauss map on $H$ and with cosmological time in the
interval $[a,b]$).
\end{cor}\cvd

\subsection{The canonical AdS rescaling}
In this section we define a map
\[\Delta_\lambda:\Uu^0_\lambda\rightarrow\mx_{-1}\]
such that the pull-back of the Anti de Sitter metric is a rescaling of
the flat metric along the gradient of the cosmological time.  A main
difference with respect to the WR map of Section \ref{hyp} shall be
that the AdS developing map $\Delta_\lambda$ is always an {\it
embedding} onto a determined convex domain $\Pp_\lambda$ of
$\mx_{-1}$.

\paragraph{The AdS spacetime $ \Uu_\lambda^{-1}$}
We know that the image of $\varphi_\lambda$ is a spacelike
$\mathrm C^0$-surface. We can consider its {\it domain of dependence}
$\Uu_\lambda^{-1}$, that is the set of points such that every causal
curve starting from them intersects the image of $\varphi_\lambda$.

\par It is possible to characterize $\Uu_\lambda^{-1}$ in a quite easy
way:
\smallskip

{\it $\Uu_\lambda^{-1}$ is the set of points $x$ such that  the boundary of 
the dual plane $P(x)$ does not intersect $c_\lambda$.}
\smallskip

\noindent
It is not hard to see that $\Uu_\lambda^{-1}$ is a convex set.  \par
Clearly $\partial_+ K_\lambda$ is contained in $\Uu_\lambda^{-1}$.  On
the other hand, it is not difficult to show that $\partial_-K_\lambda$
is contained in the closure of $\Uu_\lambda^{-1}$, and a point
$x\in\partial_-K_\lambda$ does not lie in $\Uu_\lambda^{-1}$ if and
only if there is a null support plane of $K_\lambda$ at $x$.  Since
both the boundaries of $K_\lambda$ are contained in the closure of
$\Uu_\lambda^{-1}$, we have that $K_\lambda$ is contained too.  \par
For a point $x\in\partial K_\lambda$, the dual plane $P(x)$ is
disjoint from $\Uu_\lambda^{-1}$. It follows that $\Uu_\lambda^{-1}$
does not contain any closed timelike curve. Thus, each component of
$\partial K_\lambda$ is a Cauchy surface, and $\Uu_\lambda^{-1}$ is
globally hyperbolic.  \par On the other hand, if we take a support
plane for $\partial K_\lambda$ at $x$, its dual point is on the
boundary of $\Uu_\lambda^{-1}$.

\paragraph {The past-side of $\Uu^{-1}_\lambda$}
Let us define $\Pp_\lambda$ to be the past of the surface
$\partial_+K_\lambda=\varphi_\lambda(\mh^2)$ in $\Uu_\lambda^{-1}$.

Take a point $p\in\Pp_\lambda$ and let $\Gg^+(p)$ be the set of points
related to $p$ by a future-pointing non-spacelike geodesic of length
less than $\pi/2$. For every point $q\in\Gg^+(p)$, there exists a
unique timelike geodesic contained in $\Gg^+(p)$ joining $p$ to
$q$. We define $\tau(q)$ the length of such a segment. It is not hard
to show that $\tau$ is a strictly concave function, so its level
surfaces are convex from the past. Since $p\in\Pp_\lambda$, we have
that $\Gg^+(p)\cap\partial_+K_\lambda$ is a compact disk. So there
exists a maximum of $\tau$ restricted to such a disk. Since
$\partial_+K_\lambda$ is convex from the past, such a point is unique,
and we denote it by $\rho_+(p)$. The plane orthogonal to the segment
$[p, \rho_+(p)]$ at $\rho_+(p)$ is a support plane for $K_\lambda$ at
$\rho_+(p)$.  In fact this property characterizes the point
$\rho_+(p)$.  We denote by $\rho_-(p)$ the dual point of that
plane. We have that $p$ lies on the future-pointing timelike segment
joining $\rho_-(p)$ to $\rho_+(p)$.

\par There is a foliation of $\Pp_\lambda$ by timelike geodesics: for
every $p$ consider the geodesic $[\rho_-(p),\rho_+(p)]$. Thus we can
define the function $\tau$ on $\Pp_\lambda$ by setting $\tau(p)$ as
the proper time of the segment $[\rho_-(p),p]$.  By direct
computation, it follows that $\tau$ is a $\mathrm C^1$-submersion of
$\Pp_\lambda$ taking values on $(0,\pi/2)$, and the level surfaces are
spacelike surfaces orthogonal to the geodesic foliation (that is
obtained by taking the integral line of the gradient of $\tau$).
Moreover, the map
$\rho_+:\Pp_\lambda(\tau=k)\rightarrow\partial_+K_\lambda$ is
Lipschitzian and proper. In particular, the level surfaces are
complete so that $c_\lambda$ is the boundary of all of them.

\paragraph {The AdS developing map $\Delta_\lambda$}
The idea is to find a diffeomorphism from $\Uu^0_\lambda$ to
$\Pp_\lambda$ which sends the integral lines of the gradient of the
cosmological time onto integral lines of the gradient of $\tau$, and
level surfaces of the cosmological time onto $\tau$-level surfaces.\\
For every $p\in\Uu^0_\lambda$, we define $x_-(p)$ as the dual point of
the plane $\hat\beta_\lambda(p_0,p)(P(Id))$, and
$x_+(p)=\hat\beta_\lambda(p_0,p)(N(p))$.  Thus let us choose
representatives $\hat x_-(p)$ and $\hat x_+(p)$ in $\SL{2}{R}$ such
that $\hat x_+(p)$, the geodesic segment between $\hat x_-(p)$and
$\hat x_+(p)$, is future directed.  Let us set
\[    \Delta_\lambda(p)=[\cos\tau(p)\hat x_-(p)+\sin\tau(p)\hat x_+(p)]\]
where $\tau(p)= \arctan T(p)$ .

\begin{teo}\label{adesitt:rescaling:teo} The map
\[    \Delta_\lambda:\Uu^0_\lambda \rightarrow \mx_{-1}\]
is a $\mathrm C^1$-diffeomorphism onto $\Pp_\lambda$.  Moreover, the
pull-back of the Anti de Sitter metric is equal to the rescaling of
the flat Lorentzian metric, directed by the gradient of the
cosmological time $T$, with rescaling functions:
\[\begin{array}{ll}   \alpha =\frac{1}{(1+T^2)^2} \ , \ \  &   
\beta=\frac{1}{1+T^2} \ .\end{array}\]
\end{teo}
\Dim Clearly $\Delta_\lambda$ is continuous and takes values on
$\Pp_\lambda$.  Now, for a point $x\in\Pp_\lambda$, let us set
$\rho_+(x)$ and $\rho_-(x)$ as above. We know that there exists
$y\in\mh^2$ such that $\beta(x_0,y)I(y)=\rho_+(x)$.  If $y$ does not
lie on $L_W$, then there exists a unique support plane at $\rho_+(x)$
equal to $\beta(x_0,y)(\mh^2)$, and $\rho_-(x)$ is its dual point.
Moreover, $N^{-1}(y)$ is a single geodesic $c$ in $\Uu^0_\lambda$. The
image of $c$ through $\Delta_\lambda$ is the geodesic segment between
$\rho_-(x)$ and $\rho_+(x)$.  \par Suppose that $\rho_+(x)\in
L_W$. Then there exists a unique point $p\in\Uu^0_\lambda(1)$ such
that $\hat\beta(p_0,p)(\mh^2)$ is the plane dual to $\rho_-(x)$.
Thus, if $c$ denote the integral line of the gradient of $T$, the
image $\Delta_\lambda(c)$ is the geodesic between $\rho_-(x)$ and
$\rho_+(x)$.  Finally, we have that $\Delta_\lambda$ is surjective.
An analogous argument shows that it is injective.\\ In order to
conclude the proof it is sufficient to analyze the map
$\Delta_\lambda$ in the case when $\lambda$ is a single weighted
geodesic. In fact, if we prove the statements of the theorem in that
case, then the same result will be proved when $\lambda$ is a finite
lamination.  Finally, by using standard approximations, we can achieve
the proof of the theorem.\\ 
\begin{figure}
\begin{center}
\input{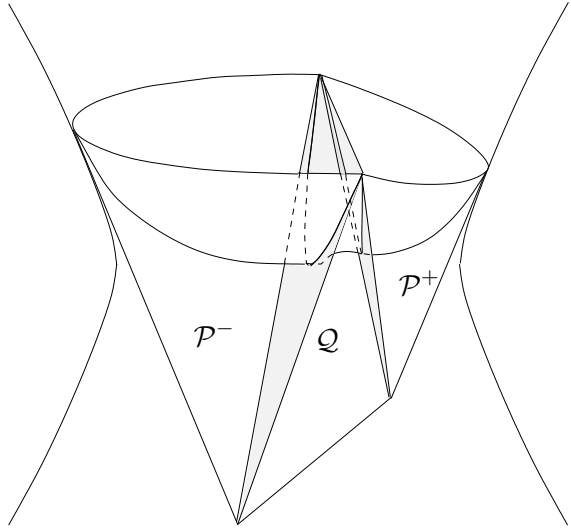}
\caption{{\small The domain $\Pp$ with its decomposition. Also the surface
    $\Pp(a)$ is shown.}}
\end{center}
\end{figure}
Let us set $\lambda=(l_0,a_0)$ and choose
a base point $x_0\in\mh^2-l_0$.  The surface $\partial_+K_\lambda$ is
simply the union of two half-planes $P_-$and $P_+$ meeting each other
along a geodesic (that, with a little abuse of notation, we will
denote by $l_0$).  We can suppose $x_0$ to be in $P_-$, and $l_0$
oriented as the boundary of $P_-$.  Let $u_\pm$ denote the dual points
of the planes containing $P_\pm$. Then the image of the map $\rho_-$ is
the segment on $l_0^*$ with end-points $u_-$ and $u_+$.  For $x
\in\Pp=\Pp_{\lambda_0}$, let $\tau(x)$ denote the Lorentzian length of
the segment $[\rho_-(x),x]$.  It is not difficult to see that $\tau$
is a $(0,\pi/2)$-valued $\mathrm C^1$-submersion and the gradient of
$\tau$ at $x$ is a past pointing unit timelike vector tangent to the
segment $[\rho_-(p),\rho_+(p)]$.  Denote by $\Pp(a)$ the level surface
$\tau^{-1}(a)$. It is a spacelike Cauchy surface. We can consider the
decomposition of $\Pp(a)$ given by
\[\begin{array}{l} \Pp^-(a)=
\Pp(a)\cap\overline{\rho_+^{-1}(\mathrm{int}P_-)} \ ,\\ \Qq(a)=
\Pp(a)\cap{\rho_+^{-1}(l_0)} \ ,\\ \Pp^+(a)=
\Pp(a)\cap\overline{\rho_+^{-1}(\mathrm{int}P_+)} \ .
\end{array}\] 
It is
not hard to show that the map $\rho_+$ restricted to $\Pp^{\pm}(a)$ is
a dilation of a factor $(\cos a)^{-1}$, whereas the map
\[   \Qq(a)\ni x\mapsto(\rho_-(x),\rho_+(x))\in[u_-,u_+]\times l_0\]
sends the metric on $\Qq(a)$ onto the metric
\[ (\sin a)^2\d u + (\cos a)^2\d t\]
where $u$ and $t$ are the natural parameters on $[u_-,u_+]$ and$l_0$.
\par Denote by $\delta_a$ the intrinsic distance on the surface
$\Pp(a)$, and by $l_a$ the boundary of $\Pp^-(a)$. For every point
$x\in\Pp(a)$, denote by $\pi(x)$ the point on $l_a$ that realizes the
distance of $l_a$ from $x$.  Let us fix a point $y_0\in l_0$, and
denote by $y_a$ the point on $l_a$ such that $\rho_+(y_a)=y_0$. Then
consider the functions
\[ \begin{array}{l} 
\sigma(x)= \eps(x)\delta_{\tau(x)}(x,l_{\tau(x)} )/\sin(\tau(x))\\     
\xi(x)=\eps(x)\delta_{\tau(x)}(\pi(x), y_{\tau(x)})/\sin(\tau(x))
\end{array}\] 
where $\eps(x)=-1$ if $x\in\Pp^-$ and $\eps(x)=1$ otherwise.  The set
of functions $(\tau,\sigma,\xi)$ furnishes $\mathrm C^1$-coordiantes on
$\Pp$.  Actually it $v_0$ denote the tangent vector at $y_0$ of $l_0$
then it is not hard to see that the induced parametrization is given
by formula
\[
  (\tau,\sigma,\xi)\mapsto\left\{\begin{array}{ll} \sin
                                  t\left(\ch\xi(\ch\sigma
                                  y_0\,+\,\sh\sigma v_0)\, +\, \sh\xi
                                  X\right)\, +\, \cos t Id \\ \textrm{
                                  if }\xi<0\ ;\\ \sin t
                                  \left(\ch\sigma y_0+\sh\sigma
                                  v_0\right)\, +\, \cos t \exp(\xi\tan
                                  t X)\\ \textrm{ if }
                                  \xi\in[0,\alpha_0/\tan t]\ ;\\ \sin
                                  t\left(\ch\xi(\ch\sigma y_0
                                  +\sh\sigma v_0)\right) +\cos t
                                  \exp(-\alpha_0 X/2)X+ \exp(-\alpha_0
                                  X_0)\\ \textrm{ otherwise}\ .
                                  \end{array}\right.
\]

By using these coordinates on
$\Pp$, and coordinates $(T,u,\zeta)$ on $\Uu_0$, we
have
\[   \begin{array}{l}    
\tau(\Delta_0(T,u,\zeta))=\arctan(T)\\    
\sigma(\Delta_0(T,u,\zeta))=u\\    
\xi(\Delta_0(T,u,\zeta))=\zeta.    
\end{array}\]
Then, by direct computation in these coordinates, it is not hard to
obtain that the pull-back on $\Uu^0_\lambda$ of the Anti de Sitter
metric is the rescaling of the flat metric, directed by the gradient
of $T$, with rescaling functions
\[   \begin{array}{lr}    
\alpha=\frac{1}{(1+T^2)^2}\ ,\ \ \beta=\frac{1}{(1+T^2)}\ .
\end{array}\] 
\cvd With a bit abuse of notation, we denote by
$\Pp_\lambda$ also the domain $\Uu^0_\lambda$ endowed with the Anti de
Sitter metric induced by $\Delta_\lambda$.  We know that $\Pp_\lambda$
is globally hyperbolic and the level surfaces
$\Uu^0_\lambda(a)=\Pp_\lambda(\tau=\arctan a)$ are Cauchy surfaces.
By using the argument of Proposition~\ref{desitter:ct:prop} we can
easily compute the cosmological time on $\Pp_\lambda$.
\begin{cor}\label{adesitt:ct:prop} The function \[   \tau=\arctan(T)\]
is the cosmological time on $\Pp_\lambda$. The initial singularities
of $\Uu^0_\lambda$ and $\Pp_\lambda$ coincide.
\end{cor}\cvd
\smallskip

\paragraph{WR: AdS towards hyperbolic geometry}
By composing $\Delta^{-1}$ with the WR constructed in Section
\ref{hyp}, we obtain a WR defined on
$\Uu^{-1}_\lambda(]\pi/4,\pi/2[)\subset \Pp_\lambda$, directed by the
gradient of the cosmological time $\tau$.  Moreover, this WR extends
continuously to the closure of $\Uu^{-1}_\lambda(]\pi/4,\pi/2[)$,
producing a local homeomorphism onto the completion
$\overline{M}_\lambda$ of the hyperbolic manifold $M_\lambda$.  In
particular, this gives us a local isometry of the surface
$\partial_+K_\lambda$ onto the hyperbolic boundary of
$\overline{M}_\lambda$, which preserves the respective bending loci.

\subsection{AdS $\Mm\Ll(\mh^2)$- spacetimes}
Recall that the flat $\Mm\Ll(\mh^2)$-spacetimes $\Uu^0_\lambda$ are
characterized as the flat regular domains in $\mx_0$ with surjective
Gauss map (Section \ref{flat}).  In this section we want to
characterize the Anti de Sitter spacetimes $\Uu^{-1}_\lambda$, arising
by performing the above canonical rescaling $\Uu^0_\lambda\to
\Pp_\lambda \subset \Uu^{-1}_\lambda$.

Take a globally hyperbolic Anti de Sitter spacetime $M$ with a {\it
complete} spacelike Cauchy surface $S$, and suppose $S$ to be simply
connected.  Then the restriction to $S$ of a developing map $d$ of $M$
is an embedding onto a spacelike surface of $\mx_{-1}$, and its
closure in $\overline \mx_{-1}$ topologically is a closed disk. In
particular, its boundary, say $c$, is a achronal curve on the
boundary of $\mx_{-1}$.  Denote by $\hat M \subset \mx_{-1}$ the
domain of dependence of $d(S)$.  It is not hard to see that $d(M)$ is
contained in $\hat M$.  Moreover, by an argument of Choquet-Bruhat and
Geroch (see ~\cite{H-E} cap.7) we have that $d$ is an embedding of $M$
into $\hat M$.  Thus $\hat M$ is a maximal globally hyperbolic
spacetime.  Notice that $\overline{\hat M}\cap\partial\mx_{-1}=c$.
Moreover, $\hat M$ is determined by $c$.  In fact a point
$x\in\mx_{-1}$ lies in $\hat M$ if and only if the boundary of its
dual plane does not intersect $c$.\\

Take now a no-where timelike curve $c$ in the boundary of $\mx_{-1}$
and denote $\hat M$ the maximal globally hyperbolic spacetime
determined by $c$ as above. Denote by $K$ the convex hull of $c$ in
$\mx_{-1}$.  For every point in the interior of $K$, its dual plane
does not intersect $c$.  Thus the interior of $K$ is contained in
$\hat M$. This is called the {\it convex core} of the AdS spacetime
$\hat M$.  The boundary of $K$ in $\mx_{-1}$ has two connected
components, the {\it past boundary} $\partial_-K$ and the {\it future
boundary} $\partial_+K$.  The past of the future boundary of $K$ in
$\hat M$ is called the {\it past side} of $\hat M$.

\par In general $\partial_+K$ is an achronal surface. We say that it
is spacelike if every support plane of $K$ touching $\partial_+K$ is
spacelike.  In this case we have that there exists an intrinsic
distance defined on $\partial_+K$. In fact, since $\partial_+K$ is the
boundary of a convex set, then given $x,y\in\partial_+K$, there exists
a Lipschitzian path $c$ joining $x$ to $y$ (Lipschitzian with respect
to the Euclidean distance on$\mx_{-1}$). Moreover, the speed of $c$ at
$t$ is defined for almost all $t$ and lies on a support plane of
$\partial_+K$ at $c(t)$ (so it is a spacelike vector).  We can define
the length of $c$, and by consequence an intrinsic distance on
$\partial_+K$, by setting $\delta(x,y)$ the infimum of the
Lipschitzian path joining $x$ to $y$.  By using that $\partial_+K$ is
spacelike, it is not difficult to see that $\delta$ is locally
equivalent to the Euclidean distance.

We can finally state the characterization of $\Mm\Ll(\mh^2)$- AdS
spacetimes.

\begin{teo}\label{MLAdS} 
The canonical rescaling (Theorem~\ref{adesitt:rescaling:teo})
produces a bijection $\Uu^0_\lambda \leftrightarrow \Uu^{-1}_\lambda$
between the isometry classes of maximal globally hyperbolic flat
spacetimes with surjective Gauss map and the isometry classes of
maximal globally hyperbolic Anti de Sitter spacetimes with {\rm
complete} spacelike future convex-core boundary.
\end{teo}
The theorem will follows from the following proposition
\begin{prop} If $\partial_+K$  is spacelike then it is
locally $\mathrm C^0$-isometric to $\mh^2$.\\ If it is a spacelike
{\rm complete} surface, then there exists a measured geodesic
lamination $\lambda$ on $\mh^2$ such that $\partial_+K$ is the image
of $\varphi_\lambda$.
\end{prop}
\Dim Take a point $x\in\partial_+K$. We want to define a sequence
$S_n$ of spacelike surfaces obtained by bending a spacelike plane $P$
along a finite number of geodesics ``converging'' to $\partial_+K$
near $x$.  Take a plane $Q$ that does not intersect the boundary curve
$c$ and let $\mr^3$ be identified to $\mathbb P^3-Q$. Consider a
countable dense subset $\{x_n\}$ in $\partial_+K$ with $x_1=x$ and,
for every $n$, choose a support plane $P_n$ touching $\partial_+K$ at
$x_n$.  Then, let $K_n$ the convex set obtained by intersecting the
past of the planes $P_1,\ldots, P_n$ in $\mx_{-1}\cap\mr^3$. In
general, its future boundary $C_n$ is not a pleated surface (it may
contain vertices).  However, take its boundary curve $c_n$ and
consider the future boundary $S_n$ of its convex hull. It is not
difficult to show that $S_n$ is a finite pleated surface. Notice that
$S_n$ is contained in the intersection of the past of $C_n$ and the
future of $\partial_+K$. In particular, $S_n$ converges to
$\partial_+K$ in a neighbourhood of $x$.  \par Every compact
neighbourhood of $x$ in $\overline P$ is contained in the past side of
the domain of dependence of $S_n$, for $n$ sufficiently large.  Denote
by $\rho_n$ the future retraction on $S_n$. We have that $S_n$ is
isometric to $\mh^2$ and we can choose the isometry $I_n$ so that
$I_n\circ\rho_n(x)$ is a fixed point of $\mh^2$ (that does not depend
on $n$).  We claim that $I_n\circ\rho_n$ converges to a local isometry
of $\partial_+K$ to $\mh^2$.  Denote by $\theta_n$ the cosmological
time of the past side $\Pp_n$ of the domain of dependence of $S_n$,
and by $\theta_\infty$ the cosmological time of the the past side
$\Pp$ of the domain of dependence of $\partial_+K$.  For every
compact set $H$ in $\Pp \cup\partial_+K$ there exists a constant
$n$ such that $H\subset\Pp_n\cup S_n$.  Moreover, $\theta_n$ converges
to $\theta_\infty$ in $\mathrm C^0(H)\cap\mathrm
C^1(H\cap\Pp)$.  By using this fact we can prove that there
exists a constant $C$ such that $I_n\circ\rho_n$'s are $C$-Lipschitz
on $H\cap\partial_+K$.  Thus they converge to a map $J_\infty$. The
same argument shows that $J_\infty$ preserves the length of
Lipschitzian paths.  Thus it is a local isometry.\\ If $\partial_+K$
is complete, then the map $J_\infty$ is an isometry.  Moreover, the
bending locus on $\partial_+K$ produces a geodesic lamination $\Ll$ on
$\mh^2$. If $\Ll_n$ denote the geodesic lamination obtained
pushing-forward the bending locus of $S_n$, we can equip $\Ll_n$ with
a measure $\mu_n$ defining the weight of $l$ to be the bending angle
of the corresponding geodesic on $S_n$. If $c$ is a geodesic arc
transverse to $\Ll$ it is not hard to see that it will be transverse
to $\Ll_n$, for $n$ sufficiently large.  We want to show that
$\lambda_n=(\Ll_n,\mu_n)$ tends to a lamination $\lambda=(\Ll,\mu)$.
The construction of the measure $\mu$ is obtained in a way similar to
the one used in~\cite{Ep-M} for the boundary of the convex hull of a
topological circle at infinity in the hyperbolic space.\\ Let us take
a path $c$ on $\partial_+K$ transverse to the bending locus.  Up to
subdivision, $c$ can be supposed intersecting every bending line at
most once. Take a dense set of points $\{x_n\}$ on $c$ and choose for
every point, a support plane $P_n$. Notice that $P_i\cap P_j$ is
either empty or a geodesic line.  If $P_k$ meets $P_i\cap P_j$, then
$P_k \cap P_i = P_i\cap P_j=P_k\cap P_j$.  Thus the future boundary
$S_n$ of the convex set obtained intersecting the past of $P_i$ in
$\mx_{-1}$ for $i=1,\ldots,n$ is a finite bent surface.  Define
$\alpha_n$ to be the sum of the bending angles along all the bending
lines of $S_n$.\\ {\it We claim that $\alpha_n$ converges to
$\alpha_\infty$ and this number does not depend on choices we made}.
Assuming this last claim, then we can define the total mass of $c$ as
$\alpha_\infty$. In this way a measure transverse to the bending locus
is defined on $\partial_+K$.  Pushing forward this measure, we obtain
a measured geodesic lamination $\lambda$ of $\mh^2$. Since
$\alpha_\infty$ does not depend on the sequence of $P_i$ the measure
$\lambda$ is the limit of $\lambda_k$.  We see that the inverse of the
isometry $J_\infty$ is $\varphi_\lambda$ (up to post-composition with an
element of $PSL(2,\R)\times PSL(2,\R)$).\\ \par Just as in the
hyperbolic case, the claim follows from the following Lemma that is
the strictly analogous of Lemma 1.10.1 (Three planes) of~\cite{Ep-M}.

\begin{lem}\label{3planes}
Let $P_1,P_2,P_3$ be three spacelike planes in $\mx_{-1}$ without a
common point of intersection, such that any two intersect
transversely. Suppose that there exists a spacelike plane between
$P_2$ and $P_3$ that does not intersect $P_1$.  Then the sum of the
angles between $P_1$ and $P_2$ and $P_1$ and $P_3$ is less than the
angle between $P_2$ and $P_3$. \end{lem}

Assuming the lemma, the claim follows immediately.  In fact the
sequence of bending angles $(\alpha_n)$ is decreasing, so it admits
limit.  Clearly the limit does not depend on the choice of planes
(otherwise we should be able to construct another sequence that does
not admit limit).

It remains to prove the lemma.

\emph{Proof of Lemma \ref{3planes}:} Denote by $x_i$ the dual point of
$P_i$. The hypothesis implies that segment between $x_i$ and $x_j$ is
spacelike, but the plane containing $x_1$, $x_2$ and $x_3$ is
Lorentzian.  The existence of a plane between $P_2$ and$P_3$ implies
that every timelike geodesic starting at $x_1$ meets the segment
$[x_2,x_3]$. Thus there exists a point $u\in[x_2,x_3]$ and a unit
timelike vector $v\in T_u\mx_{-1}$ orthogonal to $[x_2,x_3]$, such
that $x_1=\exp_u(tv)$. Choose a lift of $[x_2,x_3]$, say $[\hat
x_2,\hat x_3]$, on $\hat\mx_{-1}$ and denote by $\hat u$ the lift of
$u$ on that segment.  Denote by $l_1$ (resp. $l_2$, $l_3$) the length
of the segment $[x_2,x_3]$ (resp. $[x_1,x_3]$, $[x_1,x_2]$). We have
that\[\ch l_i=|\E{\hat x_j}{\hat x_k}|\] where $\{i,j,k\}=\{1,2,3\}$.
Now we have that $\hat x_1=\cos t \hat u+\sin t \hat v$ so that

$|\E{\hat x_1}{\hat x_i}|=\cos t |\E{\hat u}{\hat x_i}|$. 
Hence
\[ \begin{array}{ll} \ch l_2+<\ch l_2' & \ch l_3<\ch l_3' \end{array}\]
where $l_2'$ and $l_3'$ are the length of $[x_2,u]$ and $[x_3,u]$.
Finally we have $l_2+l_3<l_2'+l_3'=l_1$. \cvd
\smallskip

The proof of the proposition is now complete. \cvd

\subsection{$\Gamma$-invariant rescaling}
Assume that $\lambda\in \Mm\Ll(\mh^2)$ is invariant by the action of a
group $\Gamma<PSL(2,\R)$ acting freely and properly discontinuously
on $\mh^2$.  That is, $\lambda$ is the pull-back of a measured
geodesic lamination on the hyperbolic surface $F= \mh^2/\Gamma$.
Denote by $\beta=(\beta_-,\beta_+)$ the
$PSL(2,\R)\times PSL(2,\R)$-cocycle associated to $\lambda$.  The
cocycles $\beta_-$ and $\beta_+$ can be used to construct
representations of $\Gamma$ into $PSL(2,\R)$.  In fact, if we choose
a base point $x_0\in \mh^2$, we can define
\[   \begin{array}{ll}    
\rho_-(\gamma)=\beta_-(x_0,\gamma x_0)\circ\gamma\ , \qquad &    
\rho_+(\gamma)=\beta_+(x_0,\gamma x_0)\circ\gamma \ .   
\end{array}\]
By using the bending rule it is quite easy to check that $\rho_-$ and
$\rho_+$ are representations, and that the respective conjugacy
classes do not depend on the choice of the base point.  If we consider
the product representation $\rho=(\rho_-,\rho_+)$ taking values into
$PSL(2,\R)\times PSL(2,\R)$, then the bending map $\varphi_\lambda$
is clearly $\rho$-equivariant (the same base point is used to define
both the bending map and the representation $\rho$).  \par It follows
that the domain $\Uu_\lambda^{-1}$ is preserved by the action of
$\Gamma$ on $\mx_{-1}$ induced by $\rho$.  Moreover, as the image of
$\varphi_\lambda$ is a Cauchy surface of $\Uu_\lambda^{-1}$, the
action of $\Gamma$ on that domain is proper and free and the quotient
is a maximal globally hyperbolic spacetime homeomorphic to
$F\times\mr$.  \par Remind that a representation $f_\lambda$ of
$\Gamma$ into the isometry group of Minkowski space $\mx_0$ is defined
in Section~\ref{flat}, in such a way that the domain $\Uu^0_\lambda$
is invariant by $f_\lambda$.  Notice that by construction the
rescaling developing map
\[    \Delta_\lambda:\Uu_\lambda^0\rightarrow\mx_{-1}\]
satisfies the following condition
\[    \Delta_\lambda\circ f_\lambda(\gamma)=
\rho(\gamma)\circ\Delta_\lambda \ .\]
In particular, the rescaling on $\Uu_\lambda^{0}/f_\lambda(\Gamma)$,
directed by the gradient of the cosmological time $T$, with rescaling
functions
\[  \begin{array}{ll}   \alpha=
\frac{1}{(1+T^2)^2}\ ,\qquad & \beta=\frac{1}{1+T^2} 
\end{array}
\] 
makes it isometric to the past side of $\Uu_\lambda^{-1}/\rho(\Gamma)$.

\subsection{Earthquakes and AdS bending}
In this section we want to study more closely the boundary curve
$c_\lambda$ that determines a given AdS $\Mm\Ll(\mh^2)$-spacetime.

Mess discovered in \cite{M} a deep relationship between the
classification of AdS globally hyperbolic spacetimes with compact
Cauchy surfaces and Thurston's {\it Earthquakes Theorem} \cite{Thu2},
that we are going to summarize.

Remind first, with the notations of the present paper, that $\lambda\in
\Mm\Ll(\mh^2)$ {\it produces a left earthquake} (having $\lambda$ as
shearing lamination) if the map
\[
    \Ee:\mh^2\ni x\mapsto\beta_+(x_0,x)x\in\mh^2
\] 
is bijective, and, in this case, there exists a unique continuous
extension on $S^1_\infty =\partial\mh^2$ that is a homeomorphism.
Similarly for the {\it right earthquake}.

Assume now that $\lambda$ is invariant under the action of a
cocompact Fuchsian group $\Gamma$. Then: 
\smallskip

(1) Representations $\rho_-$ and $\rho_+$ introduced in the above
Section are cocompact Fuchsian representations of the same genus of
$\Gamma$.  Moreover, all ordered pairs $(\rho_-,\rho_+)$ of cocompact
Fuchsian representations of a given genus $g\geq 2$ arise in this way
(\cite{M}). It is well known that there is an orientation preserving
homeomorphism of $\overline{\mh}^2 = \mh^2 \cup S^1_\infty$ which
conjugates the action of $\rho_-$ on $\overline{\mh}^2$ with the one
of $\rho_+$.  Let us denote by $u$ its restriction to $S^1_\infty$.
\smallskip

(2) Thurston proved \cite{Thu2} that there exists a unique
$\rho_-$-invariant $\lambda' \in \Mm\Ll(\mh^2)$ such that $\rho_+$ is
given by the left eartquake produced by $\lambda'$, and this last
continuously extends to $S^1_\infty$ by $u$.
\smallskip

(3) In \cite{M}, it is proved that $\lambda$ is the image of
$\lambda'/2 = (\Ll',\mu'/2)$ via the left earthquake produced by
$\lambda'/2$, and the boundary curve $c_\lambda$ coincides with the
graph of $u$.  Similar results hold for the unique right earthquake
from $\rho_+$ to $\rho_-$. This implies, in particular that $\lambda$
itself generates both left and right earthquakes, but this is not
a surprise in this case as it is a fact that any cocompact 
$\Gamma$-invariant lamination produces earthquakes.
\medskip

\noindent
Similar facts hold, for instance, when $\Gamma$ is not necessarily
cocompact, but $\lambda$ is the pull-back of a measured geodesic
lamination on $F=\mh^2/\Gamma$ with {\it compact support}. 

One might wonder the generalization of the above results to the full
range of application of Thurston's Earthquake theorem,
i.e. considering arbitrary orientation preserving homeomorphisms of
$\mh^2$ that continuously extend by homeomorphisms, $u$ say,  of the
boundary $S^1_\infty$, as in the above point (1).

We are going to see that this {\it not} exactly the case.
\smallskip

First we characterize the spacetimes $\Uu^{-1}_\lambda$ such that
the corresponding curve at infinity  $c_\lambda$ is the graph
of a homeomorphism (recovering in particular the cocompact case):
this happens iff $\lambda$ produces earthquakes.
In general, $c_\lambda$ is just an achronal curve (see also
the examples in Section \ref{3cusp}).
 
Later we will show examples of homeomorphism $u$ of $S^1$, such that
their graphs are {\it not}  $c_\lambda$ for any $\Uu^{-1}_\lambda$
(that is such a graph is associated to AdS spacetimes of more general
type). In other words, we show a lamination $\lambda'$ that produces
a left earthquake, but such that $\lambda'/2$ does {\it not}.

\begin{teo} Let $\lambda \in  \Mm\Ll(\mh^2)$. 
Then the following are equivalent statements:
\begin{enumerate}
\item The curve $c_\lambda$ is the graph of a homeomorphism of
$S^1_\infty$.
\item Both left and right earthquakes with shearing locus equal to
$\lambda$ are well-defined.\end{enumerate}
\end{teo}
\Dim If $\lambda$ is a finite lamination and $\Ee_L$ and $\Ee_R$
denote the left and right earthquake along $\lambda$, it is not hard
to see that
\[    c_\lambda=\{(\Ee_R(x),\Ee_L(x))|x\in S^1\}.\]
In particular, if $\lambda'$ is the image of $\lambda$ through
$\Ee_R$, we have that $c_\lambda$ is the graph of a left-earthquake
with shearing lamination equal to $2\lambda'$.  \par Suppose now that
left and right earthquakes along $\lambda$ are well defined.  Let
$\Delta\subset S^1_\infty$ denote the set of points that are in the
boundary of some piece of $\lambda$. Then it is easy to see that the
set
\[    \{(\Ee_R(x),\Ee_L(x))|x\in\Delta\}\]
is contained in $c_\lambda$. On the other hand, 
the closure of this set is
\[    C=\{(\Ee_R(x),\Ee_L(x))|x\in S^1_\infty\}\]
that is the graph of a homeomorphism of $S^1_\infty$. We can deduce
that $c_\lambda=C$ and one implication is proved.\\ Suppose now
$c_\lambda$ to be the graph of a homeomorphism $u$.  There exists a
left earthquake $\Ee$ along $\lambda'$ such that $\Ee|_{S^1}=u$.  On
the other hand we can approximate $\Ee$ by a sequence of left simple
earthquakes $\Ee_n$ with shearing lamination equal to $\lambda_n'$.
Denote by $\lambda_n$ the image of the lamination $\lambda_n'/2$ by a
left earthquake with shearing lamination equal to $\lambda_n'/2$.
\par We know that the future boundary of the convex hull $K_n$ of
$c_n$ is obtained by bending $\mh^2$ along $\lambda_n$.  It is not
hard to see that $u_n=\Ee_n|_{S^1}$ converges to $u$ at least on a
dense subset given by points in the boundary of pieces of
$\lambda'$. Thus if $c_n$ denote the graph of $u_n$ we have that $c_n$
converges to $c_\lambda$ in the Hausdorff topology of
$\overline\mx_{-1}$.  Thus $K_n$ converges to $K_\lambda$ in the
Hausdorff topology.  It follows that $\lambda_n$ converges to
$\lambda$.  Since $\lambda_n$ is a convergent sequence, it follows
that $\lambda'/2$ gives rise to a left earthquake and the image of
$\lambda'/2$ is $\lambda$. Thus $\lambda$ gives rise to a right
earthquake.  Moreover, since a left earthquake along $\lambda'$ is
defined, we have that the left earthquake along $\lambda$ exists.\cvd
\begin{cor}
If $c_\lambda$ is the graph of a homeomorphism, 
then \[   c_\lambda=\{(\Ee_R(x),\Ee_L(x))|x\in S^1_\infty\}.\]
\end{cor}\cvd
\smallskip
\begin{exa}\label{earthex}
{\rm We will make an example of a homeomorphism $u$ of $S^1_\infty$
such that the future boundary of the convex core of its graph $c$ is not a
complete surface.  In this way we will show an example of a
homeomorphism such that its domain of dependence is not $\Uu^{-1}_\lambda$
for any $\lambda\in\Mm\Ll(\mh^2)$.}\par
\begin{figure}
\begin{center}
\input{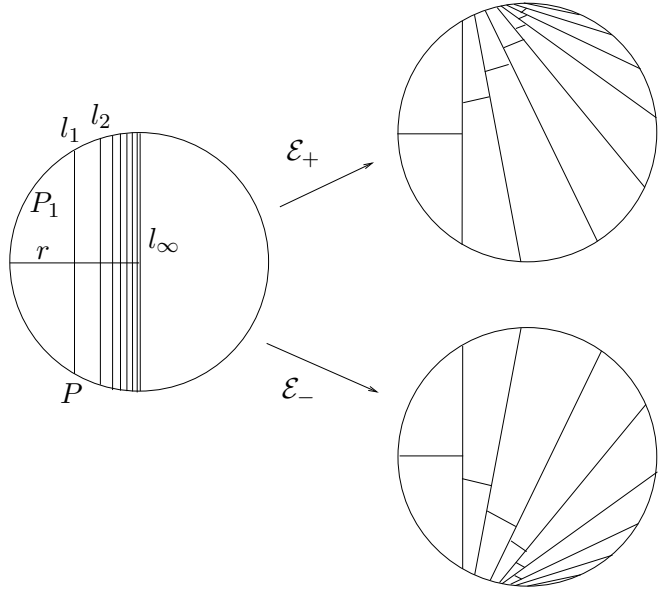}
\caption{{\small The maps $\Ee_\pm$ are injective. Since the geodesic
    $l_\infty$ escapes to infinity then they are also surjective.}}
\end{center}
\end{figure}
{\rm 
Take a geodesic ray $r$ in $\mh^2$  starting from $x_\infty$
and denote by $x_n$ the point on the ray  such that $d_\mh(x_\infty,x_n)=1/n$.
Let $l_n$ (resp. $l_\infty$) be the geodesic through $x_n$
(resp. $x_\infty)$ orthogonal to $r$ and  $P_n$ be the half-plane bounded 
by $l_n$ that does not contain $x_\infty$. Finally let $P$ be the
half-plane bounded by $l_\infty$ containing all $P_n$'s.
Orient $l_n$ as boundary of $P_n$ and denote by $X_n$
the generator of the positive translation along $l_n$.
Finally fix x a point $x_0\in P_1$. 
For every $x\in P$ the geodesic ray between $x_0$ and $x$ meets only a finite
number of $l_i$, say $l_1,\ldots,l_n$. Then define
\[
    \begin{array}{l}
    \beta_\pm(x)=\exp(\pm X_1)\circ\exp(\pm X_2)\circ\cdots\circ\exp(\pm
    X_n)\\
    \Ee_\pm(x)=\beta_\pm(x)x\ .
    \end{array}
\]
It is not hard to see that $\Ee_\pm$ are injective functions that
extend in a natural way to functions $\partial_\infty P\rightarrow
S^1_\infty$.  Moroever even if $\Ee_\pm$ are not continuous their
extension on the boundary are continuous. Now we claim that $\Ee_\pm$
is surjective.  The claim implies that the restriction of $\Ee_\pm$ on
$\partial_\infty P-\partial_\infty l_\infty$ is injective and sends
the end-points of $l_\infty$ on the same point.  Now the map}
\[
   f: P\ni x\mapsto (\beta_-(x),\beta_+(x))I(x)
\]
{\rm
(where $I:\mh^2\rightarrow P(Id)$ is the standard embedding) is an isometry
 onto  the future boundary of the convex core of the curve
}
\[
    c=\{(\Ee_-(x),\Ee_+(x))|x\in\partial_\infty P\}
\]
{\rm that is a graph of a homeomorphism.  Since $P$ is not complete so
is that boundary component.}\par {\rm Now let us sketch the proof of
claim.  If we set $g_n=\exp X_1\circ\cdots\circ\exp X_n$ we have to
prove that the sequence of geodesics $g_n(l_\infty)$ is divergent.
Now it is not hard to see the following facts}\par
\smallskip
{\rm 1. $g_n$ is a hyperbolic transformation and its axis is contained
in $P-P_1$.}\par
\smallskip
{\rm 2. The translation length of $g_n$ is greater than $n$.}\par
\smallskip
{\rm 3. The distance of $l_\infty$ from the axis of $g_n$ is greater then
$1/n$.}\par

{\rm If $z_n$ denote the point on the axis of $g_n$ from $l_\infty$,
  facts 2. and 3. imply that the distance of $g_n(l_\infty)$ from
  $z_n$ tends to $+\infty$. By point 1. we have that $z_n$ runs in a
  compact set and this imply that $g_n(l_\infty)$ is divergent.}\par

{\rm We want to remark that $\bigcup \Ee_\pm(l_i)$ are simplicial
geodesic laminations $\Ll_\pm$ of $\mh^2$. If we consider the weight
$1$ on every leaf of $\Ll_-$ we have a measured geodesic lamination
$\lambda_-$ that in a sense its the image of the initial
``lamination''.  Notice that $\lambda_-$ does not produce earthquake
whereas $2\lambda_-$ gives rise to the eartquake such that $c$ is the
graph of its extension on the boundary.  }
\end{exa}

\subsection{$T$-symmetry}\label{Tsym}
Let $(\Uu^{-1}_\lambda)^*$ the AdS spacetime obtained by reversing the
time orientation. It is not hard to verify that the effect on the
holonomy representation is just of exchanging
$$ (\rho_-, \rho_+) \leftrightarrow (\rho_+, \rho_-) \ .$$ When
$\lambda$ is invariant for a cocompact Fuchsian group $\Gamma$, it
follows from the results of the previous Section, that also
$(\Uu^{-1}_\lambda)^*$ is a AdS $\Mm\Ll(\mh^2)$-spacetime.  More
precisely
$$(\Uu^{-1}_\lambda)^*= \Uu^{-1}_{\lambda^*}$$ where $\lambda^*$ is
the bending lamination of the past boundary of $K_\lambda$.  Notice
that if $(\lambda^*)'$ produces a right earthquake conjugating
$\rho_-$ to $\rho_+$ then $\lambda^*$ is the image of $(\lambda^*)'/2$
via the earthquake along $(\lambda^*)'/2$.

So cocompact $\Gamma$-invariant AdS $\Mm\Ll(\mh^2)$-spacetimes are
invariant for this $T$-{\it symmetry}.  The same fact holds, for
instance, when $\Gamma$ is not necessarily cocompact, but we consider
laminations with compact support on the quotient hyperbolic surface.

On the other hand, the $T$-symmetry is broken for general
AdS $\Mm\Ll(\mh^2)$-spacetimes. In other words, although the past
boundary $\partial_-K$ and the future boundary $\partial_+K$ of the
convex core share the same boundary curve $c_\lambda$, it actually
happens in general that $\partial_-K$ is {\it not} complete, that is
$(\Uu^{-1}_\lambda)^*$ is {\it not} a AdS $\Mm\Ll(\mh^2)$-spacetime.
In Section \ref{3cusp} we show some simple examples illustrating such
{\it broken $T$-symmetry}.


\section{Moving along a ray of laminations}\label{der}
In this section we fix a measured geodesic lamination $\lambda$ of
$\mh^2$, invariant by a group $\Gamma<\ISO^+(\mh^2)$ that acts freely
and properly on $\mh^2$ (possibly $\Gamma=\{Id\}$).  Let us put
$F=\mh^2/\Gamma$.

The ray of ($\Gamma$-invariant) measured laminations determined by
$\lambda$ is given by $t\lambda = (\Ll,t\mu)$, $t\geq 0$.  So we have
$1$-parameter families of $\Mm\Ll(\mh^2)$-spacetimes
$\hat\Uu^k_{t\lambda}$, of constant curvature $\kappa\in\{0,1,-1\}$,
diffeomorphic to $F\times\mr_+$, having as universal covering
$\Uu^k_{t\lambda}$. We have also a family of hyperbolic $3$-manifolds
$M_{t\lambda}$, obtained via WR.  $\hat\Pp_{t\lambda}$ is contained in
$\hat\Uu^{-1}_{t\lambda}$ and is the image of the canonical rescaling
of $\hat\Uu^0_{t\lambda}$. Its universal covering is
$\Pp_{t\lambda}\subset \Uu^{-1}_{t\lambda}$.

First, we want to (give a sense and) study the ``derivatives'' at $t=0$
of the spacetimes $\hat\Uu^k_{t\lambda}$, of their holonomies
and ``spectra'' (see below).

\subsection{Derivatives at $t=0$}\label{derivate}

\paragraph{Derivative of spacetimes}
We set $\hat\Uu^k_{t\lambda}/t$ to be the spacetime of constant
curvature $t^2\kappa$ obtained by rescaling the Lorentzian metric of
$\hat\Uu^k_{t\lambda}$ by the constant factor $1/t^2$. We want to
study the limit when $t\to 0$. For the present discussion it is
important to remind that all these spacetimes are well defined only up
a Teichm\"uller-like equivalence relation. So we have to give a bit
of precision on this point. Fix a base copy of $F\times \mr_+$
and let  
\[
\varphi:F\times\mr_+\rightarrow\hat\Uu_\lambda^0
\]
be a 
{\it marked spacetime} representing the equivalence classes of 
$\hat\Uu_\lambda^0$. Denote by $k_0$ the flat Lorentzian
metric lifted on $F\times \mr_+$ via $\varphi$.
A developing map with respect to such a metric is a diffeomorphism
\[
   D:\tilde F\times\mr_+\rightarrow\Uu_\lambda^0\subset\mx_{0} \ .
\] 
Up to translation, we can suppose $0\in\Uu^0_\lambda$.
Notice that, for every $s>0$, the map
\[
   g_s:\Uu_\lambda^0\ni z\mapsto sz\in\mx_{0}
\]
is a diffeomorphism onto $\Uu_{s\lambda}^0$. Moreover, it is
$\Gamma$-equivariant, where, $\Gamma$ is supposed to act on
$\Uu_\lambda^0$ (resp. $\Uu_{s\lambda}^0$) via $f_\lambda$
(resp. $f_{s\lambda}$) as established in Section \ref{flat}.  
Thus $g_s$ induces to the quotient a
diffeomorphism
\[
   \hat g_s:\hat \Uu_\lambda^0\rightarrow\hat \Uu_{s\lambda}^0
\]
such that the pull-back of the metric is simply obtained by
multiplying the metric on $\hat\Uu_\lambda^0$ by a factor $s^2$.\par
Thus the metric $k_s=s^2k$ makes $F\times\mr_+$ isometric to
$\hat\Uu^0_{s\lambda}$.
We want to prove now a similar result for $\kappa=\pm 1$.

The cosmological time of $(F\times\mr_+,k_s)$ is
$\tau_s=s\tau$, where $\tau$ is the cosmological time of
$(F\times\mr_+,k_0)$. It follows that the gradient with respect to $k_s$
of $\tau_s$ does not depend on $s$ and we denote by $X$ this field.
Now suppose $\kappa=-1$ and denote by $h_s$ the metric obtained by
rescaling $k_s$ around $X$ with rescaling functions
\[
   \begin{array}{ll}
   \alpha=\frac{1}{(1+\tau_s^2)^2} \ ,\ \  & \beta=\frac{1}{1+\tau_s^2} \ .
   \end{array}
\]  
We know that $(F\times\mr_+,h_s)$ is isometric to $\Pp_{s\lambda}$.
Moreover the metric $h_s/s^2$ is obtained by a rescaling of the metric $k_0$
along $X$ by rescaling functions
\[
 \begin{array}{ll}
   \alpha=\frac{1}{(1+s^2\tau^2)^2} \ , \ \ & \beta=\frac{1}{1+s^2\tau^2} \ .
   \end{array}
\]  
Thus we obtain $\lim h_s/s^2=k$.\par 
Finally suppose $k=1$, then we
can define a metric $h'_s$ on the subset $\Omega_s$ of $F\times\mr_+$
of points $\{x|\tau_s(x)<1\}=\{x|\tau<1/s\}$ such that $(\Omega_s,
h'_s)=\hat\Uu^1_\lambda$. In fact we can set $h'_s$ to be the metric
obtained by rescaling $k_s$ by rescaling functions
\[
   \begin{array}{ll}
   \alpha=\frac{1}{(1-\tau_s^2)^2}, & \beta=\frac{1}{1-\tau_s^2}.
   \end{array}
\]
Choose a continuous family of embeddings $u_s:F\times\mr_+\rightarrow
F\times\mr_+$ such that
\begin{enumerate}
\item
   $u_s(F\times\mr_+)=\Omega_s$;
\item
 $  u_s(x)=x \textrm{ if } \tau(x)<\frac{1}{2s}$.
\end{enumerate}
Then the family of metrics  $h_s=u_s^*(h'_s)$ works.
\medskip

We can summarize the so obtained results as follows:

\begin{prop}\label{lim:space}
For every $\kappa=0,\pm 1$, 

\[ \lim_{t\rightarrow 0}\ \frac{1}{t} \Uu^\kappa_{t\lambda}=  \Uu^0_\lambda
\ .\]
\end{prop}
\cvd

\noindent 
For $\kappa = 0$ we have indeed the strongest fact that
for every $t>0$
\[\frac{1}{t}\Uu^0_{t\lambda} =  \Uu^0_\lambda \ . \]
Note that this convergence is in fact like a convergence
of pointed-spaces; for example, the convergence of 
spacetimes $\frac{1}{t}\Uu^{-1}_{t\lambda}$ only concerns the
past side of them, while the future sides simply disappear.

\paragraph{Derivatives of representations}
We have seen that for any $\kappa \in\{0,1,-1\}$ the set of holonomies
of $\hat\Uu_{t\lambda}^\kappa$ gives rise to continuous families of
representations
\[
    \rho_t^\kappa:\Gamma\rightarrow\ISO(\mx_{\kappa})
\]
We want to compute the derivative of such families at $t=0$.  The
following lemma contains the formula we need.  In fact this lemma is
proved in~\cite{McM}, and we limit ourselves to a sketch
of proof.
\begin{lem}\label{var:der:lem}
Let $\lambda$ be a {\rm complex-valued} measured geodesic lamination 
on $\mh^2$, and denote by
$E_\lambda$ the Epstein-Marden {\rm bending-quake} cocycle.
Fix two points $x,y\in\mh^2$ then the function
\[
    u_\lambda:\mc\ni z\mapsto E_{z\lambda}(x,y)\in\PSL{2}{C}
\]
is holomorphic.\par Moreover, if $\lambda_n\rightarrow\lambda$ on a
neighbourhood of $[x,y]$, then $u_{\lambda_n}\rightarrow u_\lambda$ \\
in $\Oo(\mc;\PSL{2}{C})$.
\end{lem}
\Dim The statement is obvious when $\lambda$ is a finite lamination.
On the other hand, for every $\lambda$ there exists a sequence of
standard approximations $\lambda_n$. Now it is not hard to see that
$u_{\lambda_n}$ converges to $u_\lambda$ in the compact-open topology
of $\mathrm C^0(\mc;\PSL{2}{C})$. Since uniform limit of holomorphic
functions is holomorphic the first statement is achieved. The same
argument proves the last part of the lemma.  \cvd 
\smallskip

By using the lemma we can easily compute the derivative of $u_\lambda$
at $0$.
\par Notice that $\sG\lG(2,\mc)$ is the complexified of
$\sG\lG(2,\mr)$ that is
\[
    \sG\lG(2,\mc)=\sG\lG(2,\mr)\oplus i\sG\lG(2,\mr) \ .
\]
Now if $l$ is an oriented geodesic denote by $X(l)\in\sG\lG(2,\mr)$
the unitary generator of positive translations along $l$. It is not
hard to show that $iX(l)/2$ is the standard generator of positive
rotation around $l$.  Thus if $\lambda$ is a finite lamination and
$l_1,\ldots,l_n$ are the geodesics between $x$ and $y$ with respective
weights $a_1,\ldots,a_n\in\mc$ we have that
\[
    \frac{\d E_{z\lambda}(x,y)}{\,z}|_0=\frac{1}{2}\sum_{i=1}^{n} a_iX_i \ .
\]
The following statement is a corollary of this formula and
Lemma~\ref{var:der:lem}.
\begin{prop}\label{var:der:prop}
If $\lambda=(\Ll,\mu)$ is a complex-valued measured geodesic
lamination and $x,y$ are fixed points in $\mh^2$ then
\begin{equation}\label{var:der:eq1}
    \frac{\d E_{z\lambda}(x,y)}{\,z}|_0=\frac{1}{2}\int_{[x,y]} X(t)\d\mu(t)
\end{equation}
where $X(t)$ is so defined:
\[
\left\{\begin{array}{ll}
        X(t)=X(l) & \textrm{if } t\in\Ll\textrm{ and }l\textrm{ is the leaf
        through } t\\
        X(t)=0    & \textrm{otherwise \ .}
       \end{array}\right.
\]
\end{prop}
\cvd 
\smallskip

Now we can compute the derivative of $\rho_t^\kappa$ at $0$.  There
exists a canonical linear isometry between $\sG\lG(2,\mr)$, endowed
with its killing form, and the Minkowsky space
$(\mr^3,\E{\cdot}{\cdot})$. In particular, if $l$ is an oriented
geodesic with end-points $x^-,x^+$, then that identification takes
$X(l)$ to the unit spacelike vector $v\in\mx_0$ orthogonal to $l$ such
that $x^-,x^+,v$ form a positive basis of $\mr^3$.
\begin{cor}\label{var:der:DS:cor}
The derivative of $\rho^{1}_{t\lambda}$ at $0$ is an imaginary cocycle
in
\[
\coom1_{\Ad}(\Gamma,\sG\lG(2,\mc))=\coom1_{\Ad}(\Gamma,\sG\lG(2,\mr))\oplus
i\coom1_{Ad}(\Gamma,\sG\lG(2,\mr)) \ .
\]
Moreover, up to the identification of $\sG\lG(2,\mr)$ with $\mr^3$ we
have that
\[
    \dot \rho(0)=\frac{i}{2}\tau_\lambda
\]
where $\tau_\lambda\in\coom1(\Gamma,\mr^{3})$ is the translation part of
$f_\lambda$.
\end{cor}
\cvd 
In the same way we have the following statement
\begin{cor}\label{var:der:ADS:cor}
The derivative of $\rho^{-1}_{t\lambda}$ at $t=0$ is a pair of cocycles
$(\tau-,\tau_+)\in\coom1(\Gamma,\sG\lG(2,\mr))
\oplus\coom1(\Gamma,\sG\lG(2,\mr))$. In particular, if $\tau_\lambda$ 
is the translation part of $f_\lambda$, then
\[
\begin{array}{l}
\tau_-=-\frac{1}{2}\tau_\lambda \ ,\\
\tau_+=\frac{1}{2}\tau_\lambda \ .
\end{array}
\]
\end{cor}
\cvd

\paragraph{Derivatives of the spectra}
Let us denote by $\Cc$ the set of conjugacy classes
of elements of the group $\Gamma$.
For every $\kappa = 0,\pm 1$,
we want to associate to certain elements $[\gamma]\in \Cc$
two numerical ``characters'' 
$\ell_\lambda^\kappa(\langle\gamma\rangle)$ and 
$\Mm_\lambda^\kappa(\gamma)$.

\par First consider $\kappa=0$.  
If $[\gamma]\in \Cc$ is a class of isometries of $\mh^2$ of hyperbolic
type (we simply say that $[\gamma]$ is a hyperbolic element),
then define $\ell_\lambda^0(\langle[\gamma]\rangle)$ to be
the translation length of $\gamma$. $\Mm_\lambda^0([\gamma])$ was
introduced by Margulis in~\cite{Marg}. Denote by
$\tau\in\Z^1(\Gamma,\mr^3)$ the translation part of
$f_\lambda$ (obtained by fixing a base point $x_0\in\mh^2$). 
Denote by $X\in\sG\lG(2,\mr)$ the
unit positive generator of the hyperbolic group containing $\gamma$.
Let $v\in\mr^3$ be, as above, the corresponding point in the Minkowsky
space. Then we have
\[
      \Mm_\lambda^0(\langle[\gamma]\rangle)=\E{v}{\tau(\gamma)} \ .
\] 
It is not hard to see that $\Mm^0_\lambda$ is well defined.\\
\smallskip

Suppose now $k=1$. Take $[\gamma]$ such that $\rho^1_\lambda(\gamma)$ is
 hyperbolic.  In this case $\ell^1_\lambda(\langle[\gamma]\rangle)$ is
 the length of the simple closed geodesic $c$ in
 $\mh^3/\rho_\lambda^1(\gamma)$. On the other hand
 $\Mm_\lambda^1([\gamma])\in[-\pi,\pi]$ is the angle formed by a tangent
 vector $v$ orthogonal to $c$ at a point $x\in
 c\subset\mh^3/\rho_\lambda^1(\gamma)$ with the vector obtained by the
 parallel transport of $v$ along $c$.  By making computation we have
 that
\[
  \tr(\rho^1_\lambda(\gamma))=
2\ch(\frac{\ell^1_\lambda([\gamma])}{2}+i\frac{\Mm^1_\lambda([\gamma])}{2}) \ .
\]
In particular it follows that $\rho^1_\lambda(\gamma)$ is conjugated
to an element of $PSL(2,\R)$ if and only if
$\Mm^1_\lambda([\gamma])=0$.\\
   
Finally suppose $k=-1$. It is not hard to see that if $\gamma$ is a
hyperbolic transformation, then $\rho_\lambda(\gamma)$ is a pair of
hyperbolic transformations too (in fact by choosing the base point on
the axis of $\gamma$ we have that the axis of $\beta_-(x_0,\gamma
x_0)$ intersects the axis of $\gamma$).  Thus if we put
$\rho_\lambda(\gamma)=(\rho_-(\gamma),\rho_+(\gamma))$ there are
exactly two spacelike lines $l_-,l_+$ fixed by $\rho_\lambda(\gamma)$.
Namely $l_+$ has endpoints
\[
\begin{array}{ll}
  p_-=(x^-_L,x^-_R)\ , & p_+=(x^+_L,x^+_R)
\end{array}
\]
and $l_-$ has endpoints
\[ 
\begin{array}{ll}
  q_-=(x^+_L,x^-_R)\ , & q_+=(x^-_L,x^+_R)
\end{array}
\]
where $x^\pm_L$ (resp. $x^\pm_R$) are the fixed points of
$\rho_-(\gamma)$ (resp. $\rho_+(\gamma)$).  Orient $l_+$ (resp. $l_-$)
from $p_-$ towards $p_+$ (resp. from $q_-$ towards $q_+$).  If $m,n$
are the translation lengths of $\rho_-(\gamma)$ and $\rho_+(\gamma)$
then it is not hard to see that $\rho_\lambda(\gamma)$ acts on $l_+$
by a positive translation of length equal to $\frac{m+n}{2}$ and on
$l_-$ by a translation of a length equal to $\frac{n-m}{2}$.  Thus let
us define
\[
\begin{array}{ll}
\ell^{-1}_\lambda([\gamma])=\frac{m+n}{2}\\
\Mm^{-1}_\lambda([\gamma])=\frac{n-m}{2} \ .
\end{array}
\]
\begin{prop}
If $\gamma$ is a hyperbolic element of $\Gamma$ then there exists
$t<1$ sufficiently small such that $\rho_{s\lambda}^1(\gamma)$ is
hyperbolic for $s<t$.  Moreover the following formulas hold
\[
\begin{array}{l}
   \frac{\d \ell_{t\lambda}^\kappa([\gamma])}{\,\d t}|_0=0\\
   \frac{\d \Mm_{t\lambda}^\kappa([\gamma])}{\,\d t}|_0=
\Mm_\lambda^0([\gamma]) \ .
\end{array}
\]
\end{prop}
\Dim For $k=0$ the statement is trivial.
\par Suppose $k=1$. Denote by
$B_t$ the cocycle associated to the lamination $\lambda_t$
\[
    \tr(B_t(x_0,\gamma x_0)\gamma)=
2\ch(\frac{\ell^1_t([\gamma])+i\Mm^1_t([\gamma])}{2} \ .)
\]
By deriving at $0$ we obtain
\[
   \frac{1}{2}\tr(i X(\gamma)\gamma)=\sh(\ell([\gamma])/2)(\dot
   \ell^1([\gamma])|_0 + i \dot\Mm^1([\gamma])|_0)
\]
where $X(\gamma)$ is the element of $\sG\lG(2,\mr)$ corresponding to
$\tau(\gamma)\in\mr^3$ (where $\tau$ is the is the translation part of
$f_\lambda$ ).  Now if $Y\in\sG\lG(2,\mr)$ is the unit generator of
the hyperbolic group containing $\gamma$ we have
\[
     \gamma=\ch(\ell([\gamma])/2) I + \sh(\ell([\gamma])/2) Y \ .
\]
Thus we obtain
\[
  \dot \ell^1([\gamma])|_0 + i \dot\Mm^1([\gamma])|_0=i\Mm^0([\gamma]) \ . 
\]
An analogous computation shows the same result in the case $\kappa=-1$.
\cvd

\subsection {More cocompact case}\label{cocompact}
\paragraph{About $(\Uu^{-1}_{t\lambda})^*$}
Let us assume now that $\Gamma$ is cocompact,
so that $F=\mh^2/\Gamma$ is compact of genus $g\geq 2$.
With the notations of Subsection \ref{Tsym}, we have also
the family of AdS spacetimes 
$$(\Uu^{-1}_{t\lambda})^* = \Uu^{-1}_{\lambda_t^*} \ .$$
We want to determine the derivative at $t=0$ of this
family.
Remind (see Section \ref{flat}) that in such a case the set
of $\Gamma$-invariant measured lamination has a natural
$\R$-linear structure, so it makes sense to consider $-\lambda$.
We have (the meaning of the notations is as above)
\begin{prop}\label{der:inv:prop}
\[ \lim_{t \rightarrow 0}\  \frac{1}{t}\Uu^{-1}_{\lambda^*_t}=  
\Uu^0_{-\lambda}
\ .\]
\end{prop}
\Dim Let $(\rho'_t,\rho''_t)$ be the holonomy of
$\ \Uu^{-1}_{t\lambda}$.  Denote by $F^*_t$ the quotient of the past
boundary of $K_{t\lambda}$ by $(\rho'_t,\rho''_t)$. Notice that
$\lambda^*_t$ is a measured geodesic lamination on $F^*_t$.  We claim
that $(F^*_t,\lambda^*_t/t)$ converges to $(F,-\lambda)$ in
$\Tt_g\times\Mm\Ll_g$ as $t\rightarrow 0$.  Let us show first how this
claim implies the proposition.  We can choose a family of developing
maps
\[
   D_t:\tilde F\times\mr\rightarrow\Uu^0_{\lambda^*_t/t}\subset\mx_0
\]
such that $D_t$ converges to a developing map $D_0$ of
$\hat\Uu^0_{-\lambda}=\Uu^0_{-\lambda}/f_{-\lambda}$ as $t\rightarrow
0$.  Denote by $k_t$ the flat Lorentzian metric on $\tilde F\times\mr$
corresponding to the developing map $D_t$. We have that $k_t$ converges
to $k_0$ as $t\rightarrow 0$.  Moreover, if $T_t$ denotes the
cosmological time on $F\times\mr$ induced by $D_t$, then $T_t$
converges to $T_0$ in $\mathrm C^1(F\times\mr)$ as $t\rightarrow 0$.
Now, as in the proof of Proposition~\ref{lim:space}, we have that
$\Pp_{\lambda^*_t}$ is obtained by a Wick Rotation directed by the gradient
of $T_t$ with rescaling functions
\[
  \begin{array}{ll}
   \alpha=\frac{t^2}{(1+(tT_t)^2)^2} \ , \ \ & \beta
   =\frac{t^2}{1+(tT_t)^2} \ .
  \end{array}
\]
Thus we easily obtain the statement.\\

In order to prove the claim, first we will show that the set
$\{(F^*_t,\lambda^*_t/t)|t\in [0,1]\}$ is pre-compact in
$\Tt_g\times\Mm\Ll_g$, and then we will see that the only possible
limit of a sequence $(F^*_{t_n},\lambda^*_{t_n}/t_n)$ is
$(F,-\lambda)$.\\

We have seen that $F'_t=\mh^2/\rho'_t(\Gamma)$
(resp. $F''_t=\mh^2/\rho''_t(\Gamma)$) is obtained by a right
(resp. left) earthquake on $F=\mh^2/\Gamma$ with shearing measured
lamination equal to $t\lambda$.  Thus, if $\lambda'_t$ is the measured
geodesic lamination of $F'_t$ corresponding to $t\lambda$ via the
canonical identification of $\Mm\Ll(F)$ with $\Mm\Ll(F'_t)$, we have
that $F''_t$ is obtained by a left earthquake on $F'_t$ along
$2\lambda'_t \ $.\par On the other hand let $(\lambda^*_t)'$ be the
measured geodesic lamination on $F'_t$ such that the right earthquake
along it sends $F'_t$ on $F''_t$.  Then we know that the quotient
$F^*_t$ of the past boundary of the convex core $K_{t\lambda}$ is
obtained by a right earthquake along $F'_t$ with shearing lamination
$(\lambda^*_t)'$. Moreover, the bending locus $\lambda^*_t$ is the
lamination on $F^*_t$ corresponding to $2(\lambda^*_t)'$.\par Clearly
$F^*_t$ runs in a compact set of $\Tt_g$ as $t$ runs in $[0,1]$.  In
order to prove that the family $\{(F^*_t,\lambda^*_t/t)|t\in [0,1]\}$
is pre-compact we will use some classical facts about $\Tt_g$. For the
sake of clearness we will remind them, referring to
~\cite{Pen,Ot} for details.\par Denote by $\Cc$ the set of
conjugacy classes of $\Gamma$.
For $\lambda\in\Mm\Ll(S)$ we denote by $\iota_\gamma(\lambda)$ the
total mass of the closed geodesic curve corresponding to $[\gamma]$
with respect to the transverse measure given by $\lambda$.  The
following facts are well-known.
\begin{enumerate}
\item
Two geodesic laminations $\lambda$ on $S$ and $\lambda'$ on $S'$ are
identified by the canonical identification
$\Mm\Ll(S)\rightarrow\Mm\Ll(S')$ if and only if
$\iota_\gamma(\lambda)=\iota_\gamma(\lambda')$ for every
$[\gamma]\in\Cc$.
\item
A sequence $(F_n,\lambda_n)$ converges to $(F_\infty,\lambda_\infty)$
in $\Tt_g\times\Mm\Ll_g$ if and only if $F_n\rightarrow F_\infty$ and
$\iota_\gamma(\lambda_n)\rightarrow\iota_\gamma(\lambda_\infty)$ for
every $[\gamma]\in\Cc$.
\item
A subset $\{(F_i,\lambda_i)\}_{i\in I}$ of $\Tt_g\times\Mm\Ll_g$ is
pre-compact if and only if the base points $\{F_i\}$ runs in a compact
set of $\Tt_g$ and for every $[\gamma]\in\Cc$ there exists a constant
$C>0$ such that
\[
     \iota_\gamma(\lambda_i)<C\qquad\textrm{ for every }i\in I.
\]
\end{enumerate}
Clearly we have $F^*_t\rightarrow F$ as $t\rightarrow 0$.  Thus in
order to show that $(F_t,\lambda^*_t)$ is pre-compact it is sufficient
to find for every $[\gamma]\in\Cc$ a constant $C>0$ such that
\[
    \iota_\gamma(\lambda^*_t)<Ct
\]
for every $t\in [0,1]$.\par  
The following lemma gives the estimate we need.
\begin{lem}\label{der:inv:lem}
For every compact set $K\subset\mh^2$ there exists a constant $C>0$
which satisfies the following statement.\par If $\lambda$ is a
measured geodesic lamination on $\mh^2$ and $\beta$ is the right
cocycle associated to $\lambda$ then
\[
   ||\beta(x,y)-Id+\frac{1}{2}
\int_{[x,y]}X_\lambda(u)\d\lambda||\leq e^{M\lambda(x,y)}-1-M\lambda(x,y)
\]
where $X_\lambda(u)$ is defined as in~(\ref{var:der:eq1}) and  $x,y\in K$.
\end{lem}
\Dim It is sufficient to prove the lemma when $\lambda$ is simplicial.
In this case denote by $l_1,\ldots,l_N$ the geodesics meeting the
segment $[x,y]$ with respective weights $a_1,\ldots,a_N$.  If
$X_i\in\sG\lG(2,\mr)$ is the unitary infinitesimal generator of the
positive translation along $l_i$ we have
\[
   \beta(x,y)= \exp(-a_1X_1/2)\circ\exp(-a_2X_2/2)\circ\cdots\circ\exp(-a_N
X_N/2)\ .
\]
Thus $\beta(x,y)$ is a real analytic function of $a_1,\ldots,a_n$.
If we write
\[
   \beta(x,y)=\sum_n A_n(a_1,\ldots,a_n)
\]
where $A_n$ is a matrix-valued homogenous polynomial in $x_1,\ldots,x_n$ of
degree $n$, then it is not difficult to see that
\[
     ||A_n||\leq (\sum_{i=1}^N a_i||X_i||)^n/n!\ .
\]
We have that
\[
   \beta(x,y)-Id+\frac{1}{2}
\int_{[x,y]}X_\lambda(u)\d\lambda=\sum_{i\geq 2} A_n(a_1,\ldots, a_N)
\]
Since the axes of trasformations generated by $X_i$
cut $K$, there exists a constant $M>0$ (depending only on $K$) such
that $||X_i||<M$.  Thus
\[
   ||\beta(x,y)-Id+\frac{1}{2}
\int_{[x,y]}X_\lambda(u)\d\lambda||\leq e^{M\sum a_i}-1-M\sum a_i\ .
\]
\cvd
Let us go back to the proof of Proposition~\ref{der:inv:prop} .
Since
$\iota_\gamma(\lambda^*_t)=\iota_\gamma((\lambda^*_t)')$, we can
replace $\lambda^*_t$ with $(\lambda^*_t)'$.  Now let us put
$\gamma_t=\rho'_t(\gamma)$. We know that $\gamma_t$ is a
differentiable path in $\PSL{2}{R}$ such that $\gamma_0=\gamma$ and
\[
   \dot\gamma(0)=-\frac{1}{2}\int_{[x,\gamma(x)]} X(u)\d\lambda(u)
\]
where $X(u)$ is defined as in~(\ref{var:der:eq1}).  On the other hand,
if $\beta_t$ is the right cocycle associated to the measured geodesic
lamination $(\lambda^*)'_t$ we have
\[
  \beta_t(x,\gamma_tx)\gamma_t=\beta_{t\lambda}(x,\gamma x)\gamma
\]  
where $\beta_{t\lambda}$ is the left cocycle associated to $t\lambda$.
Thus $\beta_t(x,\gamma_tx)$ is a differentiable path and
\begin{equation}\label{der:rev:eq1}
  \lim_{t\rightarrow 0}
  \frac{\beta_t(x,\gamma_tx)-Id}{t}=\int_{[x,\gamma(x)]}X(u)\d\lambda\ .
\end{equation}
By using Lemma~\ref{der:inv:lem} it is not hard to show that there exists a
constant $C>0$ depending only on $\gamma$ such that
\[
    ||\beta_t(x,\gamma_tx)-Id||> \frac{1}{2}||\int_{[x,\gamma_tx]}
      X(u)\d\lambda^*_t|| - C \iota_\gamma((\lambda^*_t)')2 .
\]
On the other hand, there exists a constant $L>0$ such tah
\[
||\int_{[x,\gamma_tx]}
      X(u)\d\lambda^*_t||\geq L |\int_{[x,\gamma_tx]} X(u)\d\lambda^*_t|  
\]
where $|\cdot|$ denotes the Lorentzian norm of $\sG\lG(2,\mr)$.  Since
$X(u)$ are generators of hyperbolic transformations with disjoint axes
pointing in the same direction, the following inequality holds
\[
   \eta(X(u),X(v))2\geq\eta(X(u),X(u))\eta(X(v),X(v))>0
\] 
and implies
\[
  ||\int_{[x,\gamma_tx]}
      X(u)\d(\lambda^*_t)'||\geq L\iota_\gamma((\lambda^*_t)') \ .
\]
From this inequality we easily obtain that
\[
     \left(L-C\lambda^*_t(x,\gamma_tx)\right)\iota_\gamma((\lambda^*_t)')
<||\beta_t(x,\gamma_tx)-Id|| \ .
\]
Thus $\lambda^*_t(x,\gamma_tx)$ converges to $0$. Moreveover, by
dividing by $t$ the last inequality, we obtain that
$\lambda^*_t(x,\gamma_tx)/t$ is bounded.  In particular we have proved
that $\{(F_t,\lambda^*_t)\}$ is pre-compact in
$\Tt_g\times\Mm\Ll_g$.\\

Now, let us set $\mu_t=\lambda^*_t/t$ and $\mu'_t=(\lambda^*)'_t/t$.
We have to show that if $\mu_{t_n}\rightarrow \mu_\infty$ then
$\mu_\infty=-\lambda$ in $\Mm\Ll(F)$.\par Notice that $\mu'_{t_n}$ is
convergent and its limit is $\mu_\infty$.
Applying lemma~\ref{der:inv:lem} we get 
\[
   \lim_{t\rightarrow 0}\frac{\beta_t(x,\gamma_tx)-Id}{t}=
   -\int_{[x,\gamma(x)]}X_{\mu_\infty}(t)\d\mu_\infty
\]
By equation~(\ref{der:rev:eq1}) this limit is equal
to $\int_{[x,\gamma(x)]}X\d\lambda$ and this shows that
$\mu_\infty=-\lambda$.  \cvd

Till now we have derived infinitesimal information at $t=0$.
For what concerns the behaviour along a ray for {\it big} $t$,
let us make a qualitative remark.

We have noticed that, for every $t>0$, $\frac{1}{t}\Uu^0_{t\lambda}=
\Uu^0_\lambda$. Moreover, the flat spacetimes $\ \Uu^0_{t\lambda}\ $ are
nice convex domains in $\mx_0$ which vary continously and tamely with
$t$.  So, in the flat case, {\it apparently nothing of qualitatively
new does happen when $t>0$ varies.}  Similarly, this holds also for
$\Pp_{t\lambda}\subset \Uu^{-1}_{t\lambda}$.  On the other hand,
radical qualitative changes do occur for $M_{t\lambda}$ (and
$\Uu^{1}_{t\lambda}$) when $t$ varies. As $\lambda$ is
$\Gamma$-invariant for the cocompact group $\Gamma$, when $t$ is small
enough $M_t$ is in fact the universal covering of an end $\hat M_t$ of
a {\it quasi-Fuchsian} complete hyperbolic $3$-manifold, $Y_t$ say.
In particular, the developing map is an embedding. When $t$ grows up,
we find a first value $t_0$ such that we are no longer in the
quasi-Fuchsian region, and for bigger $t$ the developing map becomes
more and more ``wild''.  We believe that this different behaviour
along a ray is conceptually important. It substantially supports the
conception of 3D gravity as an unitary body: looking only at
the flat Lorentzian sector, for example, significant critical
phenomena should be lost; on the other hand, one could consider the WR
as a kind of ``normalization'' of the hyperbolic developing map. We
believe that this should be useful in studying ending invariants of
hyperbolic $3$-manifolds.
\smallskip

Anyway, we give here a first simple application of these qualitative
considerations.

Assume that we are in the quasi-Fuchsian region.
So we have associated to $t\lambda$ three ordered pairs
of elements of the Teichm\"uller space $\Tt_g$. These are:
\smallskip

- the ``Bers parameter'' $(B_t^+,B_t^-)$
given by the conformal structure underlying the projective
asymptotic structures of the two ends of $Y_t$; 
\smallskip

- the hyperbolic structures $(C_t^+, C_t^-)$ of the boundary components
of the hyperbolic convex core of $Y_t$;
\smallskip

- the hyperbolic structures $(K_t^+,K_t^-)$ of the future and past
boundary components of the AdS convex core of $\Uu^{-1}_{t\lambda}$.
\smallskip

\noindent
It is natural to inquire about the relationship between
these pairs.
\smallskip
\par By construction we have that $K_t^+$ is isometric to $C_t^+$.  On
the other hand by a Sullivan's Theorem (see \cite{Ep-M}) we have that
the Teichm\"uller distance of $B_t^\pm$ from $C_t^\pm$ is less than
$2$.  Now it is natural to ask whether $C_t^-$ is isometric to
$K_t^-$.  Actually it is not hard to show that those spaces generally
are not isometric.  In fact let us fix a lamination $\lambda$ and let
$t_0>0$ be the first time such that the representation
$\rho^1_{t_0\lambda}$ is not quasi-Fuchsian.  By Bers Theorem
(c.f. again ~\cite{Ep-M}) we know that the family $\{(B^+_t,
B^-_t)\in\Tt_g\times\overline\Tt_g\}$ is not compact. Since $B^+_t$
converges to a conformal structure as $t$ goes to $t_0$ we have that
$\{B^-_t\}_{t\leq t_0}$ is a divergent family in $\overline\Tt_g$.  By
Sullivan's Theorem we have that $C_t^-$ is divergent too. On the other
hand $K_t^-$ is a pre-compact family. It follows that in general
$C_t^-$ is not isometric to $K_t^-$.  \cvd

\section{Examples}\label{3cusp}
In this section we show some examples that illustrate the results
obtained in this paper. Most of the following considerations apply to
any non-compact hyperbolic surface $F=\mh^2/\Gamma$ of finite area,
endowed with an ideal tringulation. However, we are going to focus the
simplest example, that is $F= \mh^2/\Gamma$ the {\it three-punctured
sphere}. Some features of this example have been sketched also
in \cite{Thu2, M}.  

It is well known that $F$ is {\it rigid}, that is the corresponding
Teichm\"uller space is reduced to one point. Hence, $F$ has no
measured geodesic laminations with compact support.

$F$ can be obtained by gluing two geodesic ideal triangles along their
edges as follows.  In any ideal triangle there exists a unique point
(the "barycenter") that is equidistant from the edges. In any edge
there exists a unique point that realizes the distance of the edge
from the barycenter. Such a point is called the mid-point of the edge.
Now the isometric gluing is fixed by requiring that mid-points of
glued edges matches (and that the so obtained surface is topologically
a three-punctured sphere - by a different pattern of identifications
we can obtain a $1$-punctured torus).  It is easy to see that the
resulting hyperbolic structure is complete, with a cusp for any
puncture, and equipped by construction with an {\it ideal
triangulation}.

The three edges of this triangulation form a geodesic lamination
$\Ll_F$ of $F$. A transverse measure $\mu_F = \mu_F(a_1,a_2,a_3)$ on
such a lamination consists in giving each edge a positive weight
$a_i$.  The ideal triangles in $F$ lift to a tesselation of the
universal cover $\mh^2$ by ideal triangles.  The $1$-skeleton $\Ll$ of
such a tesselation is the pull-back of $\Ll_F$; a measure $\mu_F$
lifts to a $\Gamma$-invariant measure $\mu$ on $\Ll$. So, we will
consider the $\Gamma$-invariant measured laminations $\lambda =
(\Ll,\mu)$ on $\mh^2$ that arise in this way.

 \paragraph {Blind flat Lorentzian holonomy}
Varying the weights $a_i$, we get a $3$-parameters family of flat
spacetimes $\Uu^0_{\lambda}= \Uu^0_{\lambda(a_1,a_2,a_3)}$,
with associated quotient spacetimes  $\hat \Uu^0_{\lambda} =
\Uu^0_{\lambda}/f_\lambda(\Gamma)$, where $f_\lambda(\Gamma)$ denotes
the flat Lorentzian holonomy (Section \ref{flat}).

The spacetimes $\Uu^0_{\lambda}$ have {\it homeomorphic} initial
singularities $\Sigma_\lambda$.  A topological model for them is given
by the {\it simplicial} tree (with $3$-valent vertices) which forms
the $1$-skeleton of the cell decomposition of $\mh^2$ dual to the
above tesselation by ideal triangles. The length-space metric of each
$\Sigma_\lambda$ is determined by the fact that each edge of the tree
is a geodesic arc of length equal to the weight of its dual edge of
the triangular tesselation. In fact, every $\Sigma_\lambda$ is
realized as a spacelike tree embedded into the frontier of
$\Uu^0_{\lambda}$ in the Minkowski space, and $f_\lambda(\Gamma)$ acts
on it by isometry.

The behaviour of the {\it asymptotic states} of the cosmological time
of each $\Uu^0_{\lambda}$, is formally the same as in the cocompact
case (see the end of Section \ref{flat}). In particular, when $a\to
0$, then action of $f_\lambda(\Gamma)$ on the level surface
$\Uu_\lambda(a)$ converges to action on the initial singularity
$\Sigma_\lambda$.  The marked length spectrum of $\hat \Uu_\lambda(a)$
(which coincides with the minimal deplacement marked spectrum of the
action of $f_\lambda(\Gamma)$ on $\Uu_\lambda(a)$), converges to the
minimal deplacement marked spectrum of the isometric action on the
initial singularity.  If $\gamma_i$ $i=1,2,3$ are (the conjugacy
classes of) the parabolic elements of $\Gamma$ corresponding to the
three cusps of $F$, the last spectrum takes values $\gamma_i \to
a_i+a_{i+1}$, where we are assuming that the edges of the ideal
triangulation of $F$ with weigths $a_i$ and $a_{i+1}$ enter the
$i$th-cusp, $a_4=a_1$.  By the way, this implies that these spacetimes
are not isometric each other.

However, it follows from \cite{Ba} that:
\medskip

{\it 
(1) The flat Lorentzian holonomies $f_\lambda(\Gamma)$ are 
all conjugate (by $ISO_0(\mx_0)$) to their common 
{\rm linear part} $\Gamma$.
\smallskip

(2) There are non trivial classes in $H^1(\Gamma, \mr^3)$ that
are not realized by any flat spacetime having $\Gamma$ as linear
holonomy.}
\medskip

\noindent Hence, the flat Lorentzian holonomy is completely {\it
blind} in this case, and, on the other hand, the non trivial algebraic
$ISO_0(\mx_0)$-extensions of $\Gamma$ are not correlated to the
geometry of any spacetime.  Apparently this is a completely different
behaviour with respect to the cocompact case (see Section \ref{flat}).

\begin{remarks}
{\rm (1) The above facts would indicate that the currently claimed
  equivalence between the classical formulation of 3D gravity in terms
  of Einstein action on metrics, and the formulation via Chern-Simons
  actions on connections (see \cite{W}), should be managed instead
  very carefully outside the cocompact $\Gamma$-invariant range (see
  \cite{Mat} for similar considerations about flat spacetimes with
  particles).
\smallskip

(2) Every $\Uu^0_{\lambda}$ can be embedded in a $\Gamma$-invariant
way in the static spacetime $I^+(0)$ as the following construction
shows (this is the geometric meaning of point (1) above) :
\smallskip
Up to conjugating $f_\lambda$ by an ismoetry of $\ISO(\mx_{0})$ we can suppose
$f_\lambda(\gamma)=\gamma$ for every $\gamma\in\Gamma$.
We want to prove that $\Uu_\lambda$ is contained in $\fut(0)$. By
contradiction suppose there exists $x\in\Uu_\lambda$ outside the future of
$0$. Since $\Uu_\lambda$ is future complete $\fut(x)\cap\partial\fut(0)$
is contained in $\Uu_\lambda$. But we know there is no open set of
$\partial\fut(0)$ such that the action of $\Gamma$ on it is free. Thus we
obtain a contradiction.

There is a geometric way to recognize such domains inside $\fut(0)$.
Take a $\Gamma$-invariant set of horocircle $\{B_n\}$ in $\mh^2$
centred to points corresponding to cusps.  We know that every
horocircle $B_n$ is the intersection of $\mh^2$ with an affine null
plane orthogonal to the null-direction corresponding to the center of
$B_n$.  Now it is not difficult to see that the set
\[
    \Omega= \bigcup\fut(P_n)
\]
is a regular domain invariant by $\Gamma$.  This is clear if $B_n$ are
sufficientely small (in that case we have that
$\Omega\cap\mh^2\neq\varnothing$). For the general case denote by
$\hat B_n(a)$ the intersection of $P_n$ with the surface
$a\mh^2$. Then the map
\[
    f_a: a\mh^2\ni x\mapsto x/a\in\mh^2
\]
sends $\hat B_n(a)$ to a horocircle $B_n(a)$ smaller and smaller as
$a$ increases.  It follows that $\Omega\cap a\mh^2\neq\varnothing$ for
$a>>0$.  Since a regular domain is the intersection of the future of
its null-support planes it follows that every $\Uu_\lambda$ can be
obtained in this way.

\smallskip
}
\end{remarks}

\paragraph{Earthquake failure and broken $T$-symmetry}
Here we adopt the notations of Section \ref{AdS}.
\smallskip

{\it Earthquake failure.}  It was already remarked in \cite{Thu2} that
$\lambda$ produces neither left nor right earthquake.  Nevertheless,
the representations $\rho_L-$ and $\rho_+$ are discrete.  Both the
limit set $\Lambda_-$ of $\rho_-$ and $\Lambda_+$ of $\rho_+$) are
Cantor sets such that both quotients of the respective convex hulls
are isometric to a same hyperbolic pant with totally geodesic boundary
$\Pi_\lambda$. If $\gamma_i$ is as above, then we have that
$\rho_-(\gamma_i)$ and $\rho_+(\gamma_i)$ are hyperbolic
transformations corresponding to the holonomy of a boundary component
of $\Pi_\lambda$, with translation length equal to $a_i+a_{i+1}$.
\par Thus $\rho_-$ and $\rho_+$ are conjugated in $PSL(2,\R)$.  The
above remarks imply that there exists a spacelike plane $P$ in
$\mx_{-1}$ containing $\Lambda_-\times\Lambda_+$.  It follows that for
any sequences $(x_n)\in\partial\mx_{-1}$ and $(\delta_n)\in\Gamma$ the
accumulation points of the set $\{\rho(\delta_n)x_n\}$ are contained
in $P$.
\smallskip

{\it Broken $T$-symmetry.} Consider now the bending of $\mh^2$ along
$\lambda$. The key point is to describe the curve $c_\lambda$, that is
the curve at infinity of $\Uu^{-1}_\lambda$.

\begin{figure}
\begin{center}
\input{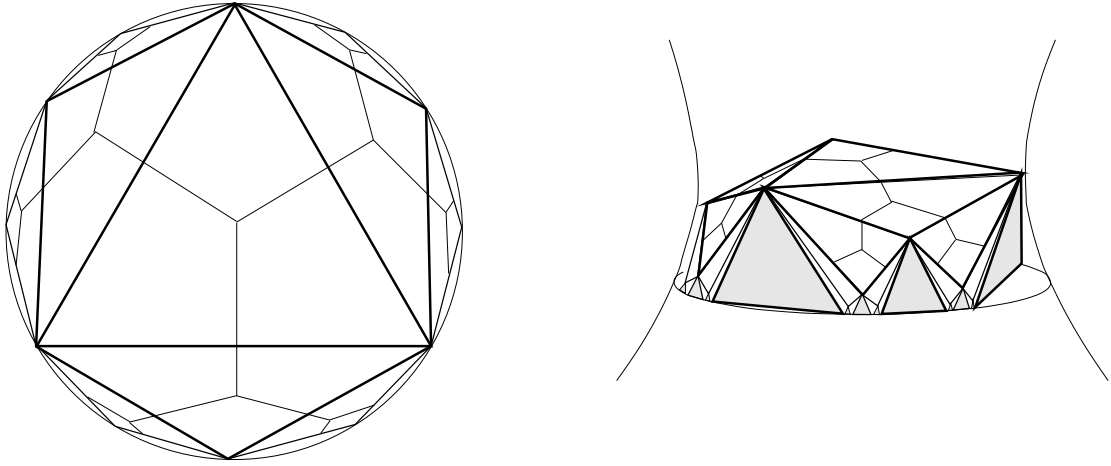}
\caption{{\small On the left the lamination $\Ll$ with its dual spine. On the
    right the bending of $\mh^2$ along $\lambda$ in $\mx_{-1}$. Grey regions
    are null components of the past boundary of $K_\lambda$.}} 
\end{center}
\end{figure}

Take a point $x\in\partial S^1$ that is vertices of a triangle $T$ of
$\lambda$. The point $x$ corresponds to a puncture of $F$ so there is
a parabolic transformation $\zeta\in\Gamma$ conjugated to one
$\gamma_i$ that fixes $x$.  Moreover, we can choose $\zeta$ in such a
way that it is conjugated to a translation $z\mapsto z+a$ with $a>0$
in $PSL(2,\R)$.  It is not hard to see that the point
$\beta_-(x_0,x)x$ is the attractive fixed point of $\rho_-(\zeta)$
whereas $\beta_+(x_0,x)$ is the repulsive fixed point of
$\rho_+(\zeta)$.  \par Choose an edge $l_0$ of $T$ with end-point at
$x$, and let $y_0$ be the end-point of $l_0$ different of $x$.  We
have that $T_n=\varphi_\lambda(\zeta^n(T))$ is a totally geodesic
ideal triangle embedded in the surfaces $\varphi_\lambda(\mh^2)$. If
we take $z\in T$, it is easy to see that a vertices of $T_n$ is the
point $u_n=\beta(x_0,\zeta^n x)(\zeta^n y_0, \zeta^n y_0)$.  Then, a
simple computation shows that
\[    u_n=(\rho_-(\zeta^n),\rho_+(\zeta^n))
(\beta_-(x_0,x)y_0, \beta_+(x_0,x)y_0).\] Thus we have that $c_n$
converges to the point $u_\infty=(x^+_L(\zeta), x^+_R(\zeta))$ where
$x^+_L(\zeta)$ and $x^+_R(\zeta)$ are respectively the attractive
fixed points of $\rho_-(\zeta)$ and $\rho_+(\zeta)$, that is to a
point in $c_\lambda$.  By exchanging $\zeta$ with $\zeta^{-1}$, we see
that the point $v_\infty=(x^-_L(\zeta), x^-_R(\zeta))$ lies in
$c_\lambda$.  By looking at the image of $T$ through $\varphi_\lambda$
we easily see that the point
$u=(x^+_L(\zeta),x^-_R(\zeta)=(\beta_-(x_0,x)x,\beta_+(x_0,x)x)$ lies
in $c_\lambda$. Since $c_\lambda$ is achronal and $u$ and $u_\infty$
are in the same left leaf, it follows that the future directed segment
on the left leaf from $v_\infty$ towards $u$ is contained in
$c_\lambda$.  In the same way we have that the future directed segment
on the right leaf from $u_\infty$ towards $u$ is contained in
$c_\lambda$.  \par Notice that $u_\infty$ and $v_\infty$ are the
vertices of a boundary component of the convex hull $H$ of
$\Lambda_L\times\Lambda_R$ in $P$.  Now take a boundary component $l$
of $H$ oriented in the natural way and let $u_-$ and $u_+$ be its
vertices. Then the left leaf through $u_-$ and the right leaf through
$u_+$ meet each other at a point $u_l$.  The future directed segments
from $u_\pm$ towards $u$ in the respective leaves is contained in
$c_\lambda$. Denote by $V_l$ the union of such segments.  We have that
the union of $V_l$, for $l$ varying in the boundary components of $H$,
is contained in $c_\lambda$. But it is not hard to see that its
closure is a closed path so that $c_\lambda$ coincides with it.

Since the description of the curve $c_\lambda$ is quite simple we can
describe also the past boundary $\partial_- K$ of the convex hull
of $c_\lambda$, i.e. of the AdS convex core of $\Uu^{-1}_\lambda$. 
Notice that $P$ is the unique spacelike support plane
touching $\partial_- K$.  Then for every component $l$ of $H$, there
exists a unique null support plane $P_l$ with dual point at $u_l$. Thus
we have that $\partial_- K$ is the union of $H$ and an infinite number of
null half-planes, each attached to a boundary component of $H$.  
It follows that $\partial_- K$ is {\it not} complete. 

It is not hard to verify that also $(\Uu^{-1}_\lambda)^*$
has canonical cosmological time $T^*$, which passes to
the quotient spacetime $(\hat \Uu{-1}_\lambda)^*$.
A level surface $\hat \Uu^{-1}_\lambda(a)$, $a<\pi/2$, 
of the quotient spacetime
is isometric to the surface obtained by rescaling the metric on the
pant $\Pi_\lambda$ by a factor $\sin a$, and attaching to every
boundary components a infinite Euclidean cylinders.
By using the results of \cite{BBo}, we can say that the {\it future}
singularity of  $\Uu^{-1}_\lambda$ is "censured" by three (static)
{\it BZT black holes}. 
\smallskip

{\it Rescaling and WR.}  The ($\Gamma$-invariant) formulas of the
(inverse of) canonical rescaling of Section \ref{AdS}, apply to
$(\Uu^{-1}_\lambda)^*$ with its cosmological time $T^*$, and produce a
flat regular domain $(\Uu^0_\lambda)^*$ in $\mx_0$.  The above pant
$\Pi_\lambda$ is the convex core $C=C(F')$ of the complete hyperbolic
surface of infinite area $F'= \mh^2/ \Gamma'$, where $\Gamma' = \rho_L
\cong \rho_R$ does not contain any parabolic element. $\Pi_\lambda$
lifts to an ideal convex domain $\widetilde{C}$ in $\mh^2$ (i.e. the
convex hull of the limit set $\Lambda_{\Gamma'} \subset S^1_\infty$).
It is not hard to verify that $\widetilde{C}$ is the image of the
Gauss map of $(\Uu^0_\lambda)^*$.  In fact $(\Uu^0_\lambda)^*$ is
obtained by formally giving each boundary component of $C$ the weight
$+\infty$, hence the same weight is given each geodesic line in the
boundary of $\widetilde{C}$.  Then, we can adapt the constructions of
Section \ref{flat}.  More precisely, let $N$ be the Gauss map of the
static flat regular domain; then each level surface of the
cosmological time of $(\Uu^0_\lambda)^*$ is obtained by taking the
intersection with $N^{-1}(\widetilde{C})$ of the corresponding level
surface of the static spacetime, and glueing to each so obtained
boundary component an infinite Euclidean half-plane.  It is clear that
such a spacetime can be considered as the limit, when $t\to +\infty$,
of $\Uu^0_{t\lambda'}$, where $\lambda'$ is the $\Gamma'$-invariant
lamination having as support the boundary of $\widetilde{C}$ and we
give the weight $1$ each leaf of $\lambda'$.

The WR formulas of Section \ref{hyp} apply to $(\Uu^0_\lambda)^*$ and
one shows that the completion of the resulting hyperbolic $3$-manifold
is a handlebody, quotient of a Schottky group (we refer to \cite{BBo}
for more details about similar examples).


\end{document}